\documentclass[a4paper,10pt]{amsart}
\usepackage[T1]{fontenc}
\usepackage{lmodern}
\usepackage{amssymb}
\usepackage{thmtools}
\usepackage{mathtools}
\usepackage{enumitem}
\usepackage{caption}
\usepackage{tabularx}
\usepackage{tikz}
\usetikzlibrary{arrows.meta,positioning,decorations.pathmorphing}
\usepackage[
	pdfauthor={},
	pdftitle={On the Orthogonal Group of Quadratic Modules},
	pdfsubject={},
	pdfkeywords={},
	pdfproducer={},
	pdfcreator={},
]{hyperref}

\theoremstyle{plain}
\newtheorem{thm}{Theorem}[section]

\newtheorem{prop}[thm]{Proposition}
\newtheorem{cor}[thm]{Corollary}
\newtheorem{lem}[thm]{Lemma}

\theoremstyle{definition}
\newtheorem{defi}[thm]{Definition}
\newtheorem{nota}[thm]{Notation}
\newtheorem{rem}[thm]{Remark}
\newtheorem{constr}[thm]{Construction}
\newtheorem{ex}[thm]{Example}

\newtheoremstyle{introTheorems}
	{3pt}
	{3pt}
	{\itshape}
	{}
	{\bfseries}
	{.}
	{5pt plus 1pt minus 1pt}
	{\thmname{#1}\textnormal{ \thmnote{#3}}}

\theoremstyle{introTheorems}
\newtheorem{introTheorem}{Theorem}

\newtheorem{introLem}{Lemma}
\newtheorem{introConstr}{Construction}
\newtheorem{introSpecial}{Sketch of Theorem}

\newtheoremstyle{introDefinitions}
	{3pt}
	{3pt}
	{\upshape}
	{}
	{\bfseries}
	{.}
	{5pt plus 1pt minus 1pt}
	{\thmname{#1}\textnormal{ \thmnote{#3}}}

\theoremstyle{introDefinitions}
\newtheorem{introDefi}{Definition}

\DeclareFontFamily{U}{matha}{\hyphenchar\font45}
\DeclareFontShape{U}{matha}{m}{n}{
	<5> <6> <7> <8> <9> <10> gen * matha
	<10.95> matha10 <12> <14.4> <17.28> <20.74> <24.88> matha12
	}{}
\DeclareSymbolFont{matha}{U}{matha}{m}{n}
\DeclareFontFamily{U}{mathx}{\hyphenchar\font45}
\DeclareFontShape{U}{mathx}{m}{n}{
	<5> <6> <7> <8> <9> <10>
	<10.95> <12> <14.4> <17.28> <20.74> <24.88>
	mathx10
	}{}
\DeclareSymbolFont{mathx}{U}{mathx}{m}{n}

\DeclareMathSymbol{\obot}{2}{matha}{"6B}
\DeclareMathSymbol{\bigobot}{1}{mathx}{"CB}

\newcommand{\Q}{\mathbb{Q}}
\newcommand{\Z}{\mathbb{Z}}
\newcommand{\C}{\mathbb{C}}
\newcommand{\tr}{\operatorname{tr}}
\newcommand{\End}{\operatorname{End}}
\newcommand{\Aut}{\operatorname{Aut}}
\newcommand{\Hom}{\operatorname{Hom}}
\newcommand{\ord}{\operatorname{ord}}
\newcommand{\id}{\operatorname{id}}
\newcommand{\GL}{\operatorname{GL}}
\newcommand{\modZ}[1]{\Z/(#1)}
\newcommand{\circlegrp}{\Q/\Z}
\newcommand{\minpoly}[1]{\mu_{#1}}
\newcommand{\charpoly}[1]{\chi_{#1}}
\newcommand{\finitefield}[1]{\mathbb{F}_{#1}}
\newcommand{\fieldextension}[2]{{#1}/{#2}}
\newcommand{\Galois}[1]{\operatorname{Gal}({#1})}
\newcommand{\units}[1]{#1^{\times}}
\newcommand{\reciprocal}[1]{#1^{\ast}}
\newcommand{\inquo}[2]{{#1}/{#2}}
\newcommand{\Adj}[1]{\operatorname{Adj}[{#1}]}
\newcommand{\rd}[1]{{#1}{}^{*}}
\newcommand{\characteristic}[1]{\operatorname{char}(#1)}
\newcommand{\ideal}[1]{({#1})}
\newcommand{\completion}[2]{\hat{#1}_{#2}}
\newcommand{\discr}{\operatorname{discr}}
\newcommand{\deriv}[1]{{#1}^{\prime}}
\newcommand{\modulo}[1]{\,(\operatorname{mod}\,{#1})}
\newcommand{\defeq}{\coloneqq}
\newcommand{\Hyp}[1]{\operatorname{Hyp}({#1})}
\newcommand{\restr}[2]{\left.{#1}\right|_{#2}}
\newcommand{\matrixring}[2]{M_{#1}(#2)}
\newcommand{\companion}[1]{C(#1)}
\newcommand{\norm}[1]{\operatorname{Norm}_{#1}}
\newcommand{\Arf}{\operatorname{Arf}}
\newcommand{\legendre}[2]{\bigl(\frac{#1}{#2}\bigr)}
\newcommand{\ph}{\,\cdot\,}
\newcommand{\rank}[2]{\operatorname{rank}_{#1}{#2}}
\newcommand{\abs}[1]{|{#1}|}
\newcommand{\annihilator}[2]{\operatorname{Ann}_{#1}({#2})}
\newcommand{\cusbigoplus}{\bigoplus}
\newcommand{\pub}[1]{{#1}^{\ast}}
\newcommand{\transpose}[1]{{#1}^{t}}
\newcommand{\extord}[2]{[{#1}\colon{#2}]}
\newcommand{\ass}{\operatorname{ass}}
\newcommand{\RN}[1]{\uppercase\expandafter{\romannumeral #1\relax}}
\newcommand{\orth}[1]{{#1}^{\bot}}
\newcommand{\orthvalues}[1]{\langle{#1}\rangle}
\newcommand{\orthgrp}[3]{\mathrm{O}^{#1}_{#2}({#3})}
\newcommand{\orthsum}{\obot}
\newcommand{\bigorthsum}{\bigobot}
\newcommand{\witt}[1]{\operatorname{Wq}({#1})}
\newcommand{\wittcls}[1]{[{#1}]}

\begin{document}

	\title[On the Orthogonal Group of Quadratic Modules]{On the Orthogonal Group of Quadratic Modules}
	\author[Hannes Knötzele and Sven Möller]{Hannes Knötzele\textsuperscript{\lowercase{a}} and Sven Möller\textsuperscript{\lowercase{a},\lowercase{b}}}
	\thanks{\textsuperscript{a}{Fachbereich Mathematik, Universität Hamburg, Hamburg, Germany}}
	\thanks{\textsuperscript{b}{Department of Mathematics and Statistics, Maynooth University, Maynooth, Ireland}}
	\thanks{Email: \href{mailto:hknoetzele.work@mailbox.org}{\nolinkurl{hknoetzele.work@mailbox.org}}, \href{mailto:math@moeller-sven.de}{\nolinkurl{math@moeller-sven.de}}}

	\begin{abstract}
		We classify the conjugacy classes of isometries in the orthogonal group of regular quadratic modules over Artinian local rings with finite residue fields, provided that the characteristic is odd and coprime to the order of the isometry. (More generally, the same classification holds over Henselian local rings with finite residue fields, such as the ring of $p$-adic integers.) The main idea is to lift results over finite fields to local rings by using the Henselian property.

		By applying these results to the residue class rings $\smash{\modZ{p^l}}$, we also obtain a classification of the conjugacy classes of isometries of discriminant forms,
		under the condition that the order of the discriminant form is odd and coprime to the order of the isometry.
	\end{abstract}

	\maketitle

	\setcounter{tocdepth}{1}
	\tableofcontents
	\setcounter{tocdepth}{3}


	\section{Introduction}

	Quadratic forms play a central role in many branches of mathematics. The methods used to study them depend to a large extent on the nature of the coefficients, which may belong to fields or commutative rings. A \textit{discriminant form} (also called finite quadratic form or metric group) $q\colon\Gamma\to\circlegrp$ is a regular (or nondegenerate) quadratic form defined on a finite abelian group $\Gamma$ with values in~$\circlegrp$. Discriminant forms naturally arise in the context of even lattices and theta series in number theory. They also appear as classifying objects in the cohomology theory of Eilenberg and Mac Lane \cite{MR0056295}, with applications in algebraic topology,
	braided fusion categories \cite{MR1250465}, and vertex operator algebras \cite{MR1233387}.

	\medskip

	In this text, we are concerned with the orthogonal group $\orthgrp{}{}{\Gamma,q}$ of a discriminant form, which consists of all the isometries of $(\Gamma,q)$. We shall enumerate these isometries up to conjugacy, under mild simplifying assumptions. This will then be used in \cite{paper2} (see also \cite{hannes-thesis}) to construct new modular fusion categories.

	While the automorphism group $\operatorname{Aut}(\Gamma)$ of a finite abelian group $\Gamma$ is of course easily accessible, it is a computationally difficult problem to determine the structure of $\orthgrp{}{}{\Gamma,q}\leq\operatorname{Aut}(\Gamma)$ for a discriminant form $(\Gamma,q)$ and, for instance, its conjugacy classes. One approach is to realise $(\Gamma,q)\cong L'/L$ concretely by an even lattice $L$ such that the natural map $\orthgrp{}{}{L}\to\orthgrp{}{}{\Gamma,q}$ is surjective \cite{MR0525944}. Another approach is to build up $\orthgrp{}{}{\Gamma,q}$ by orthogonal reflections \cite{scheithauer}.

	Here, we develop a new approach. Using the Henselian property of the local residue class rings $\smash{\modZ{p^l}}$ with residue field $\finitefield{p}$, we lift results for the orthogonal groups $\orthgrp{}{}{2n+1,q}\leq\GL(2n+1,\finitefield{q})$ and $\orthgrp{\pm}{}{2n,q}\leq\GL(2n,\finitefield{q})$ over finite fields~$\finitefield{q}$ (for the prime case $q=p$) to arbitrary discriminant forms.
	These orthogonal groups appear famously in the classification of the finite simple groups, and their conjugacy classes were described in \cite{MR150210} (see also \cite{MR591256,MR249519,MR268192}).

	As a corollary of our main result, \autoref{thm:class-conj-classes}, we obtain in \autoref{constr:decomp-discr-form2} a concrete and algorithmic enumeration of the conjugacy classes
	of isometries of discriminant forms of odd order, which can be easily implemented in common computer algebra systems (that have an efficient implementation of irreducible polynomials over finite fields).

	\medskip

	As a first step, we reduce discriminant forms to quadratic modules over residue class rings $\smash{\modZ{p^l}}$ (see \autoref{constr:decomp-discr-form1}) and eventually to quadratic spaces over finite fields $\finitefield{p}=\modZ{p}$. As is well known, given a discriminant form $(\Gamma,q)$, we can write it (uniquely) as an orthogonal direct sum of its $p$-primary components
	\begin{equation*}
		\Gamma=\bigorthsum_{p\text{ prime}}\Gamma^{(p)},
	\end{equation*}
	where $\Gamma^{(p)}$ is the maximal $p$-subgroup of $\Gamma$, which is again a discriminant form with the restriction $q^{(p)}$ of the quadratic form; see \autoref{lem:splitting-p-groups}. As an abelian $p$-group, we can further decompose $\smash{(\Gamma^{(p)},q^{(p)})}$ orthogonally into homocyclic components as
	\begin{align*}
		\Gamma^{(p)}=\bigorthsum_{l\in\Z_{\geq1}}\Gamma^{(p^l)}\quad\text{with}\quad\Gamma^{(p^l)}=\big(\modZ{p^l}\big)^{r_l}
	\end{align*}
	for some $r_l\in\Z_{\geq0}$; see \autoref{lem:splitting-hom-grps}. This decomposition is essentially\footnote{The decomposition is unique up to isomorphism if we only consider the structure of the $\smash{\Gamma^{(p^l)}}$ as abelian groups. If we consider the $\smash{\Gamma^{(p^l)}}$ together with the restrictions of the quadratic form, then the decomposition is not quite unique.} unique, but only up to isomorphism. However, given an isometry $\phi\in\orthgrp{}{}{\Gamma^{(p)},q^{(p)}}$, whose order we assume to be not divisible by $p$, it is always possible to find a concrete orthogonal decomposition of $\smash{\Gamma^{(p)}}$ as above such that each direct summand is $\phi$-invariant; see \autoref{cor:inv-abelian-decomp} and \autoref{constr:decomp-discr-form1}. We can therefore restrict our attention to the study of isometries of the homocyclic components $\smash{\Gamma^{(p^l)}}$.

	In the following, let us restrict to the case that $p$ is an odd prime. Then the restriction of the quadratic form to $\smash{\Gamma^{(p^l)}=(\modZ{p^l})^{r_l}}$ is a map
	\begin{equation*}
		q\colon\bigl(\modZ{p^l}\bigr)^{r_l}\to\inquo{\frac{1}{p^l}\Z}{\Z}\overset{\kappa_{p^l}}{\cong}\modZ{p^l},
	\end{equation*}
	and defining $\smash{\hat{q}(x)=\kappa_{p^l}(2q(x))}$ yields the pair $(M,\hat{q})$ with $\smash{M=\Gamma^{(p^l)}=(\modZ{p^l})^{r_l}}$ viewed as a module over the ring $\smash{A=\modZ{p^l}}$ and the quadratic form $\hat{q}$ with values in $A$; see \autoref{lem:discr-forms-and-reg-quadr-form}. The pair $(M,\hat{q})$ has the structure of a (regular) quadratic module over $A$; see \autoref{defi:quadratic-module} and below.

	\medskip

	By the above arguments, it suffices to study orthogonal groups $\orthgrp{}{}{M,q}$ for quadratic modules $(M,q)$ over the residue class rings $\smash{A=\modZ{p^l}}$. For the special case of $l=1$, i.e.\ for the finite fields $\finitefield{p}=\modZ{p}$, this reduces to the well-studied orthogonal groups $\orthgrp{}{}{2n+1,p}\leq\GL(2n+1,\finitefield{p})$ and $\orthgrp{\pm}{}{2n,p}\leq\GL(2n,\finitefield{p})$ mentioned above.

	Now, when $l>1$, we can instead consider the reduction of the quadratic module $(M,q)$ over $\smash{A=\modZ{p^l}}$ modulo the maximal ideal $\mathfrak{m}=\ideal{p}$, which results in a quadratic space $(\overline{M},\overline{q})$ over the residue field $k=A/\mathfrak{m}=\finitefield{p}$. We shall see that the reduction from $(M,q)$ to $(\overline{M},\overline{q})$ preserves much of the structure of $(M,q)$ as a quadratic module, in particular when studying the orthogonal group.

	Overall, we have thus reduced discriminant forms to quadratic spaces over the finite fields $\finitefield{p}$ of prime characteristic, restricting here to the case that $p$ is odd.

	\medskip

	A large part of this text is then devoted to lifting conjugacy classes of isometries (and the invariants used to classify them) over finite fields to isometries over rings such as $\smash{A=\modZ{p^l}}$. The key property is that the ring $\smash{A=\modZ{p^l}}$ is a local ring with maximal ideal $\mathfrak{m}=\ideal{p}$ and finite residue field $k=A/\mathfrak{m}=\finitefield{p}$
	such that coprime factorisations over $k$ lift to polynomial factorisations over $A$. Such rings are called Henselian and we shall use this lifting property repeatedly throughout the text:
	\begin{introDefi}[\textbf{\ref{defi:Henselian-ring}} (Henselian Local Ring)]
		Let $A$ be a local ring with residue field~$k$. Then $A$ is called \textit{Henselian} if for every polynomial $Q \in A[T]$ such that its reduction over $k$ factors as $\overline{Q} = R_0 S_0$ with $R_0, S_0 \in k[T]$ where $R_0$ is monic and $R_0$ and $S_0$ are relatively prime in $k[T]$, there exist polynomials $R, S \in A[T]$ such that $Q$ factors as $Q = RS$ with $R$ monic and $\overline{R} = R_0, \overline{S} = S_0$.
	\end{introDefi}
	For most of the text, we shall work in the general framework of Henselian local rings. (Also note that all rings in this text are assumed to be Noetherian.) In fact, we often work with Artinian local rings, which are complete and thus Henselian. While finite rings (such as the residue class rings $\smash{\modZ{p^l}}$) are automatically Artinian, there are also examples of infinite Henselian local rings such as the $p$-adic integers $\Z_p$ or the ring $\finitefield{q}[[T]]$ of formal power series over a finite field.

	\medskip

	The central objects of this text are quadratic modules:
	\begin{introDefi}[\textbf{\ref{defi:quadratic-module}} (Quadratic Module)]
		Let $A$ be a ring. A \textit{quadratic module} over $A$ is a pair $(M,q)$ consisting of a finitely generated projective module~$M$ and a map $q \colon M \to A$ satisfying
		\begin{enumerate}
			\item $q(a x) = a^2 q(x)$ for all $x \in M$ and $a \in A$,
			\item the map $\lambda_q(x,y) \defeq q(x + y) - q(x) - q(y)$ is a symmetric bilinear form.
		\end{enumerate}
	\end{introDefi}
	A quadratic module $(M,q)$ is called \textit{regular} (or \textit{nondegenerate}) if the associated (symmetric) Hermitian module $(M,\lambda_q)$ is regular. For a quadratic module over a field, we use the term quadratic space. For a quadratic module $(M,q)$ over $A$, we define the \textit{orthogonal group} $\orthgrp{}{A}{M,q}$ as the group of all isometries of $(M,q)$.

	We then show that many important properties of quadratic forms $q$ over Hensel\-ian local rings $A$ can be studied by looking at the reduced quadratic form $\overline{q}$ over the residue field $k = A/\mathfrak{m}$. We shall eventually assume that the residue field~$k$ is finite, i.e.\ $k=\finitefield{q}$, and of odd characteristic.

	\medskip

	As a first important result, we establish the existence of minimal polynomials of certain automorphisms over complete (e.g.\ Artinian) local rings $A$. When $A$ is not a field, the polynomial ring $A[T]$ is not a principal ideal domain, which means that minimal polynomials do not always exist.
	\begin{introTheorem}[\textbf{\ref{thm:existence-minimal-polynomial}} (Existence of Minimal Polynomials)]
		Let $A$ be a complete (e.g.\ Artinian) local ring with maximal ideal~$\mathfrak{m}$ and residue field~$k$ of characteristic~$p$. Let $M$ be a finitely generated projective module over $A$ and let $\phi \colon M \to M$ be an automorphism of order~$N$ such that $p \nmid N$. Then $\phi$ admits a minimal polynomial~$\minpoly{\phi}$ which divides the polynomial $T^N-1$ in $A[T]$.
	\end{introTheorem}
	This theorem will be crucial for our main result, the enumeration of conjugacy classes of isometries over $A$, where minimal polynomials appear as important invariants.

	The proof proceeds as follows: using Cohen's structure theorem for complete local rings, we show that complete local rings are quotients of complete regular local rings; see \autoref{cor:cohen-class}. For complete local rings, we show that solutions of matrix equations over the residue field can be lifted to solutions over the ring; see \autoref{prop:lifting-matrix-eqn}. Since matrices over regular local rings admit minimal polynomials, this then allows us to prove the theorem.

	\medskip

	Our main result is \autoref{thm:class-conj-classes}. In it, we classify the conjugacy classes of isometries of order coprime to~$p$ in the orthogonal groups of regular quadratic modules over Artinian local rings with finite residue field of odd characteristic~$p$ by a complete set of invariants, consisting essentially of irreducible polynomials over the residue field.

	For the proof, we essentially follow the same steps as in the well-understood case of finite fields (see e.g.\ \cite{MR591256}), using that for Henselian local rings with finite residue fields, we can lift many properties of their residue fields to the ring.

	The assumption that the ring is Artinian and that the order of the isometry is coprime to the characteristic of the residue field is only needed to ensure the existence of the minimal polynomials, which then divide $T^N - 1$ for the order~$N$ of the isometry, as explained above, and that their factors are canonical Hensel lifts; see \autoref{defi:canonical-lift}. In fact, such minimal polynomials exist over all Henselian local rings (see \autoref{rem:minpoly-Henselian}), so that the classification holds for all Henselian local rings with finite residue field of odd characteristic, such as the ring~$\Z_p$ of $p$-adic integers for an odd prime~$p$.

	The proof of \autoref{thm:class-conj-classes} is constructive. That is, for every set of invariants, we can construct a regular quadratic module with an isometry with these invariants, using the \textit{involution traces} naturally defined by self-reciprocal polynomials; see \autoref{rem:realisations-conj}.

	In the following, rather than giving the precise statement of \autoref{thm:class-conj-classes}, we describe the classification steps by means of minimal polynomials.
	\begin{introSpecial}[\textbf{\ref{thm:class-conj-classes}} (Parametrisation of Conjugacy Classes)]
		Let $A$ be an Artinian local ring with maximal ideal~$\mathfrak{m}$ and finite residue field~$k$ of odd characteristic~$p$. Let $(M,q)$ be a regular quadratic module over $A$. Then the conjugacy classes of isometries of $(M,q)$ of order coprime to $p$ can be parametrised as follows:
		\begin{enumerate}
			\item If two isometries admit a minimal polynomial and are similar, i.e.\ conjugate by an automorphism of~$M$, they have the same minimal polynomial; see \autoref{lem:conj-implies-same-minpoly}.
			\item The minimal polynomial $\minpoly{\phi}$ of an isometry $\phi$ of a regular quadratic module over a ring $A$ is self-reciprocal; see \autoref{lem:minpoly-self-reciproical}.
			\item As the ring $A$ is Henselian, we can decompose the minimal polynomial as $\smash{\minpoly{\phi} = P_1 \dots P_m Q_1 \reciprocal{Q}_1 \dots Q_n \reciprocal{Q}_n}$ with self-reciprocal polynomials $P_i \in A[T]$ (type~\RN{1}) and non-self-reciprocal polynomials $Q_j \in A[T]$ (type~\RN{2}) such that the reductions $\smash{\overline{P}_i}$, $\smash{\overline{Q}_j}$, $\smash{\reciprocal{\overline{Q}}_j}$ are irreducible and pairwise relatively prime; see \autoref{lem:decomp-module-ann}. Since $\minpoly{\phi}$ divides $T^N - 1$ for the order~$N$ of~$\phi$, these factors are the \textit{canonical Hensel lifts} of their reductions, i.e.\ the unique monic lifts dividing $T^N - 1$; see \autoref{defi:canonical-lift}.
			\item Using the above decomposition of the minimal polynomial, we find an orthogonal decomposition $\smash{M = (\bigorthsum_{i = 1}^{m} M_{i, \phi}^{(\mathrm{\RN{1}})}) \orthsum ( \bigorthsum_{j = 1}^{n} M_{j,\phi}^{(\mathrm{\RN{2}})} )}$ of the quadratic module $(M,q)$; see \autoref{lem:min-poly-orthgrp}.
			\item For two isometries with the same self-reciprocal minimal polynomial whose reduction is irreducible, we can show that they are conjugate using the classification of regular Hermitian modules over Henselian local rings; see \autoref{lem:ex-hermitian-from} and \autoref{lem:red-isom-hermitian-form}.
			\item We prove in \autoref{thm:equivalence-Hermitian-forms} that over a Henselian local ring with finite residue field, two regular Hermitian modules of the same rank are isometric if the reduction of the involution is nontrivial.
			\item In \autoref{sec:appendix}, we use Galois theory for Henselian local rings with finite residue fields to determine the Witt class of a regular quadratic module with an isometry whose minimal polynomial is self-reciprocal of degree at least~$2$ with irreducible reduction; see \autoref{lem:witt-cls-nontrivial}.
			\item Consequently, the conjugacy class of an isometry~$\phi$ of $(M,q)$ of order coprime to~$p$ is determined by the irreducible factors of~$\smash{\minpoly{\overline{\phi}}}$ over~$k$, their multiplicities in the characteristic polynomial of~$\smash{\overline{\phi}}$, and the isometry classes of the quadratic modules $\ker(\phi \mp \id)$. Conversely, all such data satisfying the rank and Witt class conditions in \autoref{thm:class-conj-classes} are realised by an explicit construction; see \autoref{rem:realisations-conj}. In particular, applying this also over~$k$ and using \autoref{lem:witt-iso-Henselian}, reduction modulo~$\mathfrak{m}$ induces a bijection between the conjugacy classes of isometries of order coprime to~$p$ of $(M,q)$ and of $(\overline{M},\overline{q})$; see also \autoref{rem:Pomfret}.
		\end{enumerate}
	\end{introSpecial}

	Since irreducible polynomials over finite fields of fixed degree and order can be efficiently enumerated, the above enumeration procedure can be implemented in common computer algebra systems. We can easily produce small examples; see \autoref{ex:isometries-type-one}.

	In view of certain constructions of modular fusion categories in \cite{paper2} based on the results in this text, we point out that the distinction between type~\RN{1} and type~\RN{2} will turn out to be important in that context.

	\medskip

	Finally, we combine \autoref{thm:class-conj-classes} with the decomposition of discriminant forms as described above (see \autoref{constr:decomp-discr-form1}) to obtain the main corollary, an algorithmic classification of conjugacy classes of isometries of discriminant forms. Omitting some details, this can be stated as follows:
	\begin{introConstr}[\textbf{\ref{constr:decomp-discr-form2}}]
		Let $(\Gamma,q)$ be a discriminant form with $|\Gamma|$ odd and $\phi$ an isometry of $(\Gamma,q)$ of order~$N$ coprime to $|\Gamma|$. By applying the orthogonal and $\phi$-invariant decomposition into discriminant forms
		\begin{align*}
			(\Gamma, q) = \bigorthsum_{p \mid \abs{\Gamma}} \bigorthsum_{l \in \Z_{\geq 1}} (\Gamma^{(p^l)}, q^{(p^l)})
		\end{align*}
		from \autoref{constr:decomp-discr-form1}, and by applying the classification of \autoref{thm:class-conj-classes} in terms of minimal polynomials to the regular quadratic modules $\smash{M=\Gamma^{(p^l)}=(\modZ{p^l})^{r_l}}$ over the Artinian local rings $\smash{A=\modZ{p^l}}$, we can enumerate all such isometries $\phi$ up to conjugacy.
	\end{introConstr}

	The classification procedure in the construction is illustrated in \autoref{ex:isometries-discr-form}, and again can be easily implemented in computer algebra systems.

	\medskip

	Our motivation to study isometries of discriminant forms comes from the fact that we use them in \cite{paper2} to build new braided (or modular) fusion categories. Indeed, similarly to how discriminant forms classify
	pointed braided fusion categories up to equivalence \cite{MR1250465}, their isometries describe braided $G$-crossed fusion categories, whose equivariantisations are again braided fusion categories, now typically not pointed.

	In that setting, as one extreme case, it is natural to consider generalisations of (equivariantisations of) Tambara-Yamagami categories; see e.g.\ \cite{MR4968024,GLM24a,GLM26} and \cite{paper2,hannes-thesis}, which correspond to isometries generating fixed-point free groups.

	Hence, in preparation for \cite{paper2}, we specialise the main results of this work, \autoref{thm:class-conj-classes} and \autoref{constr:decomp-discr-form2}, to isometries generating fixed-point free groups. Here, an isometry $\phi \colon (M,q) \to (M,q)$ is \textit{fixed-point free} if $\phi(x) \neq x$ for all $0 \neq x \in M$, and it generates a fixed-point free group if all its nontrivial powers are fixed-point free. The main theorems then specialise to \autoref{thm:class-fpf-isos} and \autoref{constr:decomp-discr-form3}.

	\medskip

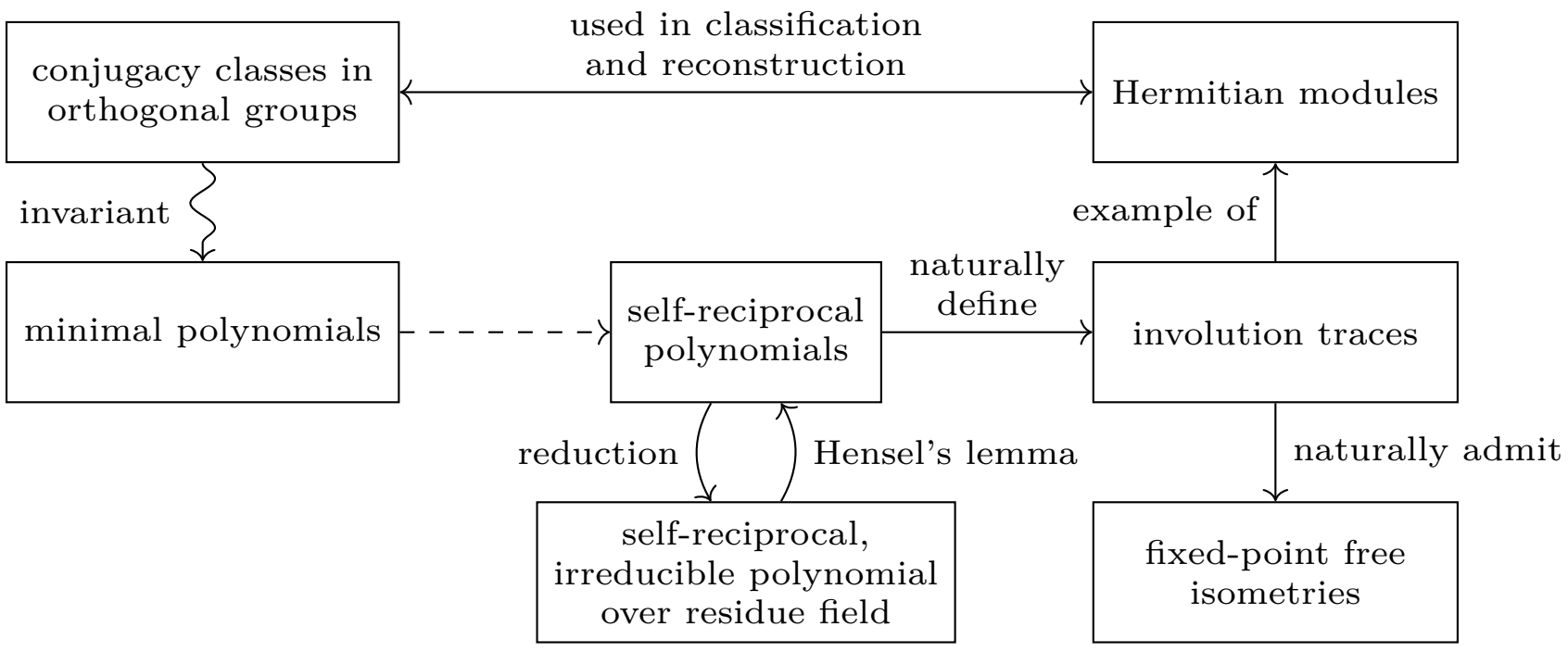
\begin{figure}[ht]
		\begin{tikzpicture}[
				font=\scriptsize,
				node distance=7mm and 15mm,
				concept/.style={draw, align=center, minimum height=10mm},
				left concept/.style={concept, minimum width=28mm},
				right concept/.style={concept, minimum width=26mm},
			]
			\node[left concept] (minpoly) {minimal polynomials};
			\node[left concept, above=of minpoly] (conj) {conjugacy classes in\\orthogonal groups};
			\node[concept, right=of minpoly] (selfrec) {self-reciprocal\\polynomials};
			\node[concept, below=of selfrec] (irred) {self-reciprocal,\\irreducible polynomial\\over residue field};
			\node[right concept, right=of selfrec] (invtrace) {involution traces};
			\node[right concept, above=of invtrace] (herm) {Hermitian modules};
			\node[right concept, below=of invtrace] (fpf) {fixed-point free\\isometries};
			\draw[<->] (conj) -- node[above, align=center] {used in classification\\and reconstruction} (herm);
			\draw[->, decorate, decoration={snake}]
				(conj) -- node[left=1mm] {invariant} (minpoly);
			\draw[->, dashed] (minpoly) -- (selfrec);
			\draw[->] (selfrec) -- node[above, align=center] {naturally\\define} (invtrace);
			\draw[->] (invtrace) -- node[left] {example of} (herm);
			\draw[->] (invtrace) -- node[right] {naturally admit} (fpf);
			\draw[->] ([xshift=-2.5mm]selfrec.south) to[bend right] node[left] {reduction} ([xshift=-2.5mm]irred.north);
			\draw[->] ([xshift=2.5mm]irred.north) to[bend right] node[right] {Hensel's lemma} ([xshift=2.5mm]selfrec.south);
		\end{tikzpicture}
		\caption{Main mathematical concepts in this work.}
		\label{fig:concepts}
	\end{figure}

	We summarise the main mathematical notions appearing in this work and how they are related in \autoref{fig:concepts}.


	\subsection*{Outline}

	In \autoref{sec:discriminant-forms}, we recall the definitions of Hermitian and quadratic modules over general rings.
	We recall some properties of quadratic modules over Henselian local rings. We also note that for a Henselian local ring~$A$ with finite residue field~$k$ of odd characteristic, its Witt ring $\witt{A}$ is isomorphic to the Witt ring $\witt{k}$ of its residue field.
	Then, we decompose a discriminant form into quadratic modules over the rings $\modZ{p^l}$.

	In \autoref{sec:minimal-polynomial}, we establish the existence of minimal polynomials of certain automorphisms over complete (e.g.\ Artinian) local rings.

	In \autoref{sec:orthogonal-grps}, we classify conjugacy classes of isometries in orthogonal groups of quadratic modules over Artinian local rings with finite residue fields.

	In \autoref{sec:fpf-isometries}, we describe the conjugacy classes of fixed-point free isometries and of isometries generating fixed-point free groups.

	Finally, in \autoref{sec:appendix}, we determine the Witt class of the quadratic module corresponding to an involution trace, which enters the Witt class condition in \autoref{thm:class-conj-classes}.


	\subsection*{Notation}

	All rings are associative, commutative, and unital. Moreover, we assume that all rings are Noetherian. For a ring $A$, let $\units{A}$ be the multiplicative group consisting of the units of $A$. For a field $k$, we denote its characteristic by $\characteristic{k}$. We denote the (residue class) ring of integers modulo a positive integer $N$ by $\modZ{N}$.

	Often, $A$ will be a Henselian local ring with finite residue field. As a special case, we frequently take $A$ to be an Artinian local ring with finite residue field. We recall that for a Noetherian local ring $A$, the following implications hold:
	\[
		\begin{tikzpicture}
			\node (A1) at (0,0) {finite};
			\node (B1) [right=of A1] at (0,0) {Artinian};
			\node (C1) [right=of B1] {complete};
			\node (D1) [right=of C1] {Henselian};
			\draw [double equal sign distance, -Implies] (A1) -- (B1);
			\draw (B1) edge [bend right = 45, double equal sign distance, -Implies] node[anchor=south] {finite residue field} (A1);
			\draw [double equal sign distance, -Implies] (B1) -- (C1);
			\draw [double equal sign distance, -Implies] (C1) -- (D1);
		\end{tikzpicture}
	\]

	Over these rings, we shall always work with finitely generated projective modules. For a local ring $A$ and a finitely generated projective module $M$, we denote its \textit{rank} (the cardinality of a minimal generating set) by $\smash{\rank{A}{M}}$. Since all our rings are local, these modules are always free modules of finite rank, i.e.\ have a basis with finitely many elements, and the cardinality of any basis is equal to the rank. For an $A$-module $M$, we denote the \textit{dual $A$-module} by $\smash{\rd{M} \defeq \Hom_{A}(M,A)}$.

	We shall often study a module, polynomial, matrix, etc.\ over a local ring by looking at its reduction to the residue field. Such a reduction will be denoted by $\smash{\overline{\ph}}$.


	\subsection*{AI Use Statement}

	No generative AI or AI-assisted tools were used for producing the results in this work.


	\subsection*{Acknowledgements}

	We thank César Galindo, Terry Gannon, Simon Lentner, Andrew Riesen, Nils Scheithauer, and Christoph Schweigert for helpful discussions. This paper (together with \cite{paper2}) is based on the PhD thesis of the first author at Universität Hamburg \cite{hannes-thesis}. We thank Terry Gannon and Christoph Schweigert for being referees for that thesis. The first author would also like to thank Karim Ritter von Merkl, Isabela Recio, and Jonas Linßen.

	Both authors acknowledge support from the DFG through the Emmy Noether Programme and the CRC 1624 \textit{Higher Structures, Moduli Spaces and Integrability}, project numbers 460925688 and 506632645.


	\section{From Discriminant Forms to Artinian Rings} \label{sec:discriminant-forms}

	In this section, we recall the definitions of discriminant forms and quadratic modules over local rings. We then use group-theoretical results to show how to reduce the study of discriminant forms to that of regular quadratic modules over the finite local rings $\modZ{p^l}$ for primes~$p$ and positive integers $l \in \Z_{\geq 1}$; see \autoref{lem:splitting-p-groups}, \autoref{lem:splitting-hom-grps}, and \autoref{lem:discr-forms-and-reg-quadr-form}.

	We recall basic definitions and conventions for Hermitian and quadratic modules; see \autoref{defi:quadratic-module}. Many of the properties of quadratic spaces over fields can be generalised to quadratic modules over local rings; see e.g.\ \cite{MR252382}. For example, Witt's cancellation theorem still holds; see \autoref{lem:witt-cancel}.

	We will work with quadratic modules over \textit{Henselian} local rings; see \autoref{subsec:henselian-rings}. In the end, we are primarily interested in the finite rings $\modZ{p^l}$ for a prime~$p$ and a positive integer $l \in \Z_{\geq 1}$, but, as we will see, the majority of the results generalise to all Henselian local rings; for the definition see \autoref{defi:Henselian-ring}. In essence, this is due to the fact that Henselian local rings have Galois-theoretic properties similar to those of their residue fields. For example, solutions of polynomial equations can be lifted from the residue field to the ring; see \autoref{lem:lifting-of-zeros}.
	Fundamentally, this is what allows us to show, using Witt's cancellation theorem, that the Witt ring of a Henselian local ring with finite residue field of odd characteristic is isomorphic to the Witt ring of its residue field; see \autoref{lem:witt-iso-Henselian}.

	\medskip


	\subsection{Discriminant Forms} \label{subsec:discr-form}
	We recall the definitions of discriminant forms, of their isometries, and of orthogonal groups. Note that we will use a similar notation for quadratic modules in \autoref{subsec:quadr-Henselian-rings}. Later, in \autoref{subsec:decomp-discr-form}, we will see why this is justified.
	\medskip

	We begin by defining \textit{discriminant forms}. The main results of the following sections also apply to quadratic modules that are not related to discriminant forms; see e.g.\ \autoref{thm:class-conj-classes}. However, the application to discriminant forms is one of the main motivations for us.
	\begin{defi}[Discriminant Form, {see e.g.\ \cite{MR1773561}}] \label{defi:discr-form}
		Let $\Gamma$ be a finite abelian group and $q \colon \Gamma \to \circlegrp$ a map such that
		\begin{enumerate}
			\item $q(a g) = a^2 q(g)$ for all $a \in \Z$ and $g \in \Gamma$,
			\item the map $\lambda_q \colon \Gamma \times \Gamma \to \circlegrp$ defined by
			\begin{align*}
				\lambda_q(g,h) \defeq q(g+h) - q(g) - q(h)
			\end{align*}
			for all $g,h \in \Gamma$ is a regular symmetric bilinear form.
		\end{enumerate}
		Then, the pair $(\Gamma,q)$ is called a \textit{discriminant form}.
	\end{defi}
	Here, the bilinear form $\lambda_q$ is called \textit{regular} (or \textit{nondegenerate}) if the map
	\begin{align*}
		\Adj{\lambda_q} \colon & \Gamma \to \Hom(\Gamma,\circlegrp), \\
		& g \mapsto \lambda_q(g, \ph)
	\end{align*}
	is an isomorphism.

	We call $g,h \in \Gamma$ orthogonal if $\lambda_q(g,h) = 0$. For a subgroup $H \leq \Gamma$, we denote its \textit{orthogonal complement} (or \textit{annihilator}) by $\orth{H}$. A subgroup $H \leq \Gamma$ is called \textit{isotropic} if $q$ restricts to $0 + \Z$ on $H$. An isotropic subgroup $H$ is called \textit{Lagrangian} if $H = \orth{H}$. An isometry between discriminant forms $(\Gamma_1, q_1)$ and $(\Gamma_2,q_2)$ is a group isomorphism $f \colon \Gamma_1 \to \Gamma_2$ such that $q_2 \circ f = q_1$. The group of all isometries from $(\Gamma,q)$ to $(\Gamma,q)$ is called the \textit{orthogonal group}; we will denote it by $\orthgrp{}{}{\Gamma,q}$.

	\medskip

	We remark that although discriminant forms are often studied as invariants of even lattices, our results are all independent of lattice realisations. They are called \textit{finite quadratic modules} in \cite{MR3101818} or \textit{finite quadratic forms} in \cite{MR0525944}.

	Discriminant forms are also equivalent to \textit{metric groups}; see \cite[Definition~2.40]{MR2609644}. Indeed, since the order of every element in the image of a quadratic map $q \colon \Gamma \to \units{\C}$ divides $\smash{\abs{\Gamma}^2}$, we can identify a quadratic map $\Gamma \to \units{\C}$ uniquely with a quadratic map $q \colon \Gamma \to \circlegrp$.


	\subsection{The Local Ring \texorpdfstring{$\modZ{p^l}$}{Z/pl} and Henselian Local Rings} \label{subsec:henselian-rings}

	In the following sections, we prove results for quadratic modules over Henselian local rings with finite residue field; sometimes, we additionally require that the ring is Artinian. The main example to keep in mind is the finite local ring $\modZ{p^l}$ for a prime~$p$ and a positive integer $l \geq 1$ since this ring will appear in the decomposition of discriminant forms; see \autoref{subsec:decomp-discr-form}.

	\medskip

	For a prime~$p$ and a positive integer $l \geq 1$, the quotient ring $\smash{A = \modZ{p^l} = \Z/p^l\Z}$ is a local ring with maximal ideal $\mathfrak{m} = \ideal{p}$ and finite residue field $k = A/\mathfrak{m} = \finitefield{p}$. For an element $x \in A$, we denote its image in the residue field by $\overline{x} \in k$.

	As a finite ring, $A$ is in particular Artinian. The maximal ideal $\smash{\mathfrak{m}=\ideal{p}}$ is nilpotent, so $A$ is separated and complete in the $\mathfrak{m}$-adic topology, i.e.\ $\smash{A \cong \completion{A}{\mathfrak{m}}}$, where we denote the completion of $A$ in the $\mathfrak{m}$-adic topology by $\completion{A}{\mathfrak{m}}$.

	\medskip

	In \autoref{subsec:quadr-Henselian-rings}, we will recall results on quadratic modules over local rings. We will put a special emphasis on Henselian local rings with finite residue fields. These rings share many of the properties of their residue fields; see e.g.\ \autoref{lem:witt-iso-Henselian} or \autoref{thm:equivalence-Hermitian-forms}. They generalise the finite rings $\modZ{p^l}$ above.
	\begin{defi}[Henselian Local Ring] \label{defi:Henselian-ring}
		Let $A$ be a local ring with residue field~$k$.
		Then $A$ is called \textit{Henselian}\footnote{Note that for a local ring $A$ with maximal ideal~$\mathfrak{m}$ that is not Noetherian, the $\mathfrak{m}$-adic topology is not necessarily separated. Since our rings are usually even Artinian and local, this will not be a problem for us.}
		if for every polynomial $Q \in A[T]$ whose reduction over~$k$ factors as $\overline{Q} = R_0 S_0$ with $R_0, S_0 \in k[T]$ where $R_0$ is monic and $R_0$ and $S_0$ are relatively prime in $k[T]$, there exist polynomials $R, S \in A[T]$ such that $Q$ factors as $Q = RS$ with $R$ monic and $\overline{R} = R_0, \overline{S} = S_0$.
	\end{defi}
	All Artinian local rings are complete and thus Henselian. In particular, the finite ring $\modZ{p^l}$ is Henselian for all primes~$p$ and positive integers $l \in \Z_{\geq 1}$.

	\smallskip

	The above lifting property will be used repeatedly throughout this paper; see e.g.\ the proof of \autoref{lem:witt-iso-Henselian}. We note that there is an extensive literature on Henselian rings; a very detailed account of Henselian rings can be found in \cite{MR277519}. We now list some equivalent definitions of Henselian rings which illustrate some of their main properties.
	\begin{lem}[Characterisations of Henselian Rings, {\cite[Theorem~\RN{1}.4.2]{MR559531}}]
		Let $A$ be a local ring with residue field~$k$.
		The following are equivalent.
		\begin{enumerate}
			\item $A$ is Henselian.
			\item Any finite $A$-algebra $B$ is a direct product of local rings.
			\item Let $P_1, \dots, P_n \in A[T_1, \dots, T_n]$ be polynomials. If there exists an $x \in k^n$ such that $\overline{P}_i(x) = 0$ for all $i = 1, \dots, n$ and $\smash{\det ( \inquo{\partial \overline{P}_i}{\partial T_j}(x) )_{i,j} \neq 0}$, then there exists a $y \in A^n$ such that $\overline{y} = x$ and $P_i(y) = 0$ for all $i = 1, \dots, n$.
		\end{enumerate}
	\end{lem}

	Here, we denoted the reduction of the coefficients of a polynomial $Q \in A[T]$ to the residue field by $\overline{Q} \in k[T]$. The polynomials $R$ and $S$ are uniquely determined with this property; see \cite[Remark~\RN{1}.4.1]{MR559531}.\footnote{To our knowledge, Azumaya named this property after Kurt Hensel in \cite{MR40287}. This property was actually already used by Theodor Schönemann in 1846; see \cite{MR1578516}.} In fact, any separated and complete local ring is Henselian; see \cite[Theorem~10.1.10]{MR3618889}.\footnote{For non-Noetherian rings, Hensel's lemma takes a slightly different form; see e.g.\ \cite[Theorem~10.1.4]{MR3618889}.}

	\smallskip

	First, we note that for Henselian local rings, we can lift solutions of polynomial equations, in particular of quadratic equations.
	\begin{lem}[{cf.\ \cite[Lemma~4]{MR40287}}] \label{lem:lifting-of-zeros}
		Let $A$ be a Henselian local ring with maximal ideal~$\mathfrak{m}$ and residue field~$k$. Let $Q$ be a polynomial in $A[T]$ and $\overline{Q}$ its reduction to $k[T]$. Suppose that $\overline{Q}$ has a simple root $\overline{a}$ in $k$. Then $Q$ has exactly one root $a$ in $A$ that is a representative of $\overline{a}$.
	\end{lem}
	Here, $Q$ need not be monic: the factorisation $\overline{Q} = (T - \overline{a}) S_0$ with $S_0(\overline{a}) \neq 0$ lifts to a factorisation $Q = (T - a) S$ by \autoref{defi:Henselian-ring}, and $a$ is unique by the uniqueness of such lifts; see below \autoref{defi:Henselian-ring}. Note that local rings are called completely primary in \cite{MR40287}.

	\smallskip

	As a direct consequence of the above lemma, we obtain:
	\begin{cor} \label{lem:square-red}
		Let $A$ be a Henselian local ring with maximal ideal~$\mathfrak{m}$ and residue field~$k$ of characteristic not~$2$. Then an element $x \in \units{A}$ is a perfect square if and only if its reduction $\overline{x} \in \units{k}$ is a perfect square.
	\end{cor}
	\begin{proof}
		If $x$ is a perfect square, then so is $\overline{x}$. Conversely, assume that $\overline{x} = \overline{x}_0^2$ for some $x_0 \in \units{A}$. Then the polynomial $Q(T) = T^2 - x$ has the reduction
		\begin{align*}
			T^2-\overline{x}_0^2 = (T-\overline{x}_0)(T+\overline{x}_0),
		\end{align*}
		and $\overline{x}_0$ is a simple root since $\overline{x}_0 \neq -\overline{x}_0$, as the characteristic of $k$ is not~$2$. Thus, by \autoref{lem:lifting-of-zeros}, there exists a unique element $m \in \mathfrak{m}$ such that
		\begin{align*}
			0 = Q(x_0 + m) = (x_0 + m)^2 - x,
		\end{align*}
		i.e.\ $x$ is a perfect square.
	\end{proof}

	\smallskip

	As pointed out, the main example of a Henselian local ring we have in mind in this paper is the finite local ring $\modZ{p^l}$ for an odd prime~$p$ and a positive integer $l \in \Z_{\geq 1}$. However, there are of course more examples of Henselian local rings.
	\begin{ex} \label{ex:Henselian-local}
		There exist Henselian local rings with finite residue field that are not finite, for example the ring of $p$-adic integers $\Z_p$ or the ring $\finitefield{q}[[T]]$ of formal power series over a finite field where $q$ is a prime power. This follows from the fact that separated and complete local rings are Henselian; see \cite[Theorem~10.1.10]{MR3618889}.

		More examples of finite local rings can be found in \cite{MR364218} (note the error pointed out in \cite{MR1180736}) or \cite{MR1778166}. For example, the quotient ring $\inquo{\finitefield{q}[T]}{\ideal{T^l}}$ is a finite local ring for all positive integers $l \in \Z_{\geq 1}$.
	\end{ex}


	\subsection{Hermitian and Quadratic Modules over Henselian Rings} \label{subsec:quadr-Henselian-rings}
	We now define Hermitian and quadratic modules over local rings. As previously mentioned, we will see how discriminant forms are related to quadratic modules over finite local rings in \autoref{subsec:decomp-discr-form}. However, we point out that our results in \autoref{sec:orthogonal-grps} may also be relevant to readers not interested in discriminant forms but in results on quadratic modules over Henselian local rings, like those in \autoref{ex:Henselian-local}.

	We will use Hermitian forms in the proof of our main theorem, \autoref{thm:class-conj-classes}; see \autoref{thm:equivalence-Hermitian-forms} and in particular \autoref{subsubsec:type-one}. As we already mentioned, this section has similarities to \autoref{subsec:discr-form}; this will be explained in \autoref{subsec:decomp-discr-form}.

	\medskip

	Recall that an \textit{involution} on a ring $A$ is a ring endomorphism $\sigma$ of $A$ of order at most $2$. Sesquilinear forms generalise bilinear forms; they are linear in the first variable and linear up to an involution in the second variable; see below. We will always use sesquilinear forms that are Hermitian, defined in \autoref{defi:Hermitian-form}.
	\begin{defi}[Sesquilinear Form, {\cite[Definition~\RN{1}.2.2]{MR1096299}}]
		Let $A$ be a ring with involution $\sigma$. Let $M$ be a finitely generated projective module over $A$. Then a \textit{$\sigma$-sesquilinear form} is a map $s \colon M \times M \to A$ satisfying
		\begin{align*}
			&& s(x+y,z) & = s(x,z) + s(y,z), & s(x, y+z) & = s(x,y) + s(x,z), \\
			&& s(ax,y) & = a s(x,y), & s(x,ay) & = \sigma(a) s(x,y)
		\end{align*}
		for all $x, y \in M$ and $a \in A$.
	\end{defi}
	The pair $(M,s)$ is called a \textit{$\sigma$-sesquilinear module}. An \textit{isometry} of $\sigma$-sesquilinear modules $(M_1, s_1)$ and $(M_2, s_2)$ over $A$ is an $A$-linear isomorphism $\phi \colon M_1 \to M_2$ such that
	\begin{align*}
		s_2(\phi(x_1), \phi(y_1)) = s_1(x_1,y_1)
	\end{align*}
	for all $x_1, y_1 \in M_1$.

	A $\sigma$-sesquilinear form $s$ on $M$ induces an $A$-linear map
	\begin{align*}
		\Adj{s} \colon M \to \rd{M}_{\sigma}, \quad \Adj{s}(x) = s(x,\ph),
	\end{align*}
	where $\rd{M}_{\sigma}$ denotes the $A$-module of $\sigma$-semilinear maps $f \colon M \to A$, i.e.\ additive maps with $f(ay) = \sigma(a) f(y)$ for all $a \in A$ and $y \in M$. For $\sigma = \id_{A}$, this is the dual module $\rd{M}$. We call $\Adj{s}$ the \textit{adjoint} of $s$. The $\sigma$-sesquilinear form $s$ is called \textit{regular} (or \textit{nondegenerate}) if its adjoint $\Adj{s}$ is an isomorphism.

	If $E$ is a commutative ring extension of $A$, an involution of $E$ with respect to $A$ is an involution $\sigma$ of $E$ that restricts to the identity on $A$, i.e.\ $\smash{\restr{\sigma}{A} = \id_{A}}$.

	\smallskip

	\textit{Hermitian forms} generalise symmetric bilinear forms. We will use \textit{Hermitian modules}, i.e.\ modules with Hermitian forms, in our classification of conjugacy classes in orthogonal groups; see e.g.\ \autoref{lem:red-isom-hermitian-form}. In particular, we will prove that all regular Hermitian modules of the same rank over a fixed Henselian local ring with finite residue field are isometric if the reduction of the involution is nontrivial; see \autoref{thm:equivalence-Hermitian-forms}. This generalises the corresponding result for finite fields; see \cite[Section~9, Remark~3]{MR43785}.
	\begin{defi}[$\epsilon$-Hermitian Form, {\cite[Definition~\RN{1}.2.3]{MR1096299}}] \label{defi:Hermitian-form}
		Let $A$ be a ring with involution $\sigma$ and $\epsilon \in A$ with $\epsilon \, \sigma(\epsilon) = 1$. Let $(M, \lambda)$ be a $\sigma$-sesquilinear module. Then $\lambda \colon M \times M \to A$ is called \textit{$\epsilon$-Hermitian} if
		\begin{align*}
			\lambda(x,y) = \epsilon \, \sigma \big( \lambda(y,x) \big)
		\end{align*}
		for all $x,y \in M$.
	\end{defi}
	In this case, $(M, \lambda)$ is called an \textit{$\epsilon$-Hermitian $\sigma$-sesquilinear module}. We introduce Hermitian modules because they are useful for classifications; see \cite{MR1957}. Of course, a $1$-Hermitian $\id_{A}$-sesquilinear form is simply called a \textit{symmetric} form.

	\medskip

	Now, we define \textit{quadratic modules}, the central objects of this work.
	\begin{defi}[Quadratic Module, {\cite[Definition~\RN{1}.5.3.5]{MR1096299}}] \label{defi:quadratic-module}
		Let $A$ be a ring. A~\textit{quadratic module} over $A$ is a pair $(M,q)$ consisting of a finitely generated projective module $M$ with a map $q \colon M \to A$ satisfying
		\begin{enumerate}
			\item $q(a x) = a^2 q(x)$ for all $x \in M$ and $a \in A$,
			\item the map $\lambda_q \colon M \times M \to A$ defined by
			\begin{align*}
				\lambda_q(x,y) \defeq q(x + y) - q(x) - q(y)
			\end{align*}
			is a symmetric bilinear form.
		\end{enumerate}
	\end{defi}
	If $2$ is a unit in $A$, we can identify quadratic forms and symmetric bilinear forms by associating with a symmetric bilinear form $\lambda$ the quadratic form $q_{\lambda}$ defined by $q_{\lambda}(x) = \lambda(x,x)$ and associating with a quadratic form $q$ the symmetric bilinear form $\frac{1}{2} \lambda_q$. The associated bilinear form is sometimes called the \textit{polar form}.

	A quadratic module $(M,q)$ is called \textit{regular} (or \textit{nondegenerate}) if the (symmetric) Hermitian module $(M,\lambda_q)$ is regular.

	For a quadratic module over a field, we will use the term \textit{quadratic space}.

	\smallskip

	An \textit{isometry} of quadratic modules $(M_1,q_1)$ and $(M_2,q_2)$ over $A$ is an $A$-linear isomorphism $\phi \colon M_1 \to M_2$ such that $q_2 \circ \phi = q_1$. For a quadratic module $(M,q)$, we denote the group of all isometries from $(M,q)$ to $(M,q)$ by $\orthgrp{}{A}{M,q}$. This group is called the \textit{orthogonal group} of the quadratic module $(M,q)$ over $A$.

	We will limit ourselves to local rings $A$ for which every finitely generated projective module is a free module of finite rank.

	For a finitely generated projective module $U$ over $A$, we define on the direct sum $U \oplus \rd{U}$ of $U$ and its dual module $\rd{U}$ a quadratic map $q_U \colon U \oplus \rd{U} \to A$ by
	\begin{align*}
		q_U(x + \phi) = \phi(x)
	\end{align*}
	for all $x \in U$ and $\phi \in \rd{U}$. The regular quadratic module $\Hyp{U} \defeq (U \oplus \rd{U}, q_U)$ is called the \textit{hyperbolic} module associated with $U$; see \cite[Section~\RN{1}.4]{MR491773}.

	We denote the Witt ring of regular quadratic modules over $A$ by $\witt{A}$, cf.\ \cite[Definition~\RN{1}.4.8]{MR491773}. For a regular quadratic module $(M,q)$ over $A$, we denote its Witt class by $\wittcls{(M,q)} \in \witt{A}$.

	For two regular quadratic modules $(M_1,q_1)$ and $(M_2,q_2)$, $\wittcls{(M_1,q_1)} = \wittcls{(M_2,q_2)}$ holds in the Witt ring if and only if there exist two finitely generated projective modules $U_1$ and $U_2$ over $A$ such that there is an isometry
	\begin{align*}
		(M_1, q_1) \orthsum \Hyp{U_1} \cong (M_2, q_2) \orthsum \Hyp{U_2}.
	\end{align*}

	In \autoref{sec:orthogonal-grps}, we will describe orthogonal groups of regular quadratic modules up to isometry, so naturally the Witt ring will play a major role in that section; see \autoref{thm:class-conj-classes}.

	Let $A$ be a local ring with maximal ideal~$\mathfrak{m}$ and residue field~$k$. For a quadratic module $(M,q)$ over $A$, we can always reduce $(M,q)$ modulo $\mathfrak{m}$ to obtain a quadratic space $(\overline{M}, \overline{q})$ over $k$.
	Moreover, the quadratic module $(M,q)$ is regular if and only if its reduction $(\overline{M},\overline{q})$ is regular; see e.g.\ \cite[Proposition~\RN{1}.1.4]{MR491773}.

	\medskip

	Before we end this section, we note that in the proof of our classification theorem of conjugacy classes, we will use the following fundamental theorem of Galois theory for Henselian local rings. The Galois theory of these rings is more or less controlled by the Galois theory of their residue fields.

	Recall that a \textit{connected} ring is a ring with no nontrivial idempotent elements (i.e.\ not equal to $1$ or $0$). In the proof, we use some general facts about Henselian local rings from \cite{MR3618889}.
	\begin{prop}[Fund.\ Theorem of Galois Theory for Henselian Local Rings] \label{prop:Galois-Henselian-rings}
		Let $A$ be a Henselian local ring with maximal ideal~$\mathfrak{m}$ and residue field~$k$. Let $C$ be a connected, separable $A$-algebra, finitely generated as an $A$-module, such that $C^G = A$ for the group $G = \Aut_{A}(C)$ of $A$-algebra automorphisms of~$C$. Then $G$ is finite, and $H \mapsto C^H$ is a bijective, order-inverting correspondence between the subgroups $H$ of $G$ and the subrings $B$ of $C$ that contain $A$ and are separable over~$A$.

		Moreover, $C$ and all these subrings $B$ are local and Henselian.
	\end{prop}
	\begin{proof}
		The first part of the proposition is the fundamental theorem of Galois theory and holds for any connected ring; see \cite[Theorem~12.5.4]{MR3618889}. We only need to prove that $C$ and the subrings $B$ are Henselian local.

		First note that since $A$ is Henselian and $C$ is finitely generated as an $A$-module, we obtain by \cite[Proposition~10.1.6]{MR3618889} that there is an isomorphism of sets of idempotents of $C$ and $C \otimes_{A} k$. Thus, as $C$ is connected, $C \otimes_{A} k$ is connected, i.e.\ it has only trivial idempotents. Since $C$ is separable over $A$, it also follows that $C \otimes_{A} k$ is a commutative separable $k$-algebra. Thus, $C \otimes_{A} k $ is a separable field extension of $k$. This implies that $C$ is local by \cite[Lemma~10.1.1]{MR3618889}. Therefore, $C$ is Henselian by \cite[\href{https://stacks.math.columbia.edu/tag/04GH}{Tag~04GH}]{stacks-project}.

		If there is a subring $B$ of $C$ containing $A$ and separable over $A$, then the ring~$B$ is connected and finitely generated as an $A$-module (as $A$ is Noetherian), so, by the above argument, it is local and Henselian.
	\end{proof}


	\subsection{Quadratic Modules over Henselian Local Rings with Finite Residue Fields} \label{subsec:quadr-modules-finite}

	We show in this section that the most important properties of quadratic forms over \textit{Henselian} local rings can be studied by looking at the reduced quadratic form over the residue field. We will work with local rings with \textit{finite} residue fields of \textit{odd} characteristic.

	\smallskip

	Let us fix some notation and recall the structure of the Witt ring of $\finitefield{q}$; see e.g.\ \cite[Corollary~\RN{2}.3.6]{MR396410}:
	\begin{align*}
		\witt{\finitefield{q}} = \begin{cases}
			\modZ{2}[\inquo{(\units{\finitefield{q}})}{(\units{\finitefield{q}})^2}], &\text{if } q \equiv 1 \modulo{4}, \\
			\modZ{4}, &\text{if } q \equiv 3 \modulo{4}.
		\end{cases}
	\end{align*}

	We consider the two maps
	\begin{align*}
		\begin{aligned}
			d \colon \\
			{}
		\end{aligned} & \begin{aligned}
			& \witt{\finitefield{q}} \to \modZ{2}, \\
			& \wittcls{(M,q)} \mapsto n \modulo{2},
		\end{aligned} \\
		\begin{aligned}
			\Arf \colon \\
			{}
		\end{aligned} & \begin{aligned}
			&\witt{\finitefield{q}} \to \inquo{(\units{\finitefield{q}})}{(\units{\finitefield{q}})^2},\\
			& \wittcls{(M,q)} \mapsto (-1)^{n(n-1)/2} \discr(M,q)
		\end{aligned}
	\end{align*}
	where $n = \smash{\rank{\finitefield{q}}{M}}$. We will refer to the latter invariant, i.e.\ the signed discriminant, as the \textit{Arf invariant}. Here, the discriminant is the determinant of the matrix $\lambda_q(x_i,x_j)$ for a basis $(x_i)$ of~$M$. It is determined up to $(\units{\finitefield{q}})^2$. Together, these two maps completely describe $\witt{\finitefield{q}}$. Similarly to \cite{MR591256}, we denote by $\Lambda$ the element of $\witt{\finitefield{q}}$ with $d(\Lambda) = 0$ and $\smash{\legendre{\Arf(\Lambda)}{\finitefield{q}} = -1}$. Here, $\smash{\legendre{\ph}{\finitefield{q}}}$ denotes the Legendre symbol for the finite field $\finitefield{q}$; see e.g.\ \cite[Section~1.14]{MR591256}. For the finite field $\finitefield{p}$ with $p$ a prime, it coincides with the Kronecker symbol.

	\medskip

	At the end of this section, we show that for a Henselian local ring $A$ with finite residue field~$k$ of odd characteristic, the reduction map induces an isomorphism of Witt rings $\witt{A} \cong \witt{k}$; see \autoref{lem:witt-iso-Henselian}.
	Thus, many properties of quadratic forms over finite fields hold for Henselian local rings with finite residue fields as well. We collect some properties of quadratic modules over Henselian rings. Our primary examples are quadratic modules over the ring $\modZ{p^l}$ for an odd prime~$p$ and some integer $l \geq 1$.

	General facts about quadratic modules can be found in \cite{MR0152507}, \cite{MR2496771}, and \cite{MR591256}.
	For quadratic modules over local rings, we refer to \cite{MR491773}. Finally, for finite fields, we use the textbooks \cite{MR1189139} and \cite{MR396410}.

	In the following, we quickly recall some basic facts relating quadratic modules over a ring with quadratic spaces over its residue field; see e.g.\ \cite{MR491773}.

	\medskip

	We recall that over finite fields of odd characteristic, every quadratic space admits an orthogonal basis. This generalises to local rings with finite residue fields of odd characteristic.
	\begin{lem}[{see e.g.\ \cite[Proposition~\RN{1}.3.4]{MR491773} or \cite[Theorem~1.6.4]{MR770063}}] \label{lem:orth-basis-quadr-module}
		Let $A$ be a local ring with residue field of \textit{odd} characteristic and let $(M,q)$ be a regular quadratic module over $A$. Then $(M,q)$ admits an orthogonal basis.
	\end{lem}
	Let us first fix some notation. For units $a_1, \dots, a_m \in \units{A}$, we denote a quadratic module of rank $m$ with orthogonal basis $(x_1, \dots, x_m)$ by $\orthvalues{a_1, \dots, a_m}$ where $a_i$ is the value of the quadratic form at $x_i$ for all $i$; see e.g.\ again \cite[Proposition~\RN{1}.3.4]{MR491773}.

	Let $(M,q)$ be a quadratic module over $A$ and $a \in \units{A}$ a unit in $A$. Then $(M,q)$ \textit{represents} $a$ if there is an element $x \in M$ such that $q(x) = a$. A quadratic module representing every unit is called \textit{universal}.

	We need the fact that regular quadratic modules of rank at least $2$ over a Henselian local ring with maximal ideal~$\mathfrak{m}$ and \textit{finite} residue field of odd characteristic are universal. Therefore, we briefly recall some facts about hyperbolic forms.

	In the study of quadratic spaces over finite fields, one uses Witt's theorem and the following lemma to show that an anisotropic quadratic space (i.e.\ a quadratic space whose quadratic form is nonzero for every nonzero vector) is of dimension $0$, $1$, or $2$; see \cite[p.~138]{MR1189139}.
	We will use it together with \autoref{lem:isotropic-reduction} to prove in \autoref{lem:isotropic-implies-universal} that isotropic regular quadratic modules over Henselian local rings with finite residue field are universal.
	\begin{lem}[{e.g.\ \cite[Theorem~\RN{1}.3.4]{MR396410} and \cite[Theorem~11.2]{MR1189139}}] \label{lem:finite-isotropic}
		Let $k$ be a finite field of odd characteristic. Let $(M,q)$ be a regular quadratic space over~$k$. Then the quadratic space $(M,q)$ over~$k$ is isotropic, i.e.\ it contains a nonzero element $x \in M$ with $q(x) = 0$, if and only if it contains a hyperbolic plane.

		If $M$ is of dimension at least $3$, $q$ is isotropic.
	\end{lem}

	We already noted that a quadratic module over a Henselian local ring is regular if and only if its reduction is regular. The following lemma establishes a similar relation for the property of being isotropic. Here, as in \cite[Lemma~\RN{5}.1.4]{MR491773}, a quadratic module $(M,q)$ over a local ring~$A$ with maximal ideal~$\mathfrak{m}$ is called \textit{isotropic} if $q(x) = 0$ for some \textit{primitive} element $x \in M$, i.e.\ $x \notin \mathfrak{m}M$; over a field, this is the notion from \autoref{lem:finite-isotropic}.
	\begin{lem}[{cf.\ \cite[Lemma~\RN{5}.1.4]{MR491773}}] \label{lem:isotropic-reduction}
		Let $A$ be a Henselian local ring with maximal ideal~$\mathfrak{m}$ and \textit{finite} residue field~$k$ of odd characteristic. Let $(M,q)$ be a regular quadratic module over $A$. Then the reduction $\smash{(\overline{M},\overline{q})}$ over $k$ is isotropic if and only if $(M,q)$ is isotropic.
	\end{lem}
	Note that the result in \cite[Lemma~\RN{5}.1.4]{MR491773} is stated for complete local rings, but the proof only uses the lifting properties of Henselian local rings; see \autoref{lem:lifting-of-zeros}. Moreover, only the implication from $\smash{(\overline{M},\overline{q})}$ to $(M,q)$ is stated there; the converse is clear since the reduction of a primitive element is nonzero. We omit the proof as it is similar to the proof of the following lemma.

	Using the lifting properties of Henselian local rings (see \autoref{lem:lifting-of-zeros}) and the description of isotropic spaces over finite fields (see \autoref{lem:finite-isotropic}), we can show:
	\begin{lem}[{cf.\ \cite[Lemma~\RN{5}.1.4]{MR491773}}] \label{lem:isotropic-implies-universal}
		Let $A$ be a Henselian local ring with maximal ideal~$\mathfrak{m}$ and \textit{finite} residue field~$k$ of odd characteristic. Let $(M,q)$ be a regular quadratic module over $A$. Then, if $(M,q)$ is isotropic, it is universal.
	\end{lem}
	\begin{proof}
		Let $a \in \units{A}$. Since $(M,q)$ is isotropic, so is its reduction $(\overline{M},\overline{q})$ over $k$ by \autoref{lem:isotropic-reduction}. By \autoref{lem:finite-isotropic}, $(\overline{M},\overline{q})$ contains a hyperbolic plane $(\overline{e},\overline{f})$. We find elements $e, f \in M$ whose residue classes modulo $\mathfrak{m}$ are $\overline{e}$ and $\overline{f}$, respectively. Then the polynomial
		\begin{align*}
			P \defeq q(e) T^2 + \lambda_q(e,f) T + q(f)- a \in A[T]
		\end{align*}
		satisfies
		\begin{align*}
			\overline{P} = T - \overline{a} \in k[T]
		\end{align*}
		and $\deriv{\overline{P}} = 1$. So, by \autoref{lem:lifting-of-zeros}, there exists an $x \in A$ such that $P(x) = 0$. Thus,
		\begin{equation*}
			q(xe + f) = q(e) x^2 + \lambda_q(e,f) x + q(f) = P(x) + a = a. \qedhere
		\end{equation*}
	\end{proof}

	To prove that regular quadratic modules of rank at least $2$ are universal, we only need the following fact. For finite fields, it is well known (see e.g.\ \cite[Corollary~\RN{1}.3.5]{MR396410}), but we generalise it to quadratic modules over Henselian rings, again (implicitly) using their lifting properties; see \autoref{lem:lifting-of-zeros}.
	\begin{lem} \label{lem:repr-isotropy}
		Let $A$ be a Henselian local ring with maximal ideal~$\mathfrak{m}$ and \textit{finite} residue field~$k$ of odd characteristic. Let $(M,q)$ be a regular quadratic module over~$A$. Then a unit $a \in \units{A}$ is represented by $q$ if and only if $q \orthsum \langle -a \rangle$ is isotropic.
	\end{lem}
	\begin{proof}
		If $a$ is represented by $q$, i.e.\ $a = q(x)$ for some $x \in M$, then $q(x)-a 1^2 = 0$, so $q \orthsum \langle -a \rangle$ is isotropic.

		Conversely, let $q \orthsum \langle -a \rangle$ be isotropic. Then there are $x \in M$ and $b \in A$, with $x \notin \mathfrak{m}M$ or $b \notin \mathfrak{m}$, such that
		\begin{align*}
			0 = q(x) - a b^2.
		\end{align*}
		If $b \in \units{A}$, then $a = q(b^{-1}x)$, i.e.\ $a$ is represented by $q$. On the other hand, if $b$ is not invertible, then $b \in \mathfrak{m}$, so that $x \notin \mathfrak{m}M$ and $\overline{q}(\overline{x}) = 0$. Then $\overline{q}$ is isotropic. By \autoref{lem:isotropic-reduction} and \autoref{lem:isotropic-implies-universal}, $q$ is universal and, in particular, $a$ is represented by~$q$.
	\end{proof}

	Finally, we are in a position to prove that all regular quadratic modules over Henselian rings with finite residue field of odd characteristic of rank at least $2$ are universal. For finite fields, this is a standard fact; see \cite{MR396410}.
	\begin{lem} \label{lem:quadr-forms-over-Henselian-rings-are-universal}
		Let $A$ be a Henselian local ring with maximal ideal~$\mathfrak{m}$ and \textit{finite} residue field of odd characteristic. Let $(M,q)$ be a regular quadratic module over $A$ of rank $\rank{A}{M} \geq 2$. Then $(M,q)$ is universal.
	\end{lem}
	\begin{proof}
		Let $a \in \units{A}$. Then $q \orthsum \langle -a \rangle$ is of rank at least $3$, so its reduction is isotropic by \autoref{lem:finite-isotropic}, and hence $q \orthsum \langle -a \rangle$ is isotropic by \autoref{lem:isotropic-reduction}. Thus, $q$ is universal by \autoref{lem:repr-isotropy}.
	\end{proof}

	\medskip

	Let us quickly recall the characterisation of quadratic separable algebras over local rings; for details we refer to \cite[Section~4.3]{MR279159} or \cite[Theorem~\RN{2}.1.6]{MR491773}. We will use quadratic extensions in our proof of \autoref{thm:equivalence-Hermitian-forms}; see \autoref{lem:norm-universal} below.
	\begin{prop}[{\cite[Theorem~\RN{2}.1.6]{MR491773}}]
		Let $A$ be a local ring. Then any quadratic separable $A$-algebra has the form
		\begin{align*}
			B = \inquo{A[T]}{\ideal{T^2-T-b}}
		\end{align*}
		for some $b \in A$ with $1 + 4b \in \units{A}$.
	\end{prop}
	We denote the image of $T$ in $B$ by $t$. Then, the (canonical) involution of $B$ keeping $A$ fixed is given by
	\begin{align*}
		\sigma \colon B \to B, \quad t \mapsto 1-t.
	\end{align*}
	Note that
	\begin{align*}
		\sigma(t^2) = 1 - t + b = \sigma(t+b).
	\end{align*}
	We define the \textit{norm map} $n \colon B \to A$ by $n(x) \defeq x \sigma(x)$ for all $x \in B$. The quadratic module $(B,n)$ is called the \textit{norm space}.

	We already proved in \autoref{lem:quadr-forms-over-Henselian-rings-are-universal} that regular quadratic modules over Henselian local rings are universal if the rank is at least $2$. In particular, the above norm space is universal. We will use this fact in \autoref{subsec:equiv-Hermitian-modules}; see \autoref{thm:equivalence-Hermitian-forms}.
	\begin{lem} \label{lem:norm-universal}
		Let $A$ be a Henselian local ring with maximal ideal~$\mathfrak{m}$ and \textit{finite} residue field of odd characteristic, and let $B$ be a quadratic separable $A$-algebra. Then the norm space $(B,n)$ over $A$ is regular and universal.
	\end{lem}
	\begin{proof}
		The regularity of $(B,n)$ follows from the regularity of the reduction. The universality follows from \autoref{lem:quadr-forms-over-Henselian-rings-are-universal}.
	\end{proof}

	\medskip

	We will end this section by proving that the Witt ring of a Henselian local ring with finite residue field of odd characteristic is isomorphic to the Witt ring of the residue field; see \autoref{lem:witt-iso-Henselian}.

	\smallskip

	The proof will use the fact that there is an orthogonal basis (see \autoref{lem:orth-basis-quadr-module}) and the famous Witt cancellation property, which we recall now.
	\begin{lem}[Witt Cancellation over Local Rings] \label{lem:witt-cancel}
		Let $A$ be a local ring, and let $E$, $F$, and $G$ be regular quadratic modules over $A$. Then, if there exists an isometry $E \orthsum F \cong E \orthsum G$, this implies $F \cong G$.
	\end{lem}
	We refer to \cite[Corollary~\RN{3}.4.3]{MR491773} or \cite{MR252382}.

	\medskip

	Now, we prove the aforementioned isomorphism of Witt rings. This lemma already hints at the close connection between quadratic modules over a Henselian ring and those over its residue field. For Artinian local rings, we will investigate this relationship further in \autoref{sec:orthogonal-grps} when we study the orthogonal group of quadratic modules over an Artinian local ring. Recall that an Artinian local ring is Henselian local.
	\begin{lem}[see e.g.\ {\cite[Theorem~6.2.4]{MR770063}}] \label{lem:witt-iso-Henselian}
		Let $A$ be a Henselian local ring with finite residue field~$k$ of odd characteristic. Then the reduction
		induces an isomorphism of Witt rings
		\begin{align*}
			\witt{A} \to \witt{k}.
		\end{align*}
	\end{lem}
	In \cite{MR770063}, the lemma is stated for symmetric bilinear forms over discrete valuation rings, but the proof carries over to our setting; for complete local rings, see also \cite[Theorem~1.6.13]{MR770063}. For the convenience of the reader, we reproduce the proof here.
	\begin{proof}
		Since every regular quadratic module over $A$ or $k$ admits an orthogonal basis by \autoref{lem:orth-basis-quadr-module}, and since the reduction of a hyperbolic module $\Hyp{U}$ is the hyperbolic space $\Hyp{\overline{U}}$, it suffices to show that two regular quadratic modules over $A$ with isometric reductions are isometric. We proceed by induction on the rank of the quadratic modules. Consider two quadratic modules $\orthvalues{a_1, \dots, a_m}$ and $\orthvalues{b_1, \dots, b_m}$ of rank $m$, and suppose there is an isometry $\smash{\orthvalues{\overline{a}_1, \dots, \overline{a}_m} \cong \orthvalues{\overline{b}_1, \dots, \overline{b}_m}}$ of quadratic spaces over $k$. Thus, there exist $x_i \in A$ such that
		\begin{align*}
			\overline{a}_1 \overline{x}_1^2 + \dots + \overline{a}_m \overline{x}_m^2 = \overline{b}_1
		\end{align*}
		and at least one of the $\overline{x}_i$ is nonzero and thus $x_i \in \units{A}$. Without loss of generality, we suppose $x_1 \in \units{A}$. Then the polynomial
		\begin{align*}
			T^2 - a_1^{-1} (b_1-a_2 x_2^2-\dots-a_m x_m^2)
		\end{align*}
		has a (simple) root $x \in A$ by \autoref{lem:lifting-of-zeros}. Thus, if $(e_1, \dots, e_m)$ is an orthogonal basis of $\orthvalues{a_1, \dots, a_m}$ with $q(e_i) = a_i$, the element $v \defeq x e_1 + x_2 e_2 + \dots + x_m e_m$ satisfies $q(v) = b_1$, and (the free module of rank~$1$) $A v$ is a regular submodule of $\orthvalues{a_1,\dots,a_m}$. We obtain an isometry
		\begin{align*}
			\orthvalues{a_1, \dots, a_m} \cong \orthvalues{b_1, c_2, \dots, c_m}.
		\end{align*}
		By Witt's cancellation theorem for $k$, \autoref{lem:witt-cancel}, $\orthvalues{\overline{c}_2, \dots, \overline{c}_m} \cong \smash{\orthvalues{\overline{b}_2, \dots, \overline{b}_m}}$. Hence, $\orthvalues{c_2, \dots, c_m} \cong \orthvalues{b_2, \dots, b_m}$ by the induction hypothesis. It follows that $\orthvalues{a_1, \dots, a_m} \cong \orthvalues{b_1, \dots, b_m}$.
	\end{proof}
	The above lemma holds for $2$-Henselian local rings (called \textit{quadratisch henselsch} in \cite{MR252382}).


	\subsection{Decomposition of Discriminant Forms} \label{subsec:decomp-discr-form}
	In this section, we decompose discriminant forms, using results by Maschke and Schur, to obtain regular quadratic modules over the local rings $\modZ{p^l}$ for primes~$p$ and positive integers $l \in \Z_{\geq 1}$.

	\medskip

	Of course, every discriminant form splits orthogonally into its $p$-subgroups.
	\begin{lem}[{see e.g.\ \cite{MR3101818}}] \label{lem:splitting-p-groups}
		Let $(\Gamma,q)$ be a discriminant form. Then $(\Gamma,q)$ splits orthogonally into its $p$-subgroups $(\Gamma^{(p)},q^{(p)})$.
	\end{lem}
	Here, $\smash{q^{(p)}}$ denotes the restriction of the quadratic form $q$ to $\smash{\Gamma^{(p)}}$.

	\smallskip

	By the next lemma, a discriminant form on a $p$-group splits orthogonally into homocyclic components.
	\begin{lem}[Splitting Lemma, {\cite[Lemma~8]{MR0156890}}] \label{lem:splitting-hom-grps}
		Let $(\Gamma,q)$ be a discriminant form such that the underlying finite abelian group $\Gamma$ is a $p$-group for a prime~$p$. Then there is a decomposition of $\Gamma$ as a finite abelian group into pairwise orthogonal subgroups
		\begin{align*}
			\Gamma = \cusbigoplus_{l \in \Z_{\geq 1}} \Gamma^{(p^l)},
		\end{align*}
		where $\Gamma^{(p^l)}$ is trivial or a homocyclic finite abelian group of exponent $p^l$ for a prime~$p$ and a positive integer $l \in \Z_{\geq 1}$, i.e.\ a finite abelian group
		\begin{align*}
			\Gamma^{(p^l)} = \big(\modZ{p^l}\big)^{r_l}
		\end{align*}
		for some $r_l \in \Z_{\geq 0}$.

		For every such decomposition, the restriction $(\Gamma^{(p^l)},q^{(p^l)})$ of $(\Gamma,q)$ to $\Gamma^{(p^l)}$ is a discriminant form and determines an orthogonal splitting
		\begin{align*}
			(\Gamma,q) = \displaystyle \bigorthsum_{l \in \Z_{\geq 1}} (\Gamma^{(p^l)},q^{(p^l)}).
		\end{align*}
	\end{lem}
	Note that Wall has a slightly different notion of quadratic modules. The result also appears in \cite{MR3101818}.

	The proof of \cite[Lemma~8]{MR0156890} shows the following, which we will need in \autoref{constr:decomp-discr-form1}.
	\begin{rem} \label{rem:splitting-top-component}
		Let $(\Gamma,q)$ be a discriminant form on a $p$-group $\Gamma$ of exponent $p^L$, and let $\Gamma = H \oplus H'$ be a (not necessarily orthogonal) decomposition as finite abelian groups such that $H$ is homocyclic of exponent $p^L$ and $p^{L-1} H' = 0$. Then $(H, \restr{q}{H})$ is a discriminant form, and hence $\smash{(\Gamma,q) = (H, \restr{q}{H}) \orthsum (\orth{H}, \restr{q}{\orth{H}})}$ with $\orth{H} \cong H'$ as finite abelian groups.
		Also, if $H$ is invariant under a group of isometries of $(\Gamma,q)$, then so is $\orth{H}$.
	\end{rem}

	\smallskip

	The following lemma should be understood as a change of coefficients of a quadratic map.
	\begin{lem}[{see e.g.\ \cite[Lemmas~2.7 and 2.10]{MR3101818}}] \label{lem:discr-forms-and-reg-quadr-form}
		Let $p$ be an odd prime and $l \geq 1$ a positive integer.
		There is a bijection between discriminant forms $(\Gamma, q)$ whose underlying group $\Gamma$ is homocyclic of exponent $p^l$ and regular quadratic modules $(M, \hat{q})$ over the ring $A = \modZ{p^l}$ where we denote $\Gamma$ viewed as an $A$-module by $M$. With the group isomorphism
		\begin{align*}
			\kappa_{p^l} \colon & \inquo{\frac{1}{p^l} \Z}{\Z} \to \modZ{p^l}, \\
			& \frac{x}{p^l} \modulo{1} \mapsto x \modulo{p^l},
		\end{align*}
		this bijection is defined by the map $q \mapsto \hat{q}$,
		\begin{align*}
			\hat{q}(x) = \kappa_{p^l}(2 q(x)).
		\end{align*}
	\end{lem}
	A regular quadratic map on a finite abelian group $\Gamma$ of odd order thus decomposes into regular quadratic modules of the form $(M,\hat{q})$ with $M=(\modZ{p^l})^{r}$ and quadratic form $\hat{q} \colon (\modZ{p^l})^{r} \to \modZ{p^l}$ for odd primes~$p$ and $l,r\in \Z_{\geq 1}$. Note that in the situation of the lemma, $p^{2l} q(x) = q(p^l x) = 0$ and $2 p^l q(x) = p^l \lambda_q(x,x) = \lambda_q(p^l x, x) = 0$ for all $x \in \Gamma$, so that $\smash{q(x) \in \inquo{\frac{1}{p^l}\Z}{\Z}}$ since $p$ is odd. In particular, $\hat{q}$ is well defined.

	\medskip

	By the above \autoref{lem:splitting-p-groups} and \autoref{lem:splitting-hom-grps}, a discriminant form is determined, up to isometry, by an orthogonal sum of homocyclic components of exponents $p^l$ for primes~$p$ and positive integers $l \in \Z_{\geq 1}$, each of which can be identified, by \autoref{lem:discr-forms-and-reg-quadr-form}, with a regular quadratic module over the ring $\modZ{p^l}$. By \autoref{lem:witt-iso-Henselian}, such a regular quadratic module is completely determined by its rank and Arf invariant.
	\begin{nota} \label{nota:discr-form-symbol}
		We denote a discriminant form $\smash{(\Gamma^{(p^l)}, q^{(p^l)})}$ for which the finite abelian group $\smash{\Gamma^{(p^l)}}$ is homocyclic with exponent $p^l$ and for which the regular quadratic module corresponding to $\smash{(\Gamma^{(p^l)}, q^{(p^l)})}$ under the bijection in \autoref{lem:discr-forms-and-reg-quadr-form} has rank~$r$ and Arf invariant $\smash{x \in \inquo{\units{\finitefield{p}}}{(\units{\finitefield{p}})^2}}$ such that $\smash{\legendre{x}{\finitefield{p}} = \epsilon \in \{\pm 1\}}$ by $\smash{(p^l)^{\epsilon r}}$.

		We denote the orthogonal sum of two discriminant forms $\smash{(p_1^{l_1})^{\epsilon_1 r_1}}$ and $\smash{(p_2^{l_2})^{\epsilon_2 r_2}}$ by $\smash{(p_1^{l_1})^{\epsilon_1 r_1} (p_2^{l_2})^{\epsilon_2 r_2}}$.
	\end{nota}
	If $p = p_1 = p_2$ and $l = l_1 = l_2$, the relations on the Witt ring $\smash{\witt{\modZ{p^{l}}}}$ allow us to add two symbols. The precise relations depend on the prime $p$, namely
	\begin{align*}
		\smash{(p^{l})^{\epsilon_1 r_1} (p^{l})^{\epsilon_2 r_2}} = \smash{(p^{l})^{\epsilon (r_1 + r_2)}}
	\end{align*}
	where
	\begin{align*}
		\epsilon \defeq \legendre{(-1)^{r_1 r_2} x_1 x_2}{\finitefield{p}}, \quad \text{with } \legendre{x_i}{\finitefield{p}} = \epsilon_i.
	\end{align*}
	In other words, rank modulo~$2$ and Arf invariant define a group isomorphism between $\smash{\witt{\modZ{p^l}} \cong \witt{\finitefield{p}}}$ and $\smash{\modZ{2} \times \inquo{\units{\finitefield{p}}}{(\units{\finitefield{p}})^2}}$ with the group law given by $(e,x) \cdot (e',x') = (e + e', (-1)^{ee'} xx')$, an extension of $\smash{\inquo{\units{\finitefield{p}}}{(\units{\finitefield{p}})^2}}$ by $\modZ{2}$; see \cite[Proposition~\RN{2}.2.1 and Theorem~\RN{2}.3.5]{MR396410}.

	\smallskip

	We illustrate the notation in the following example.
	\begin{ex} \label{ex:discr-symbol}
		A discriminant form with the symbol
		\begin{align*}
			(3)^{-4} (7^2)^{+2}
		\end{align*}
		has the underlying abelian group $\smash{\big(\modZ{3}\big)^{4} \oplus \big(\modZ{7^2}\big)^{2}}$.

		Note that by \autoref{lem:splitting-p-groups}, we can split this discriminant form into $\smash{(3)^{-4}}$ and $\smash{(7^2)^{+2}}$ and then write down all the orthogonal decompositions of these discriminant forms into discriminant forms of smaller rank (viewed as modules over $\modZ{3}$ and $\modZ{7^2}$, respectively).

		In this case, we obtain, up to reordering, the nine decompositions
		\begin{align*}
			(3)^{-4} & = (3)^{+3} (3)^{+1} = (3)^{-3} (3)^{-1} = (3)^{+2} (3)^{-2} \\
			& = (3)^{-2} (3)^{+1} (3)^{-1} = (3)^{+2} (3)^{-1} (3)^{-1} = (3)^{+2} (3)^{+1} (3)^{+1} \\
			& = (3)^{+1} (3)^{-1} (3)^{-1} (3)^{-1} = (3)^{-1} (3)^{+1} (3)^{+1} (3)^{+1}
		\end{align*}
		and the two decompositions
		\begin{align*}
			(7^2)^{+2} = (7^2)^{+1} (7^2)^{-1}.
		\end{align*}

		For example, the discriminant form with symmetric bilinear form with
		\begin{align*}
			\begin{bmatrix}
				\frac{2}{3} & \frac{2}{3} & \frac{2}{3} & \frac{1}{3} & \frac{2}{7^2} & \frac{2}{7^2}
			\end{bmatrix}
		\end{align*}
		on the diagonal has the associated symbol above.
	\end{ex}

	\medskip

	We briefly mention that this is a particular Jordan decomposition.
	\begin{rem}
		This decomposition is closely related to the \textit{Jordan decomposition} of discriminant forms; see e.g.\ \cite[Section~2]{MR2496771}, which is based on the diagonalisation of symmetric bilinear forms. Note that the sign of the Jordan decomposition corresponds (up to a normalisation) to the Arf invariant in our notation. However, the sign in \cite[Section~2]{MR2496771} is multiplicative while the Arf invariant is not necessarily multiplicative.
	\end{rem}

	Note that the above change of coefficients has the following small caveat.
	\begin{rem}
		There are Lagrangian subgroups of discriminant forms which do not correspond to maximally isotropic submodules of the associated quadratic module above, simply because they do not split off as direct summands. Consider, for example, the characteristic subgroup $p(\modZ{p^2})$ in a discriminant form with underlying abelian group $\modZ{p^2}$. However, every totally isotropic direct summand~$U$ of the quadratic module with $U = \orth{U}$ naturally yields a Lagrangian subgroup.
	\end{rem}

	\medskip

	A well-known result by Maschke and Schur will allow us to reduce the study of isometries of a discriminant form $(\Gamma,q)$ to those of its homocyclic components.

	\smallskip

	The following proposition is fundamental in representation theory. There are several closely related theorems which are usually stated for fields; see \cite[Theorem~\RN{1}.17.7]{MR4928289}. We will use this result for an abelian $p$-group $H$ viewed as a module over a (finite) group $\Phi$ of automorphisms of $H$ whose order is coprime to~$p$.
	\begin{prop}[{Maschke, Schur, see e.g.\ \cite[Theorem~\RN{1}.17.6]{MR4928289} and \cite{MR51838}}] \label{prop:maschke}
		Let $\Pi$ be a finite group and $M$ a $\Pi$-module. Let $\sigma \colon M \to M$ be an endomorphism (of $\Pi$-modules) such that $\sigma(\abs{\Pi}x) = x$ for all $x \in M$. Let $M = M_1 \oplus M_2$ be a direct-sum decomposition as abelian groups such that $M_1$ is a $\Pi$-submodule of $M$. Then there exists a $\Pi$-submodule $M_2' \subseteq M$ such that $M = M_1 \oplus M_2'$ as $\Pi$-modules.
	\end{prop}

	As we already mentioned, we now apply this proposition repeatedly to an abelian $p$-group $H$ viewed as a module over a (finite) group $\Phi$ of automorphisms of $H$ whose order is coprime to $p$.
	\begin{cor} \label{cor:inv-abelian-decomp}
		Let $H$ be an abelian $p$-group and $\Phi$ be a group of automorphisms of order not divisible by $p$. Then there is a decomposition into $\Phi$-invariant subgroups
		\begin{align*}
			H = \cusbigoplus_{l} H_{l}
		\end{align*}
		such that every $H_l$ is homocyclic with exponent $p^l$.
	\end{cor}
	To prove the above corollary, one applies \autoref{prop:maschke}, inductively and to suitable $\Phi$-invariant subgroups that are direct summands of~$H$ as abelian groups, to the morphism of $\Phi$-modules
	\begin{align*}
		\sigma \colon & H \to H, \\
		& x \mapsto \abs{\Phi}^{-1} x.
	\end{align*}
	Here, we use that the order~$\abs{\Phi}$ is coprime to $p$. For details, we refer to the references below.

	\smallskip

	This statement appears as \cite[Exercises~\RN{1}.68--69]{MR4928289}. A proof can be found in \cite[Proposition~6.5]{mayr-thesis}; for an elementary proof of this fact, we refer to \cite[Lemmas~2.1 and 6.1]{MR27759}.

	As we are just changing the coefficients of the quadratic map, there is a bijection between the isometries of the discriminant form $(\Gamma, q)$ and those of $(M, \hat{q})$. This is very helpful for us as it will allow us to reduce the study of isometries of (odd) discriminant forms to isometries of regular quadratic modules over rings $\modZ{p^l}$ for odd primes~$p$ and positive integers $l \geq 1$.

	\begin{rem}
		We note that a similar decomposition based on the classification results of \cite[Section~6]{MR0156890}, later completed in \cite{MR0743731}, appeared in \cite[Lemma~3.2]{MR0480333}; see also \cite[Theorem~4.7]{MR1990929}. Since we want an orthogonal decomposition invariant under a group of isometries, we need a slight modification using the above proposition due to Maschke and Schur together with \autoref{rem:splitting-top-component}.
	\end{rem}

	\medskip

	We now describe how \autoref{cor:inv-abelian-decomp} and \autoref{rem:splitting-top-component} yield an orthogonal decomposition into invariant subgroups which are all homocyclic.
	\begin{constr} \label{constr:decomp-discr-form1}
		Let $(\Gamma, q)$ be an odd-order discriminant form with a faithful action~$\tau$ of the cyclic group $\modZ{N}$ by isometries such that $N$ is coprime to the order of $\Gamma$. We set $\phi \defeq \tau[1]$. We denote the subgroup of $\orthgrp{}{}{\Gamma,q}$ generated by $\phi$ by~$\Phi$; it is of course a cyclic group of order~$N$.

		By \autoref{lem:splitting-p-groups}, we can decompose $(\Gamma, q)$ orthogonally into $p$-groups
		\begin{align*}
			(\Gamma, q) = \bigorthsum_{p \mid \abs{\Gamma}} (\Gamma^{(p)}, q^{(p)})
		\end{align*}
		such that $\smash{\Gamma^{(p)}}$ is a $p$-group for the primes~$p$ dividing the order of $\Gamma$.

		Since $N$ is coprime to the order of $\Gamma$, for every prime~$p$, we can apply \autoref{cor:inv-abelian-decomp} to the $p$-group $\smash{\Gamma^{(p)}}$ and (the restriction of) $\Phi$ to obtain a decomposition (as finite abelian groups)
		\begin{align*}
			\Gamma^{(p)} = \cusbigoplus_{l \in \Z_{\geq 1}} H_l
		\end{align*}
		into $\Phi$-invariant subgroups $H_l$ that are trivial or homocyclic of exponent $p^l$. This decomposition is not necessarily orthogonal. However, if $p^L$ denotes the exponent of $\smash{\Gamma^{(p)}}$, then by \autoref{rem:splitting-top-component}, applied to $H = H_L$ and $H' = \cusbigoplus_{l < L} H_l$, the subgroup $\smash{\Gamma^{(p^L)} \defeq H_L}$ with the restriction $\smash{q^{(p^L)}}$ of $q$ is a discriminant form and splits off orthogonally, i.e.\ $\smash{\Gamma^{(p)} = \Gamma^{(p^L)} \orthsum \orth{\Gamma^{(p^L)}}}$, where we take the orthogonal complement in $\smash{\Gamma^{(p)}}$. Since $\Phi$ acts by isometries, $\smash{\orth{\Gamma^{(p^L)}}}$ is again $\Phi$-invariant, and its exponent divides $p^{L-1}$. Repeating this argument with $\smash{\orth{\Gamma^{(p^L)}}}$ in place of $\smash{\Gamma^{(p)}}$, we obtain for every prime~$p$ an orthogonal decomposition into discriminant forms
		\begin{align*}
			(\Gamma^{(p)}, q^{(p)}) = \bigorthsum_{l \in \Z_{\geq 1}} (\Gamma^{(p^l)}, q^{(p^l)}),
		\end{align*}
		where $\smash{\Gamma^{(p^l)}}$ is trivial or a homocyclic group of exponent $p^l$, i.e.\ a finite direct sum of the cyclic group $\modZ{p^l}$, such that \textit{every $\smash{\Gamma^{(p^l)}}$ is $\Phi$-invariant}. Only finitely many of the summands $\smash{\Gamma^{(p^l)}}$ are nontrivial.

		By changing the coefficients with \autoref{lem:discr-forms-and-reg-quadr-form}, the restriction of $\tau$ to $\smash{(\Gamma^{(p^l)},q^{(p^l)})}$ yields an action of $\modZ{N}$ by isometries on a regular quadratic module $\smash{(M^{(p^l)}, \hat{q}^{(p^l)})}$ over $\smash{\modZ{p^l}}$. Let us denote the generator of this action by $\smash{\phi^{(p^l)}}$. This action need not be faithful: the order~$\smash{N^{(p^l)}}$ of $\smash{\phi^{(p^l)}}$ divides~$N$, and $N$ is the least common multiple of the $\smash{N^{(p^l)}}$ since $\tau$ is faithful; the induced action of $\smash{\modZ{N^{(p^l)}}}$ is faithful. (In view of \autoref{constr:decomp-discr-form3} below, we note that if $\tau$ is fixed-point free, then $\smash{N^{(p^l)}} = N$ for every nontrivial $\smash{\Gamma^{(p^l)}}$, and the restricted actions are again fixed-point free.)
	\end{constr}
	Our later results will, in particular, apply to the regular quadratic module $\smash{(M^{(p^l)}, \hat{q}^{(p^l)})}$ over $\modZ{p^l}$ with isometry $\smash{\phi^{(p^l)}}$; see \autoref{thm:class-conj-classes} and, if $\tau$ is fixed-point free, \autoref{thm:class-fpf-isos}, for a prime~$p$ and a positive integer $l \in \Z_{\geq 1}$.

	\smallskip

	We will come back to this construction after our main theorem, \autoref{thm:class-conj-classes}, in \autoref{constr:decomp-discr-form2}, and to its specialisation to fixed-point free actions in \autoref{constr:decomp-discr-form3}.


	\section{Minimal Polynomials for Automorphisms over Complete
	Rings} \label{sec:minimal-polynomial}

	In this section, our first main result concerns the existence of minimal polynomials over complete (e.g.\ Artinian) local rings. Over a local ring $A$ that is not a field, the polynomial ring $A[T]$ is not a principal ideal domain. Therefore, minimal polynomials of endomorphisms do not always exist.

	By Cohen's structure theorem for complete local rings, we can show that complete local rings are quotients of complete regular local rings; see \autoref{cor:cohen-class}. For complete local rings, we show that solutions of matrix equations over the residue field can be lifted to solutions over the ring; see \autoref{prop:lifting-matrix-eqn}. Since matrices over regular local rings admit minimal polynomials, this will allow us to prove our first main result, \autoref{thm:existence-minimal-polynomial}, establishing the existence of minimal polynomials over complete local rings.

	The existence of minimal polynomials proved in this section is necessary for all the results in \autoref{sec:orthogonal-grps}; see e.g.\ \autoref{lem:conj-implies-same-minpoly} and \autoref{lem:orth-decom-minpoly}.

	\medskip


	\subsection{Minimal Polynomials and Separable Polynomials} \label{subsec:sep-poly}
	In this section, we recall the definition of the minimal polynomial of an endomorphism. Furthermore, we also recall the definition of a separable polynomial. We will later work with automorphisms of finite order coprime to the characteristic of the residue field, so, if they admit minimal polynomials, these are separable polynomials; see \autoref{thm:existence-minimal-polynomial}.

	\medskip

	We primarily work with Artinian local rings and repeatedly use their completeness. Our first goal is to show that certain automorphisms over complete local rings admit \textit{minimal polynomials}.

	Let $M$ be a finitely generated projective module over a ring $A$ and let $\phi \colon M \to M$ be an $A$-endomorphism. Then, the kernel of the ring map
	\begin{align*}
		\theta \colon & A[T] \to \End_A(M), \\
		& T \mapsto \phi
	\end{align*}
	is an ideal. If this ideal is a principal ideal, it has a monic generator, the \textit{minimal polynomial} $\minpoly{\phi}$ of $\phi$; see \cite[Lemma~4.1]{MR3262438}.

	If $A$ is not a field, the existence of minimal polynomials is not guaranteed since $A[T]$ is not a principal ideal domain. For a discussion of this problem, we refer to \cite{MR3262438}, where the authors extend the theory of endomorphisms of vector spaces to endomorphisms of finitely generated projective modules.

	\medskip

	We recall the following characterisation of separable polynomials. We repeatedly use the \textit{discriminant} of a polynomial. We will explicitly work with discriminants in \autoref{sec:appendix}.
	\begin{defi}[Discriminant of a Polynomial, {see \cite[Définition~\RN{4}.6.2]{MR643362}}]
		Let $A$ be a ring. Let $Q \in A[T]$ be a monic polynomial of degree $d$ and $E = \inquo{A[T]}{\ideal{Q}}$. We denote the image of $T$ in $E$ by $t$. Then, the \textit{discriminant} of $Q$ is defined as the discriminant
		\begin{align*}
			\discr(Q) \defeq \discr_{\fieldextension{E}{A}}(1,t,\dots,t^{d-1}) \defeq \det (\tr_{\fieldextension{E}{A}} t^{i+j})_{i,j = 0, \dots, d-1}.
		\end{align*}
	\end{defi}

	Discriminants can be calculated in many computer algebra systems.
	For quadratic polynomials there is the well-known formula:
	\begin{ex}
		Let $A$ be a ring and $Q = T^2 + a T + b \in A[T]$ a quadratic polynomial. Then the discriminant of $Q$ is given by $\discr(Q) = a^2-4b$.
	\end{ex}

	\smallskip

	Using the above discriminant, we obtain a computational criterion to check whether a polynomial is separable.
	\begin{lem}[Characterisation of Separable Polynomials]
		Let $A$ be a ring. Let $Q \in A[T]$ be a monic polynomial and $E = \inquo{A[T]}{\ideal{Q}}$. We denote the image of $T$ in $E$ by $t$. Then the following are equivalent.
		\begin{enumerate}
			\item The algebra $E$ is étale over $A$.
			\item $\deriv{Q}(t)$ is invertible in $E$.
			\item The discriminant $\discr(Q)$ of $Q$ is invertible in $A$.
		\end{enumerate}
	\end{lem}
	A monic polynomial $Q$ satisfying the above conditions is called \textit{separable}.
	Details can be found in e.g.\ \cite[Proposition~2.8]{MR3262438}. We also refer to \cite[Section~4.6.1]{MR3618889} or \cite{MR280479}.

	\smallskip

	We quickly record the following helpful fact, which follows directly from the calculation of the discriminant of a product of polynomials.
	\begin{lem} \label{lem:div-separable}
		Let $A$ be a local ring with maximal ideal~$\mathfrak{m}$. Let $Q \in A[T]$ be a separable polynomial. If $R \in A[T]$ is a monic polynomial dividing $Q$, then $R$ is also separable.
	\end{lem}
	\begin{proof}
		As $Q$ is separable, $\discr(Q) \in \units{A}$. Since $A$ is (Noetherian and) local, this is equivalent to $\discr(Q) \notin \mathfrak{m}$. As $R$ divides $Q$, $\discr(R) \mid \discr(Q)$; see e.g.\ \cite[Proposition~\RN{4}.6.11, Corollaire~1]{MR643362}. Therefore, $\discr(R) \notin \mathfrak{m}$.
	\end{proof}

	Furthermore, since taking discriminants commutes with change of rings (see e.g.\ \cite[Proposition~\RN{4}.6.11]{MR643362}), we obtain:
	\begin{lem} \label{lem:div-sep}
		Let $A$ be a local ring with maximal ideal~$\mathfrak{m}$ and residue field~$k$. Let $Q \in A[T]$ be a monic polynomial with reduction $\overline{Q} \in k[T]$.
		Then $Q$ is separable if and only if $\overline{Q}$ is separable.
	\end{lem}

	We will use the following criterion for minimal polynomials later. Note that for endomorphisms of a finitely generated projective module, the characteristic polynomial always exists.
	\begin{lem} \label{lem:min-separable}
		Let $A$ be a local ring. Let $M$ be a finitely generated projective module over $A$ and let $\phi \colon M \to M$ be an $A$-endomorphism. Let $Q \in A[T]$ be a separable polynomial dividing the characteristic polynomial of $\phi$ such that $Q(\phi) = 0$. Then $Q$ is the minimal polynomial of $\phi$.
	\end{lem}
	A proof can be found in \cite[Proposition~4.2]{MR3262438}.


	\subsection{Lifting Matrix Equations for Complete Local Rings} \label{subsec:lift-matrix-eqn}
	To prove the existence of the minimal polynomial of an automorphism of finite order over a complete local ring in \autoref{subsec:cohen-min-poly}, we will need a generalisation of (a variant of) Hensel's lemma for matrix equations; see \autoref{prop:lifting-matrix-eqn}. To ensure convergence and thus the existence of a solution, we will work with complete local rings in this section. Recall, however, that every Artinian local ring is complete.

	\medskip

	We follow the ideas of \cite{MR220757}. We will show the existence of a minimal polynomial by lifting matrix equations and using the completeness of our local ring. We recall a fact about matrix equations, which we will use repeatedly. It is also stated in \cite[before Lemma~1]{MR220757}.
	\begin{lem} \label{lem:separable-invertible}
		Let $k$ be a field and $Q \in k[T]$ a polynomial. Let $X \in \matrixring{n}{k}$ be a matrix such that $Q(X) = 0$. Then, if $Q$ is separable over $k$,
		\begin{align*}
			\det{\deriv{Q}(X)} \neq 0.
		\end{align*}
	\end{lem}
	\begin{proof}
		Let $L$ be a field such that $Q$ splits in $L$. For every eigenvalue $\lambda$ of $X$,
		\begin{align*}
			Q(\lambda) = 0, \quad \deriv{Q}(\lambda) \neq 0.
		\end{align*}
		Thus, $\det{\deriv{Q}(X)} \neq 0$.
	\end{proof}

	\medskip

	We will see that we can solve matrix equations over complete rings by lifting solutions over the residue field; see \autoref{prop:lifting-matrix-eqn} below. Note that if $A$ is an Artinian local ring with maximal ideal~$\mathfrak{m}$ and $\mathfrak{a} \subseteq \mathfrak{m}$ is an ideal in $A$, then $\mathfrak{a}$ is nilpotent, so $A$ is complete (and separated) in the $\mathfrak{a}$-adic topology.

	The following proposition is based on \cite[Section~2]{MR220757}. The proof presented there is for the $p$-adic integers $\Z_p$. We will apply this proposition to Noetherian complete regular local rings coming from Cohen's structure theorem for complete local rings; see \autoref{cor:cohen-class}.
	\begin{prop}[Lifting Solutions of Matrix Equations] \label{prop:lifting-matrix-eqn}
		Let $A$ be a complete local ring with maximal ideal~$\mathfrak{m}$ and residue field~$k$. Let $\mathfrak{a} \subseteq \mathfrak{m}$ be an ideal in $A$ such that $A$ is complete in the $\mathfrak{a}$-adic topology.

		Let $Q \in A[T]$ be a separable polynomial. Let $X_0 \in \matrixring{n}{\inquo{A}{\mathfrak{a}}}$ be a matrix over the quotient ring $\inquo{A}{\mathfrak{a}}$ such that
		\begin{align*}
			Q(X_0) \equiv 0 \modulo{\mathfrak{a}}.
		\end{align*}
		Then there exists a matrix $X \in \matrixring{n}{A}$ such that
		\begin{align*}
			Q(X) = 0, \quad X \equiv X_0 \modulo{\mathfrak{a}}.
		\end{align*}
	\end{prop}
	\begin{proof}
		We write the polynomial $Q$ as $\smash{Q = \sum_{i = 0}^d a_i T^i}$.

		Note that $\smash{\deriv{Q}(X_0) \modulo{\mathfrak{a}}}$ viewed as a matrix over the quotient ring $\smash{\inquo{A}{\mathfrak{a}}}$ is invertible if and only if its reduction $\smash{\deriv{\overline{Q}}(\overline{X}_0)}$ over $k$ is invertible (which follows by considering determinants). Since $Q$ is separable, so is its reduction $\smash{\overline{Q}}$ modulo $\mathfrak{m}$. Therefore, $\smash{\deriv{\overline{Q}}(\overline{X}_0)}$ is invertible by \autoref{lem:separable-invertible}.

		Now, we lift $X_0$ to a matrix in $\matrixring{n}{A}$, which we again denote by $X_0$, and apply the same construction as in the proof of \cite[Lemma~1]{MR220757}, namely, we define
		\begin{align*}
			X_1 \defeq X_0 - Q(X_0) (\deriv{Q}(X_0))^{-1} \modulo{\mathfrak{a}^2}.
		\end{align*}
		Since $Q(X_0) (\deriv{Q}(X_0))^{-1}$ has entries in~$\mathfrak{a}$ and commutes with~$X_0$, as it is a polynomial in~$X_0$ by the Cayley-Hamilton theorem, it follows that
		\begin{align*}
			Q(X_1) & = Q(X_0 - Q(X_0)(\deriv{Q}(X_0))^{-1}) \\
			& = \sum_{i = 0}^d a_i (X_0 - Q(X_0)(\deriv{Q}(X_0))^{-1})^i \\
			& \equiv a_0 + \sum_{i = 1}^d a_i (X_0^i - i X_{0}^{i-1}Q(X_0)(\deriv{Q}(X_0))^{-1}) \modulo{\mathfrak{a}^2} \\
			& = Q(X_0) - Q(X_0) = 0.
		\end{align*}
		Similarly, for $i \geq 1$, we define
		\begin{align*}
			X_{i+1} \defeq X_i - Q(X_i) (\deriv{Q}(X_i))^{-1}.
		\end{align*}
		This is well defined since, for all $i \in \Z_{\geq 1}$, $\deriv{Q}(X_i) \equiv \deriv{Q}(X_0) \modulo{\mathfrak{m}}$ is invertible.
		It follows that
		\begin{align*}
			X_{i+1} \equiv X_i \modulo{\mathfrak{a}^{i+1}}, \quad Q(X_{i+1}) \equiv 0 \modulo{\mathfrak{a}^{i+2}}.
		\end{align*}
		Hence, for fixed $r,s$, the sequence of entries $(X_i)_{rs}$ converges to some $\smash{x_{rs} \in \completion{A}{\mathfrak{a}} \cong A}$. The matrix $X = (x_{rs})$ satisfies the above requirements.
	\end{proof}

	The above proposition allows us to construct matrices with the same order coprime to the characteristic of the residue field, which are solutions of a particular polynomial equation.
	\begin{ex}
		For the ring $A = \modZ{7^3}$ and the ideal $\ideal{7}$ in $A$, the matrix $\left[\begin{smallmatrix} 0 & 6 \\ 1 & 4 \end{smallmatrix}\right]$ in $\matrixring{2}{\finitefield{7}}$ of order~$8$ has the minimal polynomial $T^2 + 3 T + 1 \in \finitefield{7}[T]$, whose canonical Hensel lift is $Q \defeq T^2 + 108 T + 1 \in A[T]$; see \autoref{defi:canonical-lift}.
		Applying the above algorithm, we obtain the matrices $X_1 = \left[\begin{smallmatrix} 21 & 20 \\ 36 & 18 \end{smallmatrix}\right] \in \matrixring{2}{\modZ{7^2}}$ and $X_2 = \left[\begin{smallmatrix} 217 & 69 \\ 330 & 18 \end{smallmatrix}\right]$ in $\matrixring{2}{A}$ of order~$8$ such that $Q(X_1) \equiv 0 \modulo{7^2}$ and $Q(X_2) = 0$.
	\end{ex}


	\subsection{Cohen's Structure Theorem and Minimal Polynomials for Automorphisms over Complete Local Rings} \label{subsec:cohen-min-poly}
	Finally, we are in a position to prove the existence of the minimal polynomial of an automorphism of finite order over a complete local ring. We will use Cohen's celebrated structure theorem for complete local rings \cite{MR16094} and the lifting properties proved in the previous section; see \autoref{prop:lifting-matrix-eqn}.

	\medskip

	We prove that automorphisms of finite order not divisible by the characteristic of the residue field admit minimal polynomials. Using Cohen's structure theorem for complete local rings, we reduce the statement to regular local rings. We follow \cite[Section~28]{MR266911}; for a nice proof using formal smoothness we refer to \cite[Théorème~19.8.8]{MR173675}.
	Let $A$ be a complete local ring with maximal ideal~$\mathfrak{m}$ and residue field~$k$ of characteristic $p = \characteristic{k}$. Then, a \textit{coefficient ring} $R \subset A$ is a subring such that
	\begin{enumerate}
		\item $R$ is a Noetherian complete local ring with maximal ideal $R \cap \mathfrak{m}$,
		\item the residue field of $R$ is isomorphic to $k$ (by the canonical map),
		\item $R \cap \mathfrak{m} = p R$.
	\end{enumerate}
	With this notation, the structure theorem for complete local rings takes the following form.
	\begin{prop}[Cohen's Structure Theorem] \label{prop:cohen-structure}
		Let $A$ be a complete local ring with residue field~$k$ of characteristic~$p$. Then $A$ has a coefficient ring $R$.

		More precisely, $R$ can be chosen as the image of a local ring homomorphism $V \to A$ inducing an isomorphism of residue fields, where $V = k$ if $p = 0$, and $V$ is a \textit{Cohen ring} if $p > 0$, i.e.\ a complete discrete valuation ring whose maximal ideal is generated by~$p$.
	\end{prop}

	In fact, we will use the following consequence.
	\begin{cor} \label{cor:cohen-class}
		Let $A$ be a complete local ring whose maximal ideal~$\mathfrak{m}$ is
		generated by $x_1, \dots, x_n$, and let $V \to A$ be as in \autoref{prop:cohen-structure}. Then the extension
		\begin{align*}
			V[[T_1, \dots, T_n]] & \to A,\\
			T_i & \mapsto x_i,
		\end{align*}
		of $V \to A$ is a surjective map of local rings. In particular, $A$ is a quotient of the complete regular local ring $V[[T_1, \dots, T_n]]$.
	\end{cor}
	We refer to \cite[Section~28.P, Corollary~1]{MR266911}; a proof of \autoref{prop:cohen-structure} and \autoref{cor:cohen-class} can be found in \cite[\href{https://stacks.math.columbia.edu/tag/0323}{Tag~0323}]{stacks-project}. Here, the ring $V[[T_1, \dots, T_n]]$ is regular by induction on~$n$ since $V$ is regular and $T_n$ is a nonzerodivisor with $V[[T_1, \dots, T_n]]/(T_n) \cong V[[T_1, \dots, T_{n-1}]]$; see \cite[\href{https://stacks.math.columbia.edu/tag/00NU}{Tag~00NU}]{stacks-project}.

	\begin{rem}
		The ring $\modZ{p^l}$ has characteristic $p^l$ while its residue field has characteristic~$p$. Its coefficient ring is $R = \modZ{p^l}$ itself, which is not regular for $l \geq 2$, and in \autoref{prop:cohen-structure} we can take for~$V$ the ring $\Z_{p}$ of $p$-adic integers, which is a Cohen ring with $\modZ{p^l} \cong \inquo{\Z_{p}}{\ideal{p^l}}$. For the general theory, we refer to \cite[Section~\RN{2}.6]{MR354618}.
	\end{rem}

	\medskip

	By \cite[Theorem~5]{MR103906}, regular local rings are unique factorisation domains and unique factorisation domains are integrally closed (see e.g.\ \cite[Example~5.0]{MR242802}; the same argument applies to all unique factorisation domains). By \cite[Corollary]{MR2099718} every square matrix with coefficients in an integrally closed domain admits a minimal polynomial with coefficients in that domain. Thus, we get:
	\begin{lem} \label{lem:minpoly-regular-local}
		Let $A$ be a regular local ring. Let $X \in \matrixring{n}{A}$ be a matrix. Then $X$ admits a minimal polynomial $\minpoly{X} \in A[T]$.
	\end{lem}

	\medskip

	Now, we are in a position to prove the existence of minimal polynomials over complete local rings. Note that an Artinian local ring is in particular complete.
	\begin{thm}[Existence of Minimal Polynomials] \label{thm:existence-minimal-polynomial}
		Let $A$ be a complete local ring with maximal ideal~$\mathfrak{m}$ and residue field~$k$ of characteristic~$p$. Let $M$ be a finitely generated projective module over $A$ and let $\phi \colon M \to M$ be an automorphism of order~$N$ such that $p \nmid N$. Then $\phi$ admits a minimal polynomial $\minpoly{\phi}$, which divides the polynomial $T^N-1$ in $A[T]$.
	\end{thm}
	\begin{proof}
		First, the polynomial $T^N - 1$ is separable over $k$ since $p \nmid N$.
		Note that since $A$ is local, the module $M$ is free, so we can represent $\phi$ as a matrix $X_0 \in \matrixring{n}{A}$ with $X_0^N - 1 = 0$. By \autoref{cor:cohen-class}, there is a complete regular local ring $D$ with a surjective local map $\pi \colon D \to A$. We denote the kernel of this map by $\mathfrak{a}$ and the maximal ideal of $D$ by $\mathfrak{n}$. Then, we get by \cite[\href{https://stacks.math.columbia.edu/tag/090T}{Tag~090T}]{stacks-project} that $D$ is complete in the $\mathfrak{a}$-adic topology. Now, we can apply \autoref{prop:lifting-matrix-eqn} and find a matrix $X \in \matrixring{n}{D}$ that satisfies
		\begin{align*}
			X^N - 1 = 0 \text{ and } X \equiv X_0 \modulo{\mathfrak{a}}.
		\end{align*}

		Since matrices over regular local rings have minimal polynomials by \autoref{lem:minpoly-regular-local}, $X$ has a minimal polynomial $\minpoly{X} \in D[T]$, which divides $T^N - 1$ since $X^N - 1 = 0$. Since $T^N-1$ is separable, so is $\minpoly{X}$ by \autoref{lem:div-separable}.

		By reducing the coefficients, we get a monic, separable polynomial
		\begin{align*}
			Q \defeq \pi(\minpoly{X}) \in A[T]
		\end{align*}
		that satisfies $Q(X_0) = 0$.

		Since determinants commute with change of rings, we obtain the equality of characteristic polynomials $\charpoly{X_0} = \pi(\charpoly{X})$ and therefore $Q \mid \charpoly{X_0}$. Thus, the statement follows from \autoref{lem:min-separable}.
	\end{proof}

	\smallskip

	Alternatively, the theorem can be deduced directly from Hensel's lemma.
	\begin{rem} \label{rem:minpoly-Henselian}
		\autoref{thm:existence-minimal-polynomial} holds more generally for every Henselian local ring~$A$ with residue field~$k$ of characteristic~$p$, and the reduction of $\smash{\minpoly{\phi}}$ is the minimal polynomial of~$\smash{\overline{\phi}}$. This can be shown directly, without Cohen's structure theorem and \autoref{prop:lifting-matrix-eqn}.

		Indeed, let $\mu_0 \in k[T]$ be the minimal polynomial of~$\smash{\overline{\phi}}$. Since $T^N - 1$ is separable over~$k$, its divisor~$\mu_0$ is relatively prime to $S_0 \defeq (T^N - 1)/\mu_0$. Hence, as $A$ is Henselian, there are $R, S \in A[T]$ with $T^N - 1 = R S$, $R$ monic, $\smash{\overline{R}} = \mu_0$, and $\smash{\overline{S}} = S_0$. The reduction of $S(\phi)$ is $\smash{S_0(\overline{\phi})}$, which is invertible because $S_0$ and $\mu_0$ are relatively prime. Thus, the determinant of $S(\phi)$ is a unit, i.e.\ $S(\phi)$ is invertible, and $R(\phi) S(\phi) = \phi^N - \id = 0$ implies $R(\phi) = 0$. Conversely, let $F \in A[T]$ with $F(\phi) = 0$, and write $F = G R + H$ with $G, H \in A[T]$ and $\deg{H} < d \defeq \deg{R}$. Then $H(\phi) = 0$. The $A$-linear map
		\begin{align*}
			A^d \to \End_A(M), \quad (c_0, \dots, c_{d-1}) \mapsto \sum_{i = 0}^{d-1} c_i \phi^i
		\end{align*}
		has injective reduction since $\deg{\mu_0} = d$. So one of its $d \times d$ minors is a unit, the map is injective, and $H = 0$. Therefore, the kernel of $A[T] \to \End_A(M)$, $T \mapsto \phi$, is generated by the monic polynomial~$R$, i.e.\ $\minpoly{\phi} = R$, which divides $T^N - 1$ and has reduction $\mu_0$.

		Alternatively, the existence of $\minpoly{\phi}$ follows, even over every local ring, from \cite{MR3262438}: since $N$ is invertible, $\smash{E \defeq \inquo{A[T]}{\ideal{T^N - 1}}}$ is finite étale over~$A$, so $M$ is a projective $E$-module \cite[Proposition~2.9]{MR3262438} whose annihilator is generated by an idempotent \cite[Proposition~1.3]{MR3262438}. Hence, $\smash{A[\phi] \cong E/\annihilator{E}{M}}$ is a direct factor of~$E$ and thus projective over~$A$, so $\phi$ admits a minimal polynomial by \cite[Proposition~4.2]{MR3262438}.
	\end{rem}
	Consequently, the results of \autoref{sec:orthogonal-grps}, which use the assumption that $A$ is Artinian only through \autoref{thm:existence-minimal-polynomial}, hold for all Henselian local rings with finite residue field of odd characteristic. In contrast, \autoref{sec:fpf-isometries} uses that $A$ is finite; see \autoref{cor:fpf-reduction}.


	\section{Orthogonal Groups over Artinian Rings} \label{sec:orthogonal-grps}

	Now, we are in a position to determine the conjugacy classes of isometries in orthogonal groups over Artinian local rings with finite residue fields. We briefly summarise the main arguments.

	To classify the conjugacy classes of isometries in the orthogonal group, we will follow the general proof strategy outlined by Shinoda in \cite{MR591256}. It is based on the work of \cite{MR249519} and \cite{MR268192}; see also \cite{MR150210}.

	The reduction map for a local ring $A$ with regular quadratic module $(M,q)$,
	\begin{align*}
		\orthgrp{}{A}{M,q} \to \orthgrp{}{k}{\overline{M},\overline{q}}
	\end{align*}
	is surjective, which was shown by Knebusch; see \cite[Satz~(0.3)]{MR252382}.

	The study of conjugacy classes of isometries over a finite field naturally leads to their minimal polynomials. This is precisely the reason why we established the existence of minimal polynomials over Artinian local rings in \autoref{sec:minimal-polynomial}.

	Let $A$ be an Artinian local ring with maximal ideal~$\mathfrak{m}$ and residue field~$k$. Let $(M,q)$ be a regular quadratic module over $A$ and $\smash{\phi \in \orthgrp{}{A}{M,q}}$ an isometry of finite order~$N$ coprime to the characteristic of~$k$. Then, by \autoref{lem:minpoly-self-reciproical}, the minimal polynomial of $\phi$ (see \autoref{thm:existence-minimal-polynomial}) is self-reciprocal, i.e.\ for its degree~$d$,
	\begin{align*}
		\minpoly{\phi}^{\ast}(T) = \minpoly{\phi}(0)^{-1} T^d \minpoly{\phi}(T^{-1}) = \minpoly{\phi}(T).
	\end{align*}

	Over a Henselian ring, we can decompose the minimal polynomial into self-reciprocal factors that are either irreducible or the product of an irreducible polynomial and its reciprocal over the residue field; see \autoref{lem:decomp-module-ann}. We can decompose a regular quadratic module accordingly; see \autoref{lem:min-poly-orthgrp}.

	We will reduce the question of whether two isometries are conjugate to each other to the equivalence of Hermitian forms, and then apply \autoref{thm:equivalence-Hermitian-forms}.

	Using Galois theory for Henselian local rings with finite residue fields, we can determine the Witt classes of involution traces; see \autoref{lem:witt-cls-nontrivial}.

	This culminates in our main theorem, \autoref{thm:class-conj-classes}, where we obtain a classification of all conjugacy classes of isometries whose order is coprime to the characteristic of the residue field.

	\smallskip

	We emphasise that our proof is constructive; namely, we obtain canonical regular quadratic modules with an isometry with given invariants; see \autoref{rem:realisations-conj}. These quadratic modules can be built using involution traces; see \autoref{defi:involution-trace}.

	\smallskip

	We will use the existence of minimal polynomials proved in \autoref{thm:existence-minimal-polynomial} in the previous section. We will apply our result to classify conjugacy classes whose isometries are fixed-point free in \autoref{sec:fpf-isometries}.


	\subsection{Involution Traces and Self-Reciprocal Polynomials} \label{subsec:invol-tr}
	In this section, we recall the definition of a particular class of examples of modules with symmetric bilinear form from Galois theory, the so-called \textit{involution traces}; see \autoref{defi:involution-trace}. We will see that \textit{self-reciprocal} polynomials naturally yield involution traces; see \autoref{prop:involution-trace-self-reciprocal}. We will use involution traces in our proof of the main theorem, \autoref{thm:class-conj-classes}; see \autoref{subsubsec:type-one}. As these quadratic modules naturally admit fixed-point free isometries, they will also appear in \autoref{sec:fpf-isometries}.

	\medskip

	We start by recalling the definition of \textit{involution traces}. We will see in \autoref{ex:involution-trace-galois-fields} that separable field extensions naturally define involution traces. We will use a generalisation of this fact to local rings; see \autoref{prop:involution-trace-self-reciprocal}.
	\begin{defi}[Involution Trace, {\cite[Proposition~\RN{1}.7.2.4]{MR1096299}}] \label{defi:involution-trace}
		Let $A$ be a ring, and let $B$ be a commutative $A$-algebra that is projective and finitely generated as an $A$-module with an involution $\sigma \colon B \to B$ restricting to the identity on $A$. Then an $A$-linear map $\tau \colon B \to A$ is called an \textit{involution trace} if it satisfies
		\begin{enumerate}
			\item $\tau \circ \sigma = \tau$,
			\item the following symmetric $A$-bilinear form $\lambda_{\fieldextension{B}{A}} \colon B \times B \to A$ is regular:
			\begin{align*}
				\lambda_{\fieldextension{B}{A}}(x,y) \defeq \tau(x \sigma(y)).
			\end{align*}
		\end{enumerate}
	\end{defi}
	For a trivial involution $\sigma$, these are sometimes called \textit{Frobenius extensions}; see e.g.\ \cite[Section~\RN{1}.2.8]{MR491773}.

	The following example of an involution trace underlies the results in \cite{MR249519}.
	\begin{ex}[{\cite[\href{https://stacks.math.columbia.edu/tag/0BIL}{Tag~0BIL}]{stacks-project}}] \label{ex:involution-trace-galois-fields}
		Let $k$ be a field and $F$ a finite extension of~$k$. Then the trace map $\tr_{\fieldextension{F}{k}}$ is an involution trace for the trivial involution if and only if $F$ is separable over~$k$. If $F$ is a Galois extension of $k$ and $\sigma \colon F \to F$ is a nontrivial involution in the Galois group $\Galois{\fieldextension{F}{k}}$, then the trace map $\tr_{\fieldextension{F}{k}}$ is an involution trace for $\sigma$.
	\end{ex}

	\medskip

	We will see that separable polynomials that are self-reciprocal naturally yield involution traces. Let us quickly recall their definition.
	\begin{defi}[Reciprocal Polynomial] \label{defi:reciprocal}
		Let $A$ be a ring and $P \in A[T]$ a monic polynomial of degree $d$ with coefficients
		\begin{align*}
			P = \sum_{i = 0}^{d} a_i T^i = T^d + a_{d-1} T^{d-1} + \dots + a_1 T^1 + a_0
		\end{align*}
		such that $a_0 \in \units{A}$. Then the \textit{reciprocal polynomial} $\reciprocal{P}$ of $P$ is defined by
		\begin{align*}
			\reciprocal{P} & = a_0^{-1} T^d P(T^{-1}) = a_{0}^{-1} \sum_{i = 0}^d a_{d-i} T^i \\
			& = a_0^{-1} \big(a_0 T^d + a_1 T^{d-1} + \dots + a_{d-1} T^1 + 1 \big).
		\end{align*}
	\end{defi}
	A monic polynomial $P$ is called \textit{self-reciprocal} if it is equal to its reciprocal, i.e.\ if $P = \reciprocal{P}$ holds.

	Computing the reciprocal of a polynomial is straightforward.
	\begin{ex} \label{ex:reci-poly}
		For the ring $A = \modZ{7^3}$, the reciprocal polynomial of
		\begin{align*}
			Q = T^3 + 64 T^2 + 32 T + 342
		\end{align*}
		is given by
		\begin{align*}
			\reciprocal{Q} = T^3 + 311 T^2 + 279 T + 342.
		\end{align*}
		Note that over the residue field $k = \finitefield{7}$, the reduction $\overline{Q}$ of $Q$ is the irreducible polynomial
		\begin{align*}
			\overline{Q} = T^3 + T^2 + 4 T + 6
		\end{align*}
		with a reciprocal polynomial
		\begin{align*}
			\reciprocal{\overline{Q}} = T^3 + 3 T^2 + 6 T + 6,
		\end{align*}
		which is also irreducible. For example,
		\begin{align*}
			P = T^2 + 108 T + 1
		\end{align*}
		is a self-reciprocal polynomial over $A$.
	\end{ex}

	\medskip

	By definition, self-reciprocal polynomials canonically admit an involution. In fact, they define an involution trace:
	\begin{prop} \label{prop:involution-trace-self-reciprocal}
		Let $A$ be a local ring with maximal ideal~$\mathfrak{m}$ and finite residue field~$k$. Let $Q \in A[T]$ be a self-reciprocal
		separable polynomial whose reduction $\smash{\overline{Q} \in k[T]}$ over $k$ is irreducible. Then $\smash{E \defeq \inquo{A[T]}{\ideal{Q}}}$ is a commutative $A$-algebra, projective and finitely generated as an $A$-module with an involution $\sigma \colon E \to E$ given by $t \mapsto t^{-1}$, and the trace map $\tr_{\fieldextension{E}{A}}$ is an involution trace.
	\end{prop}
	\begin{proof}
		The trace map is $A$-linear, and it satisfies $\tr_{\fieldextension{E}{A}} \circ \sigma = \tr_{\fieldextension{E}{A}}$ because $\sigma$ is an $A$-algebra automorphism. Its reduction is the trace map of the finite field extension $\smash{\overline{E} = \inquo{k[T]}{\ideal{\overline{Q}}}}$ of~$k$, which is separable as $k$ is finite. Hence, by \autoref{lem:reduction-isomorphisms} and \autoref{ex:involution-trace-galois-fields}, the trace form $(x,y) \mapsto \tr_{\fieldextension{E}{A}}(x \sigma(y))$ is regular.
	\end{proof}
	In this result, the polynomial needs to be self-reciprocal, as otherwise $\sigma$ is not a well-defined involution on the algebra $E$.


	\subsection{Isometry of Hermitian Modules over Henselian Local Rings} \label{subsec:equiv-Hermitian-modules}
	In this section, we prove that two regular $1$-Hermitian $\sigma$-sesquilinear modules of the same rank are isometric if the reduction of the involution $\sigma$ is nontrivial. For finite fields, this property is well known. We will use this result in the proof of our main theorem, \autoref{thm:class-conj-classes}; see in particular \autoref{subsubsec:type-one}.

	\medskip

	As in \cite[Section~1.12]{MR591256}, a crucial property of finite fields necessary for the classification result is the equivalence of all regular Hermitian forms of the same rank over finite fields. We will prove the corresponding generalisation using the following result:
	\begin{prop}[{\cite[Proposition~\RN{1}.6.2.4]{MR1096299}}] \label{prop:orthogonal-basis-division-ring}
		Let $(M, \lambda)$ be an $\epsilon$-Hermitian space over a field $k$. Then $(M, \lambda)$ has an orthogonal basis in the following cases:
		\begin{enumerate}
			\item the involution of $k$ is nontrivial,
			\item the involution of $k$ is trivial, the form is symmetric ($\epsilon = 1$), and $\smash{\characteristic{k} \neq 2}$.
		\end{enumerate}
	\end{prop}

	\medskip

	Recall that over a local ring $A$, objects defined over $A$ (e.g.\ $A$-linear homomorphisms) often have the same properties as their reductions (change of rings) to the residue field~$k$ of $A$. For example, we will use the following lemma to prove in \autoref{lem:reduction-isomorphisms} that the regularity of a sesquilinear form over $A$ is equivalent to the regularity of its reduction over $k$.
	\begin{lem}[{\cite[Lemma~\RN{2}.4.5.3]{MR1096299}}] \label{lem:isos-preserved-by-change-rings}
		Let $A$ be a local ring with maximal ideal~$\mathfrak{m}$ and residue field~$k$, and let $M$ and $N$ be projective, finitely generated $A$-modules with an $A$-module homomorphism $\phi \colon M \to N$. Then $\phi$ is an isomorphism if and only if the reduction $\smash{\overline{\phi}}$ is an isomorphism.
	\end{lem}

	\smallskip

	Let $A$ be a local ring with maximal ideal~$\mathfrak{m}$ and residue field~$k$, and let $\sigma$ be an involution of $A$. Since $\sigma$ is a ring automorphism, it maps the unique maximal ideal to itself, $\sigma(\mathfrak{m}) = \mathfrak{m}$. Hence, $\sigma$ induces an involution $\overline{\sigma}$ of $k$, and the quotient map $A \to k$ is a map of rings with involution (cf.\ \cite[Section~\RN{2}.4.2.5]{MR1096299}). In particular, the reduction of a $\sigma$-sesquilinear module is a $\overline{\sigma}$-sesquilinear module.

	The following lemma can be seen as a generalisation of the fact that a quadratic module over a local ring is regular if and only if its reduction is a regular quadratic space.
	\begin{lem} \label{lem:reduction-isomorphisms}
		Let $A$ be a local ring with maximal ideal~$\mathfrak{m}$ and residue field~$k$. Let $\sigma$ be an involution of $A$. Let $E$ be a projective finitely generated $A$-module with a $\sigma$-sesquilinear form $s \colon E \times E \to A$. We denote the reduction modulo $\mathfrak{m}$ of $E$ by $\smash{\overline{E}}$ and the reduction of $s$ by $\smash{\overline{s} \colon \overline{E} \times \overline{E} \to k}$.

		Then the adjoint $\Adj{s} \colon E \to \rd{E}_{\sigma}$ is an isomorphism if and only if the adjoint $\smash{\Adj{\overline{s}} \colon \overline{E} \to \rd{\overline{E}}_{\overline{\sigma}}}$ is an isomorphism.
	\end{lem}
	\begin{proof}
		The quotient map $A \to k$ is a map of rings with involution; see above. Since a $\sigma$-semilinear map is determined by its values on a basis, $\rd{E}_{\sigma}$ is a free $A$-module with the basis dual to any basis of $E$, and its reduction is $\smash{\rd{\overline{E}}_{\overline{\sigma}}}$, so that $\Adj{\overline{s}}$ is the reduction of $\Adj{s}$. By \autoref{lem:isos-preserved-by-change-rings}, the $A$-linear map $\Adj{s}$ is an isomorphism if and only if $\Adj{\overline{s}}$ is an isomorphism.
	\end{proof}

	\medskip

	A main reason why Hermitian or quadratic modules over local rings can often be studied by reducing to the residue field is that orthogonal decompositions over the residue field can be lifted to orthogonal decompositions over the local ring.

	\begin{lem}[Lifting Orthogonal Decompositions over Local Rings] \label{lem:reduction}
		Let $A$ be a local ring with involution $\sigma$, maximal ideal~$\mathfrak{m}$, and residue field~$k$ of odd characteristic. Let $(M,\lambda)$ be a regular $\epsilon$-Hermitian $\sigma$-sesquilinear module. Then an orthogonal decomposition
		\begin{align*}
			\overline{M} = \overline{E} \orthsum \overline{F}
		\end{align*}
		over $k$, where $\overline{E}$ is a (free) regular subspace of $\overline{M}$, can be lifted to $M$, i.e.\ there exists an orthogonal decomposition $M = E \orthsum F$ of $M$ with regular $\epsilon$-Hermitian $\sigma$-sesquilinear modules $E$ and $F$ such that
		\begin{align*}
			\inquo{E}{\mathfrak{m}E} = \overline{E}, \quad \inquo{F}{\mathfrak{m}F} = \overline{F}.
		\end{align*}
	\end{lem}
	For the symmetric case, we refer to \cite[Lemma~1.1.4]{MR252382} or \cite[Corollary~\RN{1}.3.4]{MR491773}. The statement for Hermitian forms is a special case of \cite[Theorem~\RN{2}.4.6.1]{MR1096299}; note that since $2$ is invertible in $A$, every Hermitian module is even.

	\smallskip

	Using \autoref{lem:reduction} above and the fact that Hermitian spaces over fields usually admit an orthogonal basis (see \autoref{prop:orthogonal-basis-division-ring}), we directly obtain the following corollary. It will be instrumental in our proof of \autoref{thm:equivalence-Hermitian-forms}.
	\begin{cor} \label{cor:orthgonal-basis}
		Let $A$ be a local ring with involution $\sigma$ and residue field~$k$ of odd characteristic. Let $(M, \lambda)$ be a regular $\epsilon$-Hermitian $\sigma$-sesquilinear module. If either the reduced involution $\overline{\sigma}$ on $k$ is nontrivial or the reduced involution $\overline{\sigma}$ is trivial and $\epsilon = 1$, then $(M, \lambda)$ admits an orthogonal basis.
	\end{cor}
	\begin{proof}
		By \autoref{prop:orthogonal-basis-division-ring}, the reduction $\overline{M}$ over $k$ has an orthogonal basis. By the lifting result, \autoref{lem:reduction}, an orthogonal basis of an $\epsilon$-Hermitian module over $k$ can be lifted to an orthogonal basis over $A$.
	\end{proof}
	Note that the above corollary directly generalises \autoref{lem:orth-basis-quadr-module}.

	\medskip

	The proof idea for the following theorem is based on \cite[Section~9, Remark~3]{MR43785}.
	Let $A$ be a ring with a finite étale algebra $E$ over $A$. Note that by \cite[Proposition~2.9]{MR3262438}, a module over $E$ that is finitely generated and projective as a module over $A$ is also projective over $E$. This will be used in the proof of the following theorem.
	\begin{thm}[Isometry of Hermitian Forms over Henselian Rings] \label{thm:equivalence-Hermitian-forms}
		Let $A$ be a Henselian local ring with maximal ideal~$\mathfrak{m}$ and finite residue field~$k$ of odd characteristic. Let $E$ be a connected, finite étale $A$-algebra with involution $\sigma \in \Aut_{A}(E)$. Let $M$ be an $E$-module that is finitely generated and projective as an $A$-module and $s,t \colon M \times M \to E$ regular $1$-Hermitian $\sigma$-sesquilinear forms.

		If the reduction of the involution $\smash{\overline{\sigma} \in \Aut_{k}(\overline{E})}$ is nontrivial, then $s$ and $t$ are isometric.
	\end{thm}
	\begin{proof}
		Since $E$ is connected and finite over $A$, it is a Henselian local ring; see the proof of \autoref{prop:Galois-Henselian-rings}. By \cite[Proposition~2.9]{MR3262438}, $M$ is projective over~$E$ and hence free. Note that both $(M,s)$ and $(M,t)$ admit orthogonal bases by \autoref{cor:orthgonal-basis}, i.e.\ there exist bases $(e_i)$, $(e_i')$ such that
		\begin{align*}
			s(e_i, e_j) = \eta_i \delta_{i,j}, \quad t(e_i',e_j') = \xi_i \delta_{i,j}
		\end{align*}
		for some $\eta_i, \xi_i \in E$. Note that since $s$ and $t$ are $1$-Hermitian,
		\begin{align*}
			\eta_i = \sigma(\eta_i), \quad \xi_i = \sigma(\xi_i)
		\end{align*}
		for all $i$. Furthermore, let $E^{\sigma}$ be the fixed subring of $\sigma$. Since $\overline{\sigma}$ is nontrivial, there is an element $u \in E$ with $\overline{\sigma}(\overline{u}) \neq \overline{u}$. Then $\omega \defeq u - \sigma(u)$ satisfies $\sigma(\omega) = -\omega$ and $\overline{\omega} \neq 0$, so $\omega \in \units{E}$. As $2 \in \units{A}$, every $x \in E$ can be written uniquely as $x = a + b \omega$ with $a = \tfrac{1}{2}(x + \sigma(x)) \in E^{\sigma}$ and $b = \tfrac{1}{2}(x - \sigma(x)) \omega^{-1} \in E^{\sigma}$. Hence, $E = E^{\sigma} \oplus E^{\sigma} \omega \cong \inquo{E^{\sigma}[T]}{\ideal{T^2 - \omega^2}}$ with $\omega^2 \in \smash{\units{(E^{\sigma})}}$ is a quadratic separable $E^{\sigma}$-algebra, and $\sigma$ is its canonical involution. Moreover, $E^{\sigma}$ is a local ring (as $E$ is integral over $E^{\sigma}$) that is finitely generated as an $A$-module (as a direct summand of~$E$), hence Henselian by \cite[\href{https://stacks.math.columbia.edu/tag/04GH}{Tag~04GH}]{stacks-project}, and its residue field is contained in that of $E$, hence finite of odd characteristic. Thus, the norm space of $\fieldextension{E}{E^{\sigma}}$ is universal by \autoref{lem:norm-universal}.

		Since $s$ and $t$ are regular, $\eta_i, \xi_i \in \units{E}$, and hence $\eta_i, \xi_i \in \smash{\units{(E^{\sigma})}}$, for all $i$. Thus, by the universality of the norm space, there exist $\alpha_i, \beta_i \in E$ such that
		\begin{align*}
			\sigma(\alpha_i) \alpha_i = \eta_i^{-1}, \quad \sigma(\beta_i) \beta_i = \xi_i^{-1}.
		\end{align*}
		Thus, we obtain
		\begin{align*}
			s(\alpha_i e_i, \alpha_j e_j) = \delta_{i,j}, \quad t(\beta_i e_i', \beta_j e_j') = \delta_{i,j}.
		\end{align*}
		We conclude that $(M,s)$ and $(M,t)$ are isometric.
	\end{proof}
	If $\sigma$ is trivial, then $s$ and $t$ are symmetric bilinear forms over~$E$, and, as for finite fields (cf.\ \cite[Section~9, Remark~3]{MR43785}), they are isometric if and only if their determinants agree modulo $\smash{(\units{E})^2}$; this follows from (the proof of) \autoref{lem:witt-iso-Henselian} and \autoref{lem:square-red}.


	\subsection{Conjugacy Classes in Orthogonal Groups} \label{subsec:conj-cls}
	Now, we are in a position to prove the main theorem of this paper. We will determine the conjugacy classes of isometries of a regular quadratic module over an Artinian local ring with finite residue field whose orders are coprime to the characteristic of the residue field.

	We will see that conjugacy is related to minimal polynomials; see \autoref{lem:conj-implies-same-minpoly}. The minimal polynomials of isometries of a regular quadratic module are self-reciprocal polynomials; see \autoref{lem:minpoly-self-reciproical}. Again, using Hensel's lemma, we can further reduce to polynomials that are irreducible over the residue field; see \hyperref[subsubsec:type-one]{type~\RN{1}} and \hyperref[subsubsec:type-two]{type~\RN{2}} below.

	\medskip

	Let $A$ be a ring and $M$ a finitely generated projective $A$-module with two automorphisms $\phi, \psi \colon M \to M$. We will say that $\phi$ and $\psi$ are \textit{similar} if they are conjugate in $\Aut_{A}(M)$, i.e.\ if there is an automorphism $\tau \colon M \to M$ such that $\tau \circ \phi = \psi \circ \tau$. In fact, this already implies that the minimal polynomials are equal, $\minpoly{\phi} = \minpoly{\psi}$, by the following well-known fact.
	\begin{lem} \label{lem:conj-implies-same-minpoly}
		Let $A$ be a ring and $M$ a finitely generated projective $A$-module with two automorphisms $\phi, \psi \colon M \to M$. Suppose $\phi$ and $\psi$ admit minimal polynomials $\minpoly{\phi}$ and $\minpoly{\psi}$, respectively. Then, if $\phi$ and $\psi$ are similar, their minimal polynomials are equal.
	\end{lem}
	Of course, if two isometries of a regular quadratic module are conjugate in the orthogonal group of that regular quadratic module, they are also similar.

	\medskip

	Using Nakayama's lemma, we can show that automorphisms with the same minimal polynomial that is \textit{irreducible} over the residue field are similar.
	\begin{lem} \label{lem:irred-conjugated}
		Let $A$ be a local ring with maximal ideal~$\mathfrak{m}$ and finite residue field~$k$. Let $M$ be a projective finitely generated module over $A$, and let $\phi, \psi \colon M \to M$ be two $A$-linear automorphisms with monic separable minimal polynomial
		\begin{align*}
			Q \defeq \minpoly{\phi} = \minpoly{\psi} \in A[T]
		\end{align*}
		of degree $d$ whose reduction $\overline{Q} \in k[T]$ is irreducible. Then $\phi$ and $\psi$ are similar.
	\end{lem}
	\begin{proof}
		We will show that $M_{\phi}$ and $M_{\psi}$, i.e.\ $M$ viewed as an $A[T]$-module with $T$ acting by $\phi$ and by $\psi$, respectively, both are the direct sum of $m$ cyclic modules over~$A[T]$ with the same annihilator. Hence, there exist bases of the free module~$M$ such that $\phi$ and $\psi$ are represented by the same matrix. Therefore, they are conjugate; see \cite[Corollary~4.3.25]{MR1181420}. As before, we denote the reduction of the module $M$ modulo $\mathfrak{m}$ by $\overline{M}$. Let $E = \inquo{A[T]}{\ideal{Q}}$.
		By \cite[Proposition~\RN{7}.5.4]{MR643362}, there exist elements $\overline{x}_1, \dots, \overline{x}_m$ in $\overline{M}$ such that
		\begin{align*}
			\overline{M}_{\overline{\phi}} = k[T].\overline{x}_1 \oplus \dots \oplus
			k[T].\overline{x}_m
		\end{align*}
		as a $k[T]$-module with
		\begin{align*}
			k[T] \neq \annihilator{k{[T]}}{\overline{x}_1} & \supseteq \annihilator{k{[T]}}{\overline{x}_2} \supseteq \\
			& \dots \supseteq \annihilator{k{[T]}}{\overline{x}_m} = \annihilator{k{[T]}}{\overline{M}_{\overline{\phi}}} = (\overline{Q}).
		\end{align*}
		Note that the last equality holds since $\overline{Q}(\overline{\phi}) = 0$ and $\overline{Q}$ is irreducible. This already implies that
		\begin{align*}
			\annihilator{k{[T]}}{\overline{x}_1} = \annihilator{k{[T]}}{\overline{x}_2} = \dots = \annihilator{k{[T]}}{\overline{x}_m} = (\overline{Q}).
		\end{align*}
		Choose an $x_i \in M$ with reduction $\overline{x}_i$ modulo $\mathfrak{m}$. Then, the inclusion
		\begin{align*}
			\annihilator{A[T]}{x_i} \supseteq \ideal{Q}
		\end{align*}
		holds. By Nakayama's lemma, in the form of \autoref{lem:isos-preserved-by-change-rings}, the $E$-linear map $E^m \to M$, $(e_i)_i \mapsto \sum_i e_i x_i$, of finitely generated projective $A$-modules is an isomorphism, since its reduction $\smash{(\inquo{k[T]}{\ideal{\overline{Q}}})^m \to \overline{M}}$ is an isomorphism by the above. Thus, $\annihilator{A[T]}{x_i} = \ideal{Q}$, and $M$ is the direct sum of the $m$ free modules $E x_i \cong E$ of rank~$1$ over~$E$.
		The same argument applies to $\psi$.
	\end{proof}
	Note that above we have implicitly shown that there exist bases of $M$ such that $\phi$ and $\psi$ take the form $\smash{\companion{Q} \oplus \dots \oplus \companion{Q}}$ ($m$ copies) for the companion matrix $\companion{Q}$ of the monic polynomial $\smash{Q = \sum_{i = 0}^{d} a_i T^{i}}$,
	\begin{align*}
		\companion{Q} = \begin{bmatrix}
			0 & 0 & 0 & \dots & 0 & -a_0 \\
			1 & 0 & 0 & \dots & 0 & -a_1 \\
			0 & 1 & 0 & \dots & 0 & -a_2 \\
			& & & \vdots & & \\
			0 & 0 & 0 & \dots & 1 & -a_{d-1} \\
		\end{bmatrix} \in \matrixring{d}{A}.
	\end{align*}

	\medskip

	Let $A$ be an Artinian local ring with residue field~$k$ of characteristic~$p$, and let $(M,q)$ be a regular quadratic module over $A$, i.e.\ $M$ is a finitely generated projective $A$-module and $q \colon M \to A$ is a quadratic form such that the associated bilinear form $\lambda_q \colon M \times M \to A$ is regular. Since $A$ is local, $M$ is free, and we denote its rank by $n \defeq \rank{A}{M}$. Let $\phi \colon (M,q) \to (M,q)$ be an isometry of order~$N$ such that $p \nmid N$. By \autoref{thm:existence-minimal-polynomial}, $\phi$ has a separable minimal polynomial $\minpoly{\phi} \in A[T]$. By the next lemma, the minimal polynomial $\minpoly{\phi}$ is \textit{self-reciprocal}.
	\begin{lem} \label{lem:minpoly-self-reciproical}
		Let $A$ be a ring and $(M,q)$ a regular quadratic module over~$A$ with an isometry $\phi \colon (M,q) \to (M,q)$. If $\phi$ has a minimal polynomial $\minpoly{\phi} \in A[T]$, it is self-reciprocal, i.e.\
		\begin{align*}
			\minpoly{\phi}^{\ast}(T) = \minpoly{\phi}(0)^{-1} T^d \minpoly{\phi}(T^{-1}) = \minpoly{\phi}(T).
		\end{align*}
	\end{lem}
	\begin{proof}
		We denote the degree of the minimal polynomial $\minpoly{\phi}$ by $d$. For any polynomial $Q \in A[T]$, the bilinearity of the associated bilinear form $\lambda_q$ implies, since $\phi$ is an isometry, that
		\begin{align*}
			\lambda_q(Q(\phi)x,y) = \lambda_q(x,Q(\phi^{-1})y)
		\end{align*}
		for all $x,y \in M$. If we apply this equality to the minimal polynomial $\minpoly{\phi}$, we obtain that $\minpoly{\phi}(\phi^{-1})$ is the zero map since $\lambda_q$ is regular. Thus, $T^d \minpoly{\phi}(T^{-1})$ annihilates $\phi$, so that
		\begin{align*}
			\minpoly{\phi}(T) \mid T^d \minpoly{\phi}(T^{-1}).
		\end{align*}
		As $\minpoly{\phi}$ is monic of degree~$d$ and $T^d \minpoly{\phi}(T^{-1})$ has degree at most~$d$ with coefficient $a_0 \defeq \minpoly{\phi}(0)$ at $T^d$ and constant coefficient~$1$, we obtain $\smash{T^d \minpoly{\phi}(T^{-1}) = a_0 \minpoly{\phi}(T)}$ and, comparing constant coefficients, $a_0^2 = 1$. Hence, it follows that $a_0 \in \units{A}$ and $\minpoly{\phi}^{\ast} = a_0^{-1} T^d \minpoly{\phi}(T^{-1}) = \minpoly{\phi}$.
	\end{proof}

	\medskip

	We recall from \cite[Section~\RN{3}.4]{bourbaki-ac-1-4} that two elements $a,b$ of a ring $A$ are called \textit{strongly relatively prime} (or \textit{strongly coprime}) if the principal ideals $\ideal{a}$ and $\ideal{b}$ are relatively prime, i.e.\ if
	\begin{align*}
		\ideal{a} + \ideal{b} = A.
	\end{align*}
	In other words, if there exist elements $r,s \in A$ such that $ra + sb = 1$.

	\smallskip

	Strongly relatively prime polynomials admit the following well-known identity. Recall, however, that over a general ring, the polynomial ring is not necessarily a principal ideal domain.
	\begin{prop}[{\cite[Proposition~\RN{3}.4.1]{bourbaki-ac-1-4}}]
		Let $A$ be a ring and $P$ and $Q$ be two strongly relatively prime polynomials in $A[T]$ where $P$ is monic. Then every polynomial $X$ in $A[T]$ can be written \textit{uniquely} in the form
		\begin{align*}
			X = R P + S Q
		\end{align*}
		where $R, S \in A[T]$ and $\deg{S} < \deg{P}$.
		If further $\deg{X} \leq d$ for some $d \in \Z_{\geq 0}$ and $\deg{Q} \leq d - \deg{P}$, then $\deg{R} \leq d - \deg{P}$.
	\end{prop}

	\smallskip

	Using Nakayama's lemma, we can show that two polynomials over a local ring are strongly coprime if and only if their reductions are strongly coprime. Note that since, over a field $k$, the polynomial ring $k[T]$ is a principal ideal domain, being strongly relatively prime is equivalent to being relatively prime in $k[T]$.
	\begin{lem}[cf.\ {\cite[Proposition~\RN{3}.4.2]{bourbaki-ac-1-4}}] \label{lem:strongly-coprime-reduction}
		Let $A$ be a local ring with finite residue field~$k$. Let $P, Q \in A[T]$ be monic polynomials. Then, $P$ and $Q$ are strongly relatively prime in $A[T]$ if and only if their reductions $\overline{P}$ and $\overline{Q}$ in $k[T]$ are relatively prime.
	\end{lem}
	\begin{proof}
		Suppose that $P$ and $Q$ are strongly relatively prime in $A[T]$. Then there are polynomials $R$ and $S$ in $A[T]$ such that $1 = R P + S Q$. Thus,
		\begin{align*}
			1 = \overline{R} \overline{P} + \overline{S} \overline{Q}
		\end{align*}
		holds in $k[T]$.

		Conversely, suppose that $\overline{P}$ and $\overline{Q}$ are relatively prime, and hence strongly relatively prime, in the principal ideal domain $k[T]$, and let $X$ be a polynomial in $A[T]$. Then, by the above proposition applied over~$k$,
		\begin{align*}
			\overline{X} = \overline{R} \overline{P} + \overline{S} \overline{Q}
		\end{align*}
		for some polynomials $R, S \in A[T]$. Hence, $A[T] = \ideal{P,Q} + \mathfrak{m}A[T]$, where $\mathfrak{m}$ denotes the maximal ideal of~$A$. Now, $N \defeq \inquo{A[T]}{\ideal{P,Q}}$ is a finitely generated $A$-module since every polynomial $X$ is congruent modulo~$P$ to a polynomial of degree less than $\deg{P}$, and the above equality shows that
		\begin{align*}
			\mathfrak{m} N = N,
		\end{align*}
		so Nakayama's lemma implies that $N = 0$, i.e.\ $A[T] = \ideal{P,Q}$.
	\end{proof}

	\medskip

	For an $A$-module $M$ and an element $a \in A$, we define the \textit{$a$-annihilated} module by $M_a \defeq \{ m \in M : am = 0 \}$.

	\smallskip

	For products of strongly relatively prime elements, we can split annihilated modules into submodules annihilated by their factors.
	\begin{lem} \label{lem:factor-annihilator}
		Let $A$ be a ring and $M$ an $A$-module with $\annihilator{A}{M} = \ideal{a}$. Let $a = b c$ be a decomposition of $a$ into two strongly relatively prime elements $b, c \in A$. Then $M$ admits a direct-sum decomposition into the $b$- and $c$-annihilated submodules
		\begin{align*}
			M = M_b \oplus M_c.
		\end{align*}
	\end{lem}
	\begin{proof}
		There exist elements $r$ and $s$ in $A$ such that
		\begin{align*}
			1 = r b + s c.
		\end{align*}
		Thus, for every $m \in M$,
		\begin{align*}
			m = r b m + s c m = m_1 + m_2
		\end{align*}
		with $m_1 \defeq r b m$ and $m_2 \defeq s c m$. Note that
		\begin{align*}
			c m_1 = 0 = b m_2.
		\end{align*}
		In other words, $M_b + M_c = M$.

		Now, let $n \in M_b \cap M_c$, i.e.\ $b n = 0 = c n$. Then $\annihilator{A}{n} = A$, i.e.\ $n = 0$.
	\end{proof}

	\medskip

	Now, we apply the above to separable minimal polynomials over Henselian local rings with finite residue fields by first using Hensel's lemma to lift a factorisation into irreducible factors over the residue field to the Henselian local ring and then splitting the module by \autoref{lem:factor-annihilator}.

	\smallskip

	Note that a polynomial $P \in k[T]$ over a finite field $k$ with $P(0) \neq 0$ has a well-defined \textit{order}~$\ord{P}$, which is the least positive integer $e$ such that $P$ divides the polynomial $T^e - 1$; see \cite[Definition~3.2]{MR1429394}. Moreover, for $N \in \Z_{\geq 1}$, the polynomial $P$ divides $T^N - 1$ if and only if $\ord{P}$ divides~$N$; see \cite[Lemma~3.6]{MR1429394}. If $P$ is irreducible, then $\ord{P}$ divides $\abs{k}^{\deg{P}} - 1$ (see \cite[Corollary~3.4]{MR1429394}) and is thus coprime to the characteristic of~$k$.

	\smallskip

	By Hensel's lemma, this allows us to single out a lift of an irreducible polynomial over the residue field.
	\begin{defi}[Canonical Hensel Lift] \label{defi:canonical-lift}
		Let $A$ be a Henselian local ring with maximal ideal~$\mathfrak{m}$ and finite residue field~$k$, and let $P_0 \in k[T]$ be a monic irreducible polynomial with $P_0(0) \neq 0$. Then the \textit{canonical Hensel lift} of~$P_0$ is the unique monic polynomial $P \in A[T]$ with reduction $\overline{P} = P_0$ that divides $T^{\ord{P_0}} - 1$; see \autoref{lem:canonical-lift}.
	\end{defi}
	For example, the self-reciprocal polynomial $T^2 + 108 T + 1 \in \modZ{7^3}[T]$ from \autoref{ex:reci-poly} is the canonical Hensel lift of $T^2 + 3 T + 1 \in \finitefield{7}[T]$, which has order~$8$, but $T^2 + 3 T + 1$, viewed as a polynomial over $\modZ{7^3}$, does not divide $T^8 - 1$.
	\begin{lem} \label{lem:canonical-lift}
		Let $A$ be a Henselian local ring with maximal ideal~$\mathfrak{m}$ and finite residue field~$k$ of characteristic~$p$, and let $P_0 \in k[T]$ be a monic irreducible polynomial with $P_0(0) \neq 0$.
		\begin{enumerate}
			\item The canonical Hensel lift $P$ of~$P_0$ exists and is unique. Moreover, if $P_0$ divides $T^N - 1$ for some $N \in \Z_{\geq 1}$ coprime to~$p$, then $P$ is the unique monic divisor of $T^N - 1$ in $A[T]$ with reduction~$P_0$.
			\item Let $N \in \Z_{\geq 1}$ be coprime to~$p$. Then every monic divisor of $T^N - 1$ in $A[T]$ whose reduction is irreducible is the canonical Hensel lift of its reduction.
			\item The reciprocal $\reciprocal{P}$ is the canonical Hensel lift of~$\reciprocal{P_0}$. In particular, $P$ is self-reciprocal if and only if $P_0$ is.
		\end{enumerate}
	\end{lem}
	\begin{proof}
		Let $N \in \Z_{\geq 1}$ be coprime to~$p$, and let $R_0 \in k[T]$ be a monic divisor of $T^N - 1$. Since $T^N - 1$ is separable over~$k$, the polynomials $R_0$ and $(T^N - 1)/R_0$ are relatively prime. Hence, as $A$ is Henselian, there is a monic divisor $R \in A[T]$ of $T^N - 1$ with $\overline{R} = R_0$. If $R'$ is another such divisor, then $T^N - 1 = R S = R' S'$ for some $S, S' \in A[T]$ with $\smash{\overline{S} = \overline{S'} = (T^N - 1)/R_0}$, so $R = R'$ by the uniqueness of the lift; see below \autoref{defi:Henselian-ring}. Thus, there is a unique monic divisor of $T^N - 1$ with reduction~$R_0$.

		Now, $\ord{P_0}$ is coprime to~$p$, so the existence and uniqueness of $P$ follow by taking $R_0 = P_0$ and $N = \ord{P_0}$. If $P_0$ divides $T^N - 1$ for some $N$ coprime to~$p$, then $\ord{P_0}$ divides~$N$, so $P$ divides $T^{\ord{P_0}} - 1$ and hence $T^N - 1$. This proves~(1), and (2) follows from~(1), applied to $P_0 = \overline{Q}$ for a monic divisor $Q$ of $T^N - 1$ with irreducible reduction.

		For~(3), let $e \defeq \ord{P_0}$ and write $T^e - 1 = P S$ with a monic polynomial $S \in A[T]$. Then $P(0)$ and $S(0)$ are units, and by \autoref{defi:reciprocal}, $T^e - 1 = \reciprocal{(T^e - 1)} = \reciprocal{P} \reciprocal{S}$. Hence, $\reciprocal{P}$ is a monic divisor of $T^e - 1$ whose reduction $\reciprocal{P_0}$ is irreducible, so $\reciprocal{P}$ is the canonical Hensel lift of~$\reciprocal{P_0}$ by~(2). In particular, if $P_0 = \reciprocal{P_0}$, then $\reciprocal{P} = P$; the converse follows by reduction.
	\end{proof}

	\smallskip

	\begin{lem} \label{lem:decomp-module-ann}
		Let $A$ be a Henselian local ring with maximal ideal~$\mathfrak{m}$ and finite residue field~$k$ of characteristic~$p$. Let $M$ be a projective finitely generated module over $A$ and $\smash{\phi \colon M \to M}$ an $A$-linear automorphism of finite order~$N$ such that $p \nmid N$ with self-reciprocal minimal polynomial $\minpoly{\phi}$.

		Then there exists a unique factorisation, up to relabelling the polynomials $P_i$ and $Q_j$ and interchanging $Q_j$ and $\reciprocal{Q}_j$,
		\begin{align*}
			\minpoly{\phi} = P_1 \dots P_m Q_1 \reciprocal{Q}_1 \dots Q_n \reciprocal{Q}_n \in A[T]
		\end{align*}
		where $P_1, \dots, P_m, Q_1, \reciprocal{Q}_1, \dots, Q_n, \reciprocal{Q}_n$ are pairwise strongly coprime and
		\begin{enumerate}
			\item $\{P_i\}$ are polynomials over $A$ such that
			\begin{enumerate}
				\item $P_i$ is monic and self-reciprocal,
				\item the reduction $\overline{P}_i \in k[T]$ is irreducible and $P_i$ its canonical Hensel lift;
			\end{enumerate}
			\item $\{Q_j\}$ are polynomials over $A$ such that
			\begin{enumerate}
				\item $Q_j$ is monic and non-self-reciprocal,
				\item the reduction $\overline{Q}_j \in k[T]$ is irreducible and $Q_j$ its canonical Hensel lift.
			\end{enumerate}
		\end{enumerate}

		Then, the module $M$ admits the decomposition
		\begin{align*}
			M = \cusbigoplus_{i = 1}^{m} M_{P_i} \oplus \cusbigoplus_{j = 1}^{n} (M_{Q_j} \oplus M_{\reciprocal{Q_j}})
		\end{align*}
		into free modules such that the restrictions of $\phi$,
		\begin{align*}
			\restr{\phi}{M_{P_i}}, \quad \restr{\phi}{M_{Q_j}}, \quad \restr{\phi}{M_{\reciprocal{Q_j}}}
		\end{align*}
		have minimal polynomials $P_i$, $Q_j$, and $\reciprocal{Q_j}$, respectively.
	\end{lem}
	\begin{proof}
		Since $\minpoly{\phi}$ is self-reciprocal and divides $T^N - 1$, the reduction $\overline{\minpoly{\phi}}$ of $\minpoly{\phi}$ to $k[T]$ is self-reciprocal and separable (as $p \nmid N$). As taking reciprocals permutes its irreducible factors, it can be decomposed into strongly relatively prime irreducible factors
		\begin{align*}
			\overline{\minpoly{\phi}} = \overline{P}_1 \dots \overline{P}_m \overline{Q}_1 \reciprocal{\overline{Q}}_1 \dots \overline{Q}_n \reciprocal{\overline{Q}}_n
		\end{align*}
		with $\overline{P}_i$ self-reciprocal and $\overline{Q}_j \neq \reciprocal{\overline{Q}}_j$, whose orders divide $N$ and are thus coprime to~$p$. As $A$ is Henselian and $\minpoly{\phi}$ is separable, we can lift this decomposition to~$A$; the lift is unique (see \cite[Remark~\RN{1}.4.1]{MR559531}), and the lifted factors are pairwise strongly coprime by \autoref{lem:strongly-coprime-reduction}. As they divide $T^N - 1$, they are the canonical Hensel lifts of their reductions by \autoref{lem:canonical-lift}. Since $\minpoly{\phi}$ is self-reciprocal and $\reciprocal{(RS)} = \reciprocal{R} \reciprocal{S}$, every factorisation $\minpoly{\phi} = R S$ in this lift yields a factorisation $\minpoly{\phi} = \reciprocal{R} \reciprocal{S}$ lifting $\smash{\reciprocal{\overline{R}} \reciprocal{\overline{S}}}$. By uniqueness, the lifts $P_i$ of $\smash{\overline{P}_i}$ are self-reciprocal, and the lift of $\smash{\reciprocal{\overline{Q}}_j}$ is $\reciprocal{Q_j}$.

		Because $\smash{\annihilator{A[T]}{M} = \ideal{\minpoly{\phi}}}$, the module $M$ admits, by \autoref{lem:factor-annihilator}, the decomposition
		\begin{align*}
			M = \cusbigoplus_{i = 1}^{m} M_{P_i} \oplus \cusbigoplus_{j = 1}^{n} (M_{Q_j} \oplus M_{\reciprocal{Q_j}}).
		\end{align*}
		Note that, since the above decomposition is direct, and direct summands of free modules are projective, the direct summands are also free (since being projective implies being free for local rings). Finally, for each factor $\smash{R \in \{P_i, Q_j, \reciprocal{Q_j}\}}$, the summand $M_R$ is nonzero, as otherwise $\smash{\minpoly{\phi}/R}$ would annihilate~$M$. Hence, the minimal polynomial of $\smash{\restr{\phi}{M_R}}$ (see \autoref{rem:minpoly-Henselian}), a monic divisor of~$R$ whose reduction is a nonconstant divisor of the irreducible polynomial~$\smash{\overline{R}}$, is equal to~$R$.
	\end{proof}

	\medskip

	So far, we have only worked with projective, finitely generated modules and automorphisms between them. Now, we apply our result to regular quadratic modules. The decomposition in \autoref{lem:decomp-module-ann}
	can be refined to an orthogonal decomposition.

	\begin{lem}[{cf.\ \cite[Lemma~3.1]{MR249519}}] \label{lem:orth-decom-minpoly}
		Let $A$ be a Henselian local ring with maximal ideal~$\mathfrak{m}$ and finite residue field~$k$. Let $(M,q)$ be a regular quadratic module over~$A$ and $\phi \colon (M,q) \to (M,q)$ an isometry of finite order~$N$ such that $\characteristic{k} \nmid N$ with minimal polynomial $\minpoly{\phi}$.
		Let
		\begin{align*}
			\minpoly{\phi} = P_1 \dots P_m Q_1 \reciprocal{Q}_1 \dots Q_n \reciprocal{Q}_n \in A[T]
		\end{align*}
		be the factorisation from \autoref{lem:decomp-module-ann}, which applies by \autoref{lem:minpoly-self-reciproical}.

		Then, the decomposition from \autoref{lem:decomp-module-ann}
		\begin{align*}
			M = \cusbigoplus_{i = 1}^{m} M_{P_i} \oplus \cusbigoplus_{j = 1}^{n} (M_{Q_j} \oplus M_{\reciprocal{Q_j}})
		\end{align*}
		yields the orthogonal decomposition into regular quadratic modules
		\begin{align*}
			M = \bigorthsum_{i = 1}^{m} M_i^{(\mathrm{\RN{1}})} \orthsum \bigorthsum_{j = 1}^{n} M_j^{(\mathrm{\RN{2}})}
		\end{align*}
		where $\smash{M_i^{(\mathrm{\RN{1}})} \defeq M_{P_i}}$ and $\smash{M_j^{(\mathrm{\RN{2}})} \defeq M_{Q_j} \oplus M_{\reciprocal{Q_j}}}$.
	\end{lem}
	\begin{proof}
		We view $M$ as an $A[T]$-module where $T$ acts as $\phi$ on $M$. For all $x,y \in M$,
		\begin{align*}
			\lambda_q(x,T y) = \lambda_q(x, \phi(y)) = \lambda_q(\phi^{-1}x,y) = \lambda_q(T^{-1}x,y).
		\end{align*}

		Thus, for $y \in M_{P_i}$ for $P_i \in A[T]$ as in the factorisation from \autoref{lem:decomp-module-ann},
		\begin{align*}
			0 = \lambda_q(x, P_i(T)y) = \lambda_q(P_i(T^{-1})x,y),
		\end{align*}
		so $\smash{M_{P_i}}$ is orthogonal to $P_i(T^{-1})M$. For every polynomial $R \in A[T]$ strongly relatively prime to $P_i$, the map
		\begin{align*}
			M_R \to M_R, \quad x \mapsto T^{-d} P_i(T) x
		\end{align*}
		with $d = \deg{P_i}$ defines an isomorphism, and $P_i(T^{-1}) = P_i(0) T^{-d} P_i(T)$ with $P_i(0) \in \units{A}$ since $P_i$ is self-reciprocal. Hence, $M_R \subseteq P_i(T^{-1}) M$, so $\smash{M_{P_i}}$ is orthogonal to $M_R$.

		Similarly, $\smash{M_{Q_j}}$ is totally isotropic for all $Q_j \in A[T]$ in the factorisation above. We can orthogonally split off $\smash{M_{Q_j} \oplus M_{\reciprocal{Q_j}}}$.
	\end{proof}

	\smallskip

	Let $A$ be a Henselian local ring with maximal ideal~$\mathfrak{m}$ and finite residue field~$k$ of odd characteristic and $(M,q)$ a regular quadratic module over $A$. Further, let $\smash{\phi \colon (M,q) \to (M,q)}$ be an isometry of order coprime to the characteristic $\characteristic{k}$ with minimal polynomial $\minpoly{\phi}$. With the lemma above, we have reduced the problem of determining conjugacy classes to the following two cases.
	\begin{description} \label{two-types}
		\item[Type \RN{1}] $\minpoly{\phi} = P \in A[T]$ is self-reciprocal and the reduction
		$\overline{P}$ irreducible in $k[T]$.
		\item[Type \RN{2}] $\minpoly{\phi} = Q \reciprocal{Q} \in A[T]$ with $Q \neq \reciprocal{Q}$ and the reduction $\overline{Q}$ irreducible in $k[T]$.
	\end{description}

	We now study type~\RN{1} and type~\RN{2} separately.


	\subsubsection*{Type I: Isometries with Self-Reciprocal Minimal Polynomial Irreducible over the Residue Field} \label{subsubsec:type-one}
	In the following, we describe isometries with self-reciprocal minimal polynomial that is irreducible over the residue field. We will use the fact that they naturally define involution traces and the isometry of Hermitian modules for nontrivial involutions; see \autoref{prop:involution-trace-self-reciprocal} and \autoref{thm:equivalence-Hermitian-forms}.

	\medskip

	In the definition of the reciprocal polynomial for a monic polynomial, \autoref{defi:reciprocal}, we divided by the constant coefficient of the polynomial to ensure that the reciprocal of a monic polynomial is monic. The following remark explains why, in this paper, this factor is equal to $1$ for self-reciprocal polynomials of degree at least~$2$.
	\begin{rem} \label{rem:self-reciprocal-constant}
		Let $k$ be a finite field of odd characteristic. Let $P \in k[T]$ be a monic, irreducible, self-reciprocal polynomial over $k$. If $\deg{P} > 1$, then the degree of $P$ is even and $P(0) = 1$; see e.g.\ \cite[Lemma~1.7.3]{MR591256}. The same holds for a self-reciprocal polynomial $P$ over a local ring~$A$ with residue field~$k$ whose reduction is irreducible of degree greater than $1$: comparing the constant coefficients in $T^{\deg{P}} P(T^{-1}) = P(0) P(T)$ yields $P(0)^2 = 1$, and $P(0) \equiv 1$ modulo the maximal ideal, so $P(0) = 1$ as $2 \in \units{A}$.
	\end{rem}

	\medskip

	Let $A$ be a local ring with finite residue field~$k$ of odd characteristic and $(M,q)$ a regular quadratic module over~$A$. We consider an isometry $\phi \colon (M,q) \to (M,q)$ with a self-reciprocal minimal polynomial $P \in A[T]$ such that its reduction $\overline{P} \in k[T]$ is irreducible. The subalgebra of $A$-module endomorphisms of $M$ generated by $\phi$ is isomorphic to the algebra $E_{\phi} = E \defeq \inquo{A[T]}{\ideal{P}}$. We denote the image of $T$ in $E$ by~$t$. Since the polynomial $P$ is self-reciprocal, we obtain the involution
	\begin{align*}
		\sigma \colon E \to E, \quad t \mapsto t^{-1}.
	\end{align*}
	Furthermore, the algebra $E$ is separable and free over $A$ with basis $(1,t,\dots,t^{d-1})$ for $d \defeq \deg{P}$.
	By \autoref{prop:involution-trace-self-reciprocal}, the étale algebra $E$ admits an involution trace
	\begin{align*}
		\tr_{\fieldextension{E}{A}} \colon E \to A
	\end{align*}
	for the involution $\sigma$.

	We denote $M$ considered as an $E_{\phi}$-module by $M_{\phi}$, i.e.\ an element $\sum_{i = 0}^{d-1} a_i t^i \in E_{\phi}$ acts on an element $x \in M$ by
	\begin{align*}
		\sum_{i = 0}^{d-1} a_i \phi^i(x).
	\end{align*}

	\medskip

	In the following three lemmas, we show that in type~\RN{1}, conjugacy classes in orthogonal groups naturally correspond to isometry classes of certain Hermitian forms, introduced in the lemma below.
	\begin{lem}[{cf.\ \cite[Lemma~1.1]{MR249519}}] \label{lem:ex-hermitian-from}
		Let $A$ be a local ring with finite residue field~$k$ of odd characteristic and $(M,q)$ a regular quadratic module over~$A$. Further, let $\phi \colon (M,q) \to (M,q)$ be an isometry with self-reciprocal minimal polynomial $P \in A[T]$ such that its reduction $\overline{P} \in k[T]$ is irreducible. Let $E$ be the étale algebra $\inquo{A[T]}{\ideal{P}}$ over $A$ with involution $\sigma \colon E \to E,\ t \mapsto t^{-1}$ and involution trace $\tr_{\fieldextension{E}{A}} \colon E \to A$.

		There exists a unique regular $1$-Hermitian $\sigma$-sesquilinear form
		\begin{align*}
			s = s_{\phi} \colon M_{\phi} \times M_{\phi} \to E
		\end{align*}
		such that, for all $u \in E$ and all $x,y \in M$,
		\begin{align*}
			\lambda_q(u x, y) = \tr_{\fieldextension{E}{A}}(u s(x,y)).
		\end{align*}
	\end{lem}
	\begin{proof}
		For all $x, y \in M$, the map
		\begin{align*}
			E \to A, \quad u \mapsto \lambda_q(ux, y),
		\end{align*}
		is $A$-linear. Since $\tr_{\fieldextension{E}{A}}$ is an involution trace for $\sigma$ by \autoref{prop:involution-trace-self-reciprocal}, i.e.\ the trace form $\lambda_{\fieldextension{E}{A}}$ is regular, there exists a unique element $s(x,y) \in E$ such that
		\begin{align*}
			\tr_{\fieldextension{E}{A}}(u s(x,y)) = \lambda_q(ux, y)
		\end{align*}
		for all $u \in E$. Now, the map $s$ defined by
		\begin{align*}
			s \colon M \times M \to E, \quad (x,y) \mapsto s(x,y)
		\end{align*}
		is a regular $1$-Hermitian form for the involution $\sigma$:
		\begin{enumerate}
			\item The biadditivity follows from the $A$-linearity of $\tr_{\fieldextension{E}{A}}$ and the biadditivity of $\lambda_q$.
			\item The equality
			\begin{align*}
				\tr_{\fieldextension{E}{A}}(u s(v x,y)) = \lambda_q(u v x, y) = \tr_{\fieldextension{E}{A}}(u v s(x,y))
			\end{align*}
			for all $u \in E$ implies $s(v x,y ) = v s(x,y)$ for all $v \in E$ and $x,y \in M$.
			\item Similarly, the equality
			{\allowdisplaybreaks
			\begin{align*}
				\tr_{\fieldextension{E}{A}}(u s(x,y)) & = \lambda_q(ux, y) \\
				& = \lambda_q(y, ux) \\
				& = \lambda_q(\sigma(u)y, x) \\
				& = \tr_{\fieldextension{E}{A}}(\sigma(u) s(y,x)) \\
				& = \tr_{\fieldextension{E}{A}}(u \sigma(s(y,x)))
			\end{align*}}%
			for all $u \in E$ implies $s(x,y) = \sigma(s(y,x))$ for all $x,y \in M$. In particular, $s(x,vy) = \sigma(v) s(x,y)$ for all $v \in E$.
			\item The form $s$ is regular: if $s(x,\ph) = 0$, then $\lambda_q(x,\ph) = \tr_{\fieldextension{E}{A}}(s(x,\ph)) = 0$, so $x = 0$. Conversely, let $f \colon M \to E$ be $\sigma$-semilinear. Since $\lambda_q$ is regular, there is an $x \in M$ with $\lambda_q(x,y) = \tr_{\fieldextension{E}{A}}(f(y))$ for all $y \in M$. Then
			\begin{align*}
				\tr_{\fieldextension{E}{A}}(\sigma(v) s(x,y)) = \lambda_q(x, vy) = \tr_{\fieldextension{E}{A}}(\sigma(v) f(y))
			\end{align*}
			for all $v \in E$ and $y \in M$, so $s(x,\ph) = f$ since the trace form $\lambda_{\fieldextension{E}{A}}$ is regular. \qedhere
		\end{enumerate}
	\end{proof}

	\smallskip

	Now, we show how certain conjugacy classes correspond to isometry classes of the Hermitian forms introduced in \autoref{lem:ex-hermitian-from}.
	\begin{lem} \label{lem:red-isom-hermitian-form}
		Let $A$ be a local ring with finite residue field~$k$ of odd characteristic and $(M,q)$ a regular quadratic module over~$A$, and let $\phi, \psi \colon (M,q) \to (M,q)$ be two isometries with the same self-reciprocal minimal polynomial $P \in A[T]$ such that its reduction $\overline{P} \in k[T]$ is irreducible. Let $E$ be the étale algebra $\inquo{A[T]}{\ideal{P}}$ over $A$.
		Finally, let $(M_{\phi}, s_{\phi})$ be the regular $1$-Hermitian $\sigma$-sesquilinear module from \autoref{lem:ex-hermitian-from} with $\lambda_q(u x, y)= \tr_{\fieldextension{E}{A}}(u s(x,y))$ for all $u \in E$, $x,y \in M$, and $(M_{\psi}, s_{\psi})$ the corresponding module for~$\psi$.

		There is a bijection between the isometries $\smash{\tau \colon (M,q) \to (M,q)}$ with $\smash{\psi \circ \tau = \tau \circ \phi}$ and the set of isometries $\smash{f \colon (M_{\phi},s_{\phi}) \to (M_{\psi},s_{\psi})}$, i.e.\ the $E$-linear isomorphisms $\smash{f \colon M_{\phi} \to M_{\psi}}$ such that $\smash{s_{\psi}(f(x), f(y)) = s_{\phi}(x,y)}$.
	\end{lem}
	\begin{proof}
		Let $\tau \colon (M,q) \to (M,q)$ be an isometry such that $\psi \circ \tau = \tau \circ \phi$. Then $\tau$ is an $E$-linear map from $M_{\phi}$ to $M_{\psi}$ since for all elements $\smash{u = \sum_{i = 0}^{d-1} a_i t^i \in E_{\phi}}$,
		\begin{align*}
			\tau(u x) = \tau \big( \sum_{i = 0}^{d-1} a_i \phi^i x \big) = \sum_{i = 0}^{d-1} a_i \psi^i (\tau x).
		\end{align*}
		Then, for all $\smash{u = \sum_{i = 0}^{d-1} a_i t^i \in E}$ and $x,y \in M$,
		{\allowdisplaybreaks
		\begin{align*}
			\tr_{\fieldextension{E}{A}}(u s_{\psi}(\tau x, \tau y)) & = \lambda_q(\sum_{i = 0}^{d-1} a_i \psi^i \tau x, \tau y) \\
			& = \lambda_q(\tau \sum_{i = 0}^{d-1} a_i \phi^i x, \tau y) \\
			& = \lambda_q(\sum_{i = 0}^{d-1} a_i \phi^i x, y) \\
			& = \tr_{\fieldextension{E}{A}}(u s_{\phi}(x, y)).
		\end{align*}}%

		Conversely, given an $E$-linear isomorphism $f \colon M_{\phi} \to M_{\psi}$ such that
		\begin{align*}
			s_{\psi}(f(x), f(y)) = s_{\phi}(x,y)
		\end{align*}
		for all $x,y \in M_{\phi}$, we obtain for all $x \in M$,
		\begin{align*}
			f(\phi x) = f(tx) = t f(x) = \psi f(x),
		\end{align*}
		so $f \circ \phi = \psi \circ f$. Furthermore, $f$ is an isometry since
		\begin{align*}
			\lambda_q(fx, fy) & = \tr_{\fieldextension{E}{A}}(s_{\psi}(fx, fy)) \\
			& = \tr_{\fieldextension{E}{A}}(s_{\phi}(x, y)) \\
			& = \lambda_q(x,y). \qedhere
		\end{align*}
	\end{proof}

	\smallskip

	Finally, we can conclude using \autoref{thm:equivalence-Hermitian-forms}:
	\begin{lem} \label{lem:selfreci-irred-min-poly-conj}
		Let $A$ be a Henselian local ring with maximal ideal~$\mathfrak{m}$ and finite residue field~$k$ of odd characteristic. Let $(M,q)$ be a regular quadratic module over~$A$ and $\phi, \psi \colon (M,q) \to (M,q)$ two isometries with minimal polynomial
		\begin{align*}
			P = \minpoly{\phi} = \minpoly{\psi}
		\end{align*}
		for a separable self-reciprocal polynomial $P \in A[T]$ such that $\smash{\overline{P}}$ is irreducible in~$k[T]$. Then $\phi$ and $\psi$ are conjugate, i.e.\ there exists an isometry $\tau \colon (M,q) \to (M,q)$ such that $\tau \circ \phi = \psi \circ \tau$.
	\end{lem}
	\begin{proof}
		If $\deg{P} = 1$, then $P = T \mp 1$ as $P$ is self-reciprocal, so $\phi = \psi = \pm \id$, and there is nothing to prove. Otherwise, we denote $E \defeq \inquo{A[T]}{(P)}$ with involution $\sigma \colon E \to E,\ t \mapsto t^{-1}$. Since $\overline{P}$ is irreducible and $P$ is separable, $E$ is a connected, finite étale $A$-algebra, which is local and Henselian. Moreover, the reduction $\overline{\sigma}$ is nontrivial, since $\smash{\overline{t}^{2} \neq 1}$ as $\overline{P}$ is irreducible of degree at least~$2$.

		There exists an $A$-linear automorphism $\alpha \colon M \to M$ such that $\alpha \circ \phi = \psi \circ \alpha$ by \autoref{lem:irred-conjugated}. Thus, we get an $A[T]$-linear isomorphism $\alpha \colon M_{\phi} \to M_{\psi}$. The pullback $(M_{\phi}, \pub{\alpha} s_{\psi})$ is a regular $1$-Hermitian $\sigma$-sesquilinear module such that
		\begin{align*}
			\alpha \colon (M_{\phi}, \pub{\alpha} s_{\psi}) \to (M_{\psi}, s_{\psi})
		\end{align*}
		defines an isometry.

		Now, $(M_{\phi}, s_{\phi})$ and $(M_{\phi}, \pub{\alpha}s_{\psi})$ are two regular $1$-Hermitian $\sigma$-sesquilinear modules and thus isometric by \autoref{thm:equivalence-Hermitian-forms}, i.e.\ there exists an isometry
		\begin{align*}
			f \colon (M_{\phi}, s_{\phi}) \to (M_{\phi}, \pub{\alpha}s_{\psi}).
		\end{align*}
		Thus, we obtain the isometry $\alpha \circ f \colon (M_{\phi}, s_{\phi}) \to (M_{\psi}, s_{\psi})$, so $\phi$ and $\psi$ are conjugate as isometries of $(M,q)$ by \autoref{lem:red-isom-hermitian-form}.
	\end{proof}


	\subsubsection*{Type II: Isometries with Minimal Polynomial Not Irreducible over the Residue Field} \label{subsubsec:type-two}
	In the following, we describe isometries with minimal polynomial whose reduction over the residue field is the product of a non-self-reciprocal irreducible polynomial and its reciprocal over the residue field. We will see that the quadratic modules with such an isometry are hyperbolic (i.e.\ have trivial Witt class) and that the conjugacy class of the corresponding isometry is determined by (one of the) irreducible factors of the reduction of the minimal polynomial.

	\medskip

	Similarly to \autoref{lem:selfreci-irred-min-poly-conj}, we also obtain in type~\RN{2} the following lemma:
	\begin{lem} \label{lem:nonreci-irred-min-poly-conj}
		Let $A$ be a Henselian local ring with maximal ideal~$\mathfrak{m}$ and finite residue field~$k$ of odd characteristic. Let $(M,q)$ be a regular quadratic module over~$A$ and $\phi, \psi \colon (M,q) \to (M,q)$ two isometries with minimal polynomial
		\begin{align*}
			Q \reciprocal{Q} = \minpoly{\phi} = \minpoly{\psi}
		\end{align*}
		for a separable $Q \in A[T]$ such that $\overline{Q}$ is irreducible and not self-reciprocal in $k[T]$. Then $\phi$ and $\psi$ are conjugate, i.e.\ there exists an isometry $\tau \colon (M,q) \to (M,q)$ such that $\tau \circ \phi = \psi \circ \tau$.
	\end{lem}
	\begin{proof}
		First note that since $\overline{Q}$ is irreducible and not self-reciprocal, $Q$ and $\reciprocal{Q}$ are strongly relatively prime in $A[T]$ by \autoref{lem:strongly-coprime-reduction}. Hence, as an $A[T]$-module, $M$ splits as
		\begin{align*}
			M = M_{Q} \oplus M_{\reciprocal{Q}}.
		\end{align*}
		Now, $\smash{M_{Q}}$ and $\smash{M_{\reciprocal{Q}}}$ are $\phi$-invariant. The restriction $\phi'$ of $\phi$ to $\smash{M_Q}$ has minimal polynomial $Q$ since $\smash{\minpoly{\phi'} \mid Q}$, $M_Q \neq 0$, and $\smash{\overline{Q}}$ is irreducible. Since $M$ has the totally isotropic submodule $M_Q$ (see the proof of \autoref{lem:orth-decom-minpoly}) with $\smash{\orth{(M_Q)} = M_Q}$, $M$ is hyperbolic by \cite[Theorem~\RN{1}.4.6]{MR491773}. Moreover, as $M_{\reciprocal{Q}}$ is also totally isotropic, $\lambda_q$ restricts to a perfect pairing $M_Q \times M_{\reciprocal{Q}} \to A$, so we can find bases $(x_i)_{i = 1}^{m}$ of $M_Q$ and $(y_i)_{i = 1}^{m}$ of $M_{\reciprocal{Q}}$ such that
		\begin{align*}
			\lambda_q(x_i,y_j) = \delta_{i,j}
		\end{align*}
		for all $i,j$. In these bases, $\phi$ (and similarly $\psi$, in bases $(x_i')$ of $M_{Q,\psi}$ and $(y_i')$ of $M_{\reciprocal{Q},\psi}$ chosen in the same way) can be written as
		\begin{align*}
			\phi(x_i) = \sum_{j = 1}^{m} X^{\phi}_{i,j} x_j, \quad \phi^{-1}(x_i) = \sum_{j = 1}^{m} Y^{\phi}_{i,j} x_j,
		\end{align*}
		with a matrix $X^{\phi} \in \matrixring{m}{A}$ and inverse $Y^{\phi} = (X^{\phi})^{-1}$.
		Thus, it follows from
		\begin{align*}
			\lambda_q(x_i, \phi y_j) = \lambda_q(\phi^{-1} x_i, y_j) = Y^{\phi}_{i,j}
		\end{align*}
		that
		\begin{align*}
			\phi(y_i) = \sum_{j = 1}^{m} Y^{\phi}_{j,i} y_j.
		\end{align*}
		In the bases $(x_1,\dots,x_m,y_1,\dots,y_m)$ and $(x_1',\dots,x_m',y_1',\dots,y_m')$ of $M$, $\phi$ and $\psi$ thus take the forms
		\begin{align*}
			\begin{bmatrix}
				X^{\phi} & 0 \\
				0 & \transpose{(Y^{\phi})}
			\end{bmatrix}\!,
			\quad
			\begin{bmatrix}
				X^{\psi} & 0 \\
				0 & \transpose{(Y^{\psi})}
			\end{bmatrix}\!,
		\end{align*}
		respectively.
		Since the restrictions $\phi'$ and $\psi'$ of $\phi$ and $\psi$ to $M_Q$ and $M_{Q,\psi}$, respectively, have the same minimal polynomial $Q$, irreducible over the residue field, the matrices $X^{\phi}$ and $X^{\psi}$ are similar by \autoref{lem:irred-conjugated}. Hence, there exists an invertible $R \in \matrixring{m}{A}$ with $X^{\phi} R = R X^{\psi}$. Then,
		\begin{align*}
			\begin{bmatrix}
				R & 0 \\
				0 & (\transpose{R})^{-1}
			\end{bmatrix}
		\end{align*}
		defines, with respect to these two bases, an isometry $\tau \colon (M,q) \to (M,q)$ such that $\tau \circ \phi = \psi \circ \tau$.
	\end{proof}
	We also record for later use that if $\phi$ is an isometry with minimal polynomial of the above form, then the corresponding Witt class of $(M,q)$ is trivial because it is hyperbolic.


	\subsubsection*{Combining Type I and Type II}

	Combining \autoref{lem:selfreci-irred-min-poly-conj} and \autoref{lem:nonreci-irred-min-poly-conj} above, we can prove:
	\begin{lem}[Decomposing Using the Minimal Polynomial] \label{lem:min-poly-orthgrp}
		Let $A$ be an Artinian local ring with finite residue field~$k$ of odd characteristic~$p$. Let $(M,q)$ be a regular quadratic module over $A$. Then an isometry $\phi$ of finite order~$N$ coprime to $p$ has a separable minimal polynomial $\minpoly{\phi}$.

		Let
		\begin{align*}
			\minpoly{\phi} = P_1 \dots P_m Q_1 \reciprocal{Q}_1 \dots Q_n \reciprocal{Q}_n \in A[T]
		\end{align*}
		be the factorisation from \autoref{lem:decomp-module-ann}, which applies by \autoref{lem:minpoly-self-reciproical}.
		Then $M$ admits the orthogonal decomposition
		\begin{align*}
			M = \bigorthsum_{i = 1}^{m} M_{i, \phi}^{(\mathrm{\RN{1}})} \orthsum \bigorthsum_{j = 1}^{n} M_{j,\phi}^{(\mathrm{\RN{2}})}
		\end{align*}
		where $\smash{M_{i,\phi}^{(\mathrm{\RN{1}})} \defeq M_{P_i,\phi}}$ and $\smash{M_{j,\phi}^{(\mathrm{\RN{2}})} \defeq M_{Q_j,\phi} \oplus M_{\reciprocal{Q_j},\phi}}$ from \autoref{lem:orth-decom-minpoly}.

		Then the conjugacy class of $\phi$ in the orthogonal group $\orthgrp{}{}{M,q}$ is completely determined by $\minpoly{\phi}$ together with the isometry classes of the regular quadratic modules $\smash{M_{i,\phi}^{(\mathrm{\RN{1}})}}$, $\smash{M_{j,\phi}^{(\mathrm{\RN{2}})}}$ for all $i$ and $j$.
	\end{lem}
	\begin{proof}
		Let $\phi, \psi \colon (M,q) \to (M,q)$ be two isometries that are conjugate, i.e.\ there exists an isometry $\tau \colon (M,q) \to (M,q)$ such that
		$\tau \circ \phi = \psi \circ \tau$.
		Then, $\minpoly{\phi} = \minpoly{\psi}$ by \autoref{lem:conj-implies-same-minpoly}. For a polynomial $R \in A[T]$ dividing the minimal polynomial $\minpoly{\phi}$, the isometry $\tau$ maps $M_{R, \phi}$ to $M_{R,\psi}$ and $\tau^{-1}$ maps $M_{R, \psi}$ to $M_{R,\phi}$.

		Conversely, let $\phi, \psi \colon (M,q) \to (M,q)$ be two isometries such that
		\begin{align*}
			\minpoly{\phi} = \minpoly{\psi}
		\end{align*}
		and such that there exist isometries
		\begin{align*}
			\tau_i \colon M_{i,\phi}^{(\mathrm{\RN{1}})} \to M_{i,\psi}^{(\mathrm{\RN{1}})}, \quad \sigma_j \colon M_{j,\phi}^{(\mathrm{\RN{2}})} \to M_{j,\psi}^{(\mathrm{\RN{2}})}
		\end{align*}
		for all $i$ and $j$. Then,
		\begin{align*}
			\restr{\phi}{M_{i,\phi}^{(\mathrm{\RN{1}})}}, \quad \tau_i^{-1} \circ \restr{\psi}{M_{i,\psi}^{(\mathrm{\RN{1}})}} \circ \tau_i \qquad \text{and} \qquad
			\restr{\phi}{M_{j,\phi}^{(\mathrm{\RN{2}})}}, \quad \sigma_j^{-1}\circ \restr{\psi}{M_{j,\psi}^{(\mathrm{\RN{2}})}} \circ \sigma_j
		\end{align*}
		are isometries of $M_{i,\phi}^{(\mathrm{\RN{1}})}$ and $M_{j,\phi}^{(\mathrm{\RN{2}})}$ with minimal polynomials $P_i$ and $Q_j \reciprocal{Q}_j$, respectively. Thus, the assertion follows from \autoref{lem:selfreci-irred-min-poly-conj} and \autoref{lem:nonreci-irred-min-poly-conj}.
	\end{proof}

	\medskip

	Recall that the Witt class and the rank of a regular quadratic module determine its isometry class. It is possible to determine the Witt class of the quadratic module corresponding to the involution trace of a self-reciprocal polynomial. This is described in \autoref{sec:appendix}.
	\begin{introLem}[\textbf{\ref{lem:witt-cls-nontrivial}}]
		Let $A$ be a Henselian local ring with finite residue field~$k$ of odd characteristic, $(M,q)$ a regular quadratic module over $A$, and $\phi \colon (M,q) \to (M,q)$ an isometry with self-reciprocal minimal polynomial $P$ of degree at least $2$ such that $\overline{P}$ is irreducible over $k$. We define $E \defeq \inquo{A[T]}{\ideal{P}}$. Then
		\begin{align*}
			\wittcls{(M,q)} = (\rank{E}{M}) \Lambda.
		\end{align*}
	\end{introLem}
	Recall that $\Lambda$ is the element of $\witt{k}$ with $d(\Lambda) = 0$ and $\smash{\legendre{\Arf(\Lambda)}{k} = -1}$.

	There is one small caveat.
	\begin{rem} \label{rem:fixed-points}
		Let $A$ be an Artinian local ring with finite residue field~$k$ of odd characteristic~$p$. Let $(M,q)$ be a regular quadratic module over $A$. Then an isometry $\phi \colon (M,q) \to (M,q)$ of finite order~$N$ coprime to $p$ has a separable minimal polynomial $\minpoly{\phi}$.

		Let $P \in A[T]$ be a self-reciprocal polynomial of degree $1$ such that its reduction~$\overline{P}$ is irreducible. For $P = T - a_0 \in A[T]$ with $a_0 \in \units{A}$, we immediately obtain that $a_0 \in \{\pm1\}$ because $2 \in \units{A}$ and $A$ is local.

		If $P = T \mp 1$ divides the minimal polynomial, there are in general no restrictions on the Witt class $\wittcls{M_{P}}$, in the notation of \autoref{lem:min-poly-orthgrp}.
	\end{rem}


	\subsection{The Main Theorem} \label{subsec:main-thm}
	With all the work done in \autoref{subsec:conj-cls} and \autoref{sec:appendix}, we have established the main theorem of this paper, the classification of all conjugacy classes of isometries of a regular quadratic module over an Artinian local ring with finite residue field whose order is coprime to the (odd) characteristic of the residue field; see \autoref{thm:class-conj-classes}. We emphasise that this classification is constructive in the following sense: for chosen invariants, there is a natural realisation, i.e.\ a regular quadratic module with an isometry with the given invariants; see \autoref{rem:realisations-conj}.

	\medskip

	\begin{thm}[Parametrisation of Conjugacy Classes] \label{thm:class-conj-classes}
		Let $A$ be an Artinian local ring with maximal ideal~$\mathfrak{m}$ and finite residue field~$k$ of odd characteristic~$p$. Let $(M,q)$ be a regular quadratic module over $A$.

		Then the conjugacy classes of isometries of $(M,q)$ of order coprime to $p$ can be parametrised by
		\begin{align*}
			[(h_{+},r_{+}), (h_{-}, r_{-}), & (P_1, r_1^{\mathrm{\RN{1}}}), \dots, (P_m, r_m^{\mathrm{\RN{1}}}), (Q_1 \reciprocal{Q}_1, r_1^{\mathrm{\RN{2}}}), \dots, (Q_n \reciprocal{Q}_n, r_{n}^{\mathrm{\RN{2}}})]
		\end{align*}
		where $m, n \in \Z_{\geq 0}$, $\smash{r_i^{\mathrm{\RN{1}}}, r_j^{\mathrm{\RN{2}}} \in \Z_{\geq 1}}$, and
		\begin{enumerate}
			\item $h_{+}$ is a regular quadratic module over $A$ of rank $r_{+}$;
			\item $h_{-}$ is a regular quadratic module over $A$ of rank $r_{-}$;
			\item $\{P_i\}$ are pairwise strongly coprime polynomials over $A$ such that
			\begin{enumerate}
				\item $P_i$ is monic and self-reciprocal with $\deg{P_i} \geq 2$,
				\item the reduction $\overline{P}_i \in k[T]$ is irreducible and $P_i$ its canonical Hensel lift (see \autoref{defi:canonical-lift});
			\end{enumerate}
			\item $\{Q_j \reciprocal{Q}_j\}$ are pairwise strongly coprime polynomials over $A$ such that
			\begin{enumerate}
				\item $Q_j$ is monic and non-self-reciprocal,
				\item the reduction $\overline{Q}_j \in k[T]$ is irreducible and $Q_j$ its canonical Hensel lift;
			\end{enumerate}
		\end{enumerate}
		up to reordering of the polynomials $P_i$ and $Q_j \reciprocal{Q}_j$, respectively.
		The above data are required to satisfy
		\begin{align*}
			\rank{A}{M} = r_{+} + r_{-} + \sum_{i} r_i^{\mathrm{\RN{1}}} \deg{P_i} + \sum_{j} 2 r_j^{\mathrm{\RN{2}}} \deg{Q_j}
		\end{align*}
		and the equality in the Witt ring $\witt{A}$,
		\begin{align*}
			\wittcls{M} = \wittcls{h_{+}} + \wittcls{h_{-}} + \sum_{i} r_i^{\mathrm{\RN{1}}} \Lambda.
		\end{align*}
		In this parametrisation, the conjugacy class of an isometry~$\phi$ corresponds to the data with $h_{\pm} = \ker(\phi \mp \id)$, with the $P_i$ of degree at least~$2$ and the $Q_j \reciprocal{Q}_j$ from the factorisation of~$\minpoly{\phi}$ in \autoref{lem:decomp-module-ann}, and with $\smash{r_i^{\mathrm{\RN{1}}}}$ and $\smash{r_j^{\mathrm{\RN{2}}}}$ the multiplicities of $\smash{\overline{P}_i}$ and $\smash{\overline{Q}_j}$ in the characteristic polynomial of~$\smash{\overline{\phi}}$. The inverse map is given by the construction in \autoref{rem:realisations-conj}.
	\end{thm}
	\begin{proof}
		Let $\phi$ be an isometry of $(M,q)$ of order coprime to~$p$. By \autoref{lem:min-poly-orthgrp}, $\phi$ has a separable minimal polynomial~$\minpoly{\phi}$, the factors of $\minpoly{\phi}$ in \autoref{lem:decomp-module-ann} are canonical Hensel lifts, and the self-reciprocal factors of degree~$1$ are among $T \mp 1$ by \autoref{rem:fixed-points}. Accordingly, the orthogonal decomposition in \autoref{lem:min-poly-orthgrp} consists of the summands $\smash{h_{\pm} = M_{T \mp 1} = \ker(\phi \mp \id)}$, on which $\phi$ acts as $\pm \id$, the summands $\smash{M_{P_i}}$ for the self-reciprocal factors $P_i$ of degree at least~$2$, and the summands $\smash{M_{Q_j} \oplus M_{\reciprocal{Q}_j}}$. By (the proof of) \autoref{lem:irred-conjugated}, $\smash{M_{P_i}}$ is free of some rank $\smash{r_i^{\mathrm{\RN{1}}} \geq 1}$ over $\inquo{A[T]}{\ideal{P_i}}$, and $\smash{M_{Q_j}}$ and $\smash{M_{\reciprocal{Q}_j}}$ are free over $\inquo{A[T]}{\ideal{Q_j}}$ and $\inquo{A[T]}{\ideal{\reciprocal{Q}_j}}$, respectively, of the same rank $\smash{r_j^{\mathrm{\RN{2}}} \geq 1}$, since $\lambda_q$ identifies $\smash{M_{\reciprocal{Q}_j}}$ with the dual of the totally isotropic module $\smash{M_{Q_j}}$. Hence, $\smash{r_i^{\mathrm{\RN{1}}}}$ and $\smash{r_j^{\mathrm{\RN{2}}}}$ are the multiplicities of $\overline{P}_i$ and $\overline{Q}_j$ in the characteristic polynomial of~$\smash{\overline{\phi}}$, and the data of~$\phi$ satisfy (1) to~(4) and the rank condition. The Witt class condition follows from \autoref{lem:witt-cls-nontrivial}, as the summands $\smash{M_{Q_j} \oplus M_{\reciprocal{Q}_j}}$ are hyperbolic.

		Conjugate isometries $\phi$ and $\psi$ have the same data, since $\minpoly{\phi} = \minpoly{\psi}$ by \autoref{lem:conj-implies-same-minpoly} and a conjugating isometry maps the above summands for~$\phi$ isometrically onto those for~$\psi$. Conversely, if $\phi$ and $\psi$ have the same data, then $\minpoly{\phi} = \minpoly{\psi}$, and the corresponding summands are isometric: for $h_{\pm}$ by assumption, and otherwise since they have the same rank and, by the above, the same Witt class, which determine the isometry class by \autoref{lem:witt-cancel}. So, $\phi$ and $\psi$ are conjugate by \autoref{lem:min-poly-orthgrp}.

		Finally, for given data as in the theorem, \autoref{rem:realisations-conj} provides a regular quadratic module $(M',q')$ with an isometry~$\phi'$ of order coprime to~$p$ with these data. By the rank and Witt class conditions and \autoref{lem:witt-cls-nontrivial}, $(M',q')$ has the same rank and Witt class as $(M,q)$, so there is an isometry $\alpha \colon (M',q') \to (M,q)$ by \autoref{lem:witt-cancel}. Then $\alpha \circ \phi' \circ \alpha^{-1}$ is an isometry of $(M,q)$ of order coprime to~$p$ with the given data.
	\end{proof}
	For an isometry $\phi$ in a conjugacy class with the above invariants, the ranks $r_{+}$ and $r_{-}$ can be recovered as the multiplicities of $T \mp 1$, respectively, in the characteristic polynomial of~$\smash{\overline{\phi}}$. All these multiplicities are well defined since $T \mp 1$, $\overline{P}_i$, $\overline{Q}_j$, and $\smash{\reciprocal{\overline{Q}}_j}$ are pairwise distinct irreducible polynomials in~$k[T]$.

	\medskip

	The proof above is constructive: the inverse of the parametrisation in \autoref{thm:class-conj-classes} is given by the following explicit construction of a regular quadratic module with an isometry with prescribed invariants.
	\begin{rem} \label{rem:realisations-conj}
		Let $A$ be an Artinian local ring with maximal ideal~$\mathfrak{m}$ and finite residue field~$k$ of odd characteristic~$p$.
		If we are given a tuple
		\begin{align*}
			[(h_{+},r_{+}), (h_{-}, r_{-}), (P_1, r_1^{\mathrm{\RN{1}}}), \dots, (P_m, r_m^{\mathrm{\RN{1}}}), (Q_1 \reciprocal{Q}_1, r_1^{\mathrm{\RN{2}}}), \dots, (Q_n \reciprocal{Q}_n, r_{n}^{\mathrm{\RN{2}}})]
		\end{align*}
		as above, we build the regular quadratic module $(M,q)$ with an isometry with the given invariants as follows:

		For a self-reciprocal polynomial $P_i \in A[T]$ of degree~$\smash{d_i^{\mathrm{\RN{1}}}}$ at least $2$ such that its reduction $\overline{P}_i$ is irreducible in $k[T]$ and $P_i$ is its canonical Hensel lift, we define $\smash{E_{P_i} \defeq \inquo{A[T]}{\ideal{P_i}}}$. Then, $\smash{M_{P_i}}$ is the orthogonal sum of $\smash{r_i^{\mathrm{\RN{1}}}}$ copies of the quadratic module $\smash{M_{E_{P_i}}}$ corresponding to the involution trace $\smash{\tr_{\fieldextension{E_{P_i}}{A}}}$.

		For a polynomial $\smash{Q_j \in A[T]}$ with $\smash{Q_j \neq \reciprocal{Q}_j}$ of degree~$\smash{d_j^{\mathrm{\RN{2}}}}$ such that the reduction $\smash{\overline{Q}_j}$ is irreducible in $k[T]$ and $Q_j$ is its canonical Hensel lift, we define $\smash{\hat{M}_{Q_j}}$ as the orthogonal sum of $\smash{r_j^{\mathrm{\RN{2}}}}$ copies of the hyperbolic module $\smash{\Hyp{(\inquo{A[T]}{\ideal{Q_j}})}}$ of rank~$2 d_j^{\mathrm{\RN{2}}}$; here we consider $\inquo{A[T]}{\ideal{Q_j}}$ as a (finitely generated projective) $A$-module.

		Then, we define
		\begin{align*}
			M \defeq h_{+} \orthsum h_{-} \orthsum \bigorthsum_{i} M_{P_i} \orthsum \bigorthsum_{j} \hat{M}_{Q_j}.
		\end{align*}
		Now, we define an isometry $\phi \colon (M,q) \to (M,q)$ that acts as follows:
		\begin{enumerate}
			\item On $h_{+}$, $\phi$ acts as the identity.
			\item On $h_{-}$, $\phi$ acts as $-\id$.
			\item On $\smash{M_{P_i}}$, $\phi$ acts on each copy of $\smash{M_{E_{P_i}}}$ as the companion matrix $\companion{P_i}$ of~$P_i$; recall the definition from the comment after \autoref{lem:irred-conjugated}.
			\item On $\smash{\hat{M}_{Q_j}}$, $\phi$ acts on each copy of $\smash{\inquo{A[T]}{\ideal{Q_j}}}$ as the companion matrix $\companion{Q_j}$ of $Q_j$. This already fixes how $\phi$ acts on the hyperbolic module $\smash{\Hyp{(\inquo{A[T]}{\ideal{Q_j}})}}$; see \autoref{lem:nonreci-irred-min-poly-conj}.
		\end{enumerate}
		Thus, $\phi$ restricted to $M_{P_i}$ and $\smash{\hat{M}_{Q_j}}$ has the minimal polynomial $P_i$ and $\smash{Q_j \reciprocal{Q}_j}$, respectively. Since $T \mp 1$, $P_i$, $Q_j$, and $\smash{\reciprocal{Q}_j}$ divide $T^e - 1$ for the least common multiple~$e$ of~$2$ and the orders $\ord{\overline{P}_i}$ and $\ord{\overline{Q}_j}$ (see \autoref{lem:canonical-lift}), it follows that $\phi^e = \id$. In particular, the order of~$\phi$ is coprime to~$p$.
	\end{rem}

	Using the above remark, it is straightforward to write down examples of isometries of a given order~$N$.
	\begin{ex} \label{ex:isometries-type-one}
		For the prime $p = 7$, the irreducible self-reciprocal polynomials of degree~$2$ are
		\begin{align*}
			1 + T^2, \quad 1 + 3 T + T^2, \quad 1 + 4 T + T^2,
		\end{align*}
		which have order $4$, $8$, and $8$, respectively.

		Then, we construct quadratic spaces $M_P$ as in \autoref{rem:realisations-conj} for $P$ one of these polynomials:
		\begingroup
		\setlength{\tabcolsep}{12pt}
		\begin{align*}
			\begin{tabularx}{10cm}{r | X X X}
				$P$ & $1 + T^2$ & $1 + 3 T + T^2$ & $1 + 4 T + T^2$ \\ \hline
				$N$ & $4$ & $8$ & $8$ \\ \hline
				$\lambda$ & $\begin{bmatrix}
					2 & 0 \\
					0 & 2
				\end{bmatrix}$ & $\begin{bmatrix}
					2 & 4 \\
					4 & 2
				\end{bmatrix}$ & $\begin{bmatrix}
					2 & 3 \\
					3 & 2
				\end{bmatrix}$ \\[+5pt] \hline
				$C(P)$ & $\begin{bmatrix}
					0 & 6 \\
					1 & 4
				\end{bmatrix}$ & $\begin{bmatrix}
					0 & 6 \\
					1 & 3
				\end{bmatrix}$ & $\begin{bmatrix}
					0 & 6 \\
					1 & 3
				\end{bmatrix}$
			\end{tabularx}
		\end{align*}
		\endgroup
		Here, we listed the polynomial~$P$, its order~$N$, the symmetric bilinear form~$\lambda$, and the companion matrix~$C(P)$.
		The symmetric bilinear form~$\lambda = \lambda_{\fieldextension{E}{A}}$ of the involution trace (see \autoref{defi:involution-trace}) is written in the basis $(1, t)$, where $t$ is the image of $T$ in $E = \inquo{A[T]}{\ideal{P}}$.
		One can directly verify that, as quadratic spaces, they are all isometric because they have a nontrivial Arf invariant and the same rank.
	\end{ex}

	\medskip

	We describe how to use the classification in \autoref{thm:class-conj-classes} to obtain isometries of discriminant forms. This is a continuation of \autoref{constr:decomp-discr-form1}.
	\begin{constr} \label{constr:decomp-discr-form2}
		Let $(\Gamma, q)$ be an odd-order discriminant form with a faithful action~$\tau$ of the cyclic group $\modZ{N}$ by isometries such that $N$ is coprime to the order of~$\Gamma$. We follow the notation of \autoref{constr:decomp-discr-form1}, i.e.\ we set $\phi \coloneqq \tau[1]$. We denote the subgroup of $\orthgrp{}{}{\Gamma,q}$ generated by $\phi$ by $\Phi$; it is a cyclic group of order~$N$.

		In \autoref{constr:decomp-discr-form1}, we showed how to obtain an orthogonal decomposition into discriminant forms
		\begin{align*}
			(\Gamma, q) = \bigorthsum_{p \mid \abs{\Gamma}} \bigorthsum_{l \in \Z_{\geq 1}} (\Gamma^{(p^l)}, q^{(p^l)})
		\end{align*}
		where each $\smash{(\Gamma^{(p^l)}, q^{(p^l)})}$ is homocyclic of exponent $p^l$ and $\Phi$-invariant.

		We can change the coefficients using \autoref{lem:discr-forms-and-reg-quadr-form}, so that the restriction of $\tau$ to $\smash{(\Gamma^{(p^l)},q^{(p^l)})}$ yields an action of $\modZ{N}$ by isometries on a regular quadratic module $\smash{(M^{(p^l)}, \hat{q}^{(p^l)})}$ over $\modZ{p^l}$. Let us denote the generator of this action of $\modZ{N}$ on $\smash{(M^{(p^l)}, \hat{q}^{(p^l)})}$ by $\smash{\phi^{(p^l)}}$. As noted in \autoref{constr:decomp-discr-form1}, this action need not be faithful.

		Now, \autoref{thm:class-conj-classes} applies to the regular quadratic module $\smash{(M^{(p^l)}, \hat{q}^{(p^l)})}$ over $\modZ{p^l}$ and the isometry $\smash{\phi^{(p^l)}}$ for every prime~$p$ and positive integer~$l$ where $\smash{M^{(p^l)}}$ is nontrivial. We decompose the minimal polynomial of $\smash{\phi^{(p^l)}}$ and obtain a finite number of tuples for primes~$p$ and positive integers $l \in \Z_{\geq 1}$,
		\begin{align*}
			\big[(h_{+}^{(p^l)},r_{+}^{(p^l)}), (h_{-}^{(p^l)}, r_{-}^{(p^l)}), [(P_i^{(p^l)}, r_i^{\mathrm{\RN{1}}, {(p^l)}})]_i, [(Q_j^{(p^l)} \reciprocal{(Q_j^{(p^l)})}, r_j^{\mathrm{\RN{2}},(p^l)})]_j \big]_{p,l}
		\end{align*}
		where $\smash{r_i^{\mathrm{\RN{1}}, (p^l)}, r_j^{\mathrm{\RN{2}}, (p^l)} \in \Z_{\geq 1}}$ and
		\begin{enumerate}
			\item $\smash{h_{+}^{(p^l)}}$ is a regular quadratic module over $\modZ{p^l}$ of rank $\smash{r_{+}^{(p^l)}}$;
			\item $h_{-}^{(p^l)}$ is a regular quadratic module over $\modZ{p^l}$ of rank $\smash{r_{-}^{(p^l)}}$;
			\item $P_i^{(p^l)} \in \modZ{p^l}[T]$ are pairwise strongly coprime polynomials such that
			\begin{enumerate}
				\item $P_i^{(p^l)}$ is monic and self-reciprocal with $\smash{\deg{P_i^{(p^l)}} \geq 2}$,
				\item the reduction $\smash{\overline{P}_i^{(p^l)} \in \finitefield{p}[T]}$ is irreducible and $\smash{P_i^{(p^l)}}$ is its canonical Hensel lift;
			\end{enumerate}
			\item $Q_j^{(p^l)} \reciprocal{(Q_j^{(p^l)})} \in \modZ{p^l}[T]$ are pairwise strongly coprime polynomials such that
			\begin{enumerate}
				\item $Q_j^{(p^l)}$ is monic and non-self-reciprocal,
				\item the reduction $\smash{\overline{Q}_j^{(p^l)} \in \finitefield{p}[T]}$ is irreducible and $\smash{Q_j^{(p^l)}}$ is its canonical Hensel lift;
			\end{enumerate}
		\end{enumerate}
		up to reordering the polynomials $\smash{P_i^{(p^l)}}$ and $\smash{Q_j^{(p^l)} \reciprocal{(Q_j^{(p^l)})}}$, respectively.
		These data are required to satisfy
		\begin{align*}
			\rank{\modZ{p^l}}{M^{(p^l)}} = r_{+}^{(p^l)} + r_{-}^{(p^l)} + \sum_{i} r_i^{\mathrm{\RN{1}}, (p^l)} \deg{P_i^{(p^l)}} + \sum_{j} 2 r_j^{\mathrm{\RN{2}},(p^l)} \deg{Q_j^{(p^l)}}
		\end{align*}
		and the equalities
		\begin{align*}
			\wittcls{M^{(p^l)}} = \wittcls{h_{+}^{(p^l)}} + \wittcls{h_{-}^{(p^l)}} + \sum_{i} r_i^{\mathrm{\RN{1}}, (p^l)} \Lambda
		\end{align*}
		in the Witt ring $\smash{\witt{\modZ{p^l}}}$ for all primes $p$ and positive integers $l \in \Z_{\geq 1}$.

		Moreover, the order $\smash{N^{(p^l)}}$ of $\smash{\phi^{(p^l)}}$ is the least common multiple of the orders $\smash{\ord{\overline{P}_i^{(p^l)}}}$ and $\smash{\ord{\overline{Q}_j^{(p^l)}}}$, and of~$2$ if $\smash{h_{-}^{(p^l)} \neq 0}$. Thus, all these orders divide~$N$, and their least common multiple over all $p$ and $l$ equals~$N$, as $\tau$ is faithful; see \autoref{constr:decomp-discr-form1}.
	\end{constr}
	The above data do not depend on the choice of the decomposition in \autoref{constr:decomp-discr-form1}. Indeed, for every prime~$p$ and $l \in \Z_{\geq 1}$, the reduction modulo~$p$ of $\smash{(M^{(p^l)}, \hat{q}^{(p^l)})}$ with the action of $\smash{\phi^{(p^l)}}$ is equivariantly isometric to the $\finitefield{p}$-vector space $\smash{(\Gamma^{(p)}[p] \cap p^{l-1}\Gamma^{(p)}) / (\Gamma^{(p)}[p] \cap p^{l}\Gamma^{(p)})}$, where $\smash{\Gamma^{(p)}[p]}$ denotes the subgroup of elements $x$ with $px = 0$, with the action induced by~$\phi$ and the symmetric bilinear form induced by $\smash{(p^{l-1}x, y) \mapsto 2\lambda_q(x,y) \in \inquo{\frac{1}{p}\Z}{\Z} \cong \finitefield{p}}$; this only depends on $(\Gamma, q, \phi)$. By (the proof of) \autoref{lem:witt-iso-Henselian}, the reduction determines $\smash{(M^{(p^l)}, \hat{q}^{(p^l)})}$ up to isometry. Moreover, as in the proof of \cite[Theorem~1]{MR309963}, two isometries of order coprime to~$p$ of a regular quadratic module over $\modZ{p^l}$ are conjugate if and only if their reductions are conjugate, since the reduction map of orthogonal groups is surjective (see the beginning of \autoref{sec:orthogonal-grps}) and its kernel is a $p$-group. Hence, $\smash{\phi^{(p^l)}}$ is determined up to conjugacy by $(\Gamma, q, \phi)$, and so are the above data.

	\smallskip

	Let us stress that not every decomposition of a $p$-group into homocyclic components is invariant under the group $\Phi$ above, and that not every such decomposition is orthogonal. As explained in \autoref{subsec:decomp-discr-form}, we need to employ the results by Maschke and Schur, together with \autoref{rem:splitting-top-component}, to guarantee that an orthogonal and $\Phi$-invariant decomposition always exists; see \autoref{constr:decomp-discr-form1}.

	\smallskip

	We will specialise the above construction to fixed-point free actions in \autoref{constr:decomp-discr-form3}.

	\smallskip

	Also note that, while this might seem like a lot of data, the above construction is still suitable for efficient and concrete calculations, the reason being that irreducible polynomials over finite fields are well understood. For small degree, they can be listed in most computer algebra systems.

	Let us illustrate this point with an example of small order.
	\begin{ex} \label{ex:isometries-discr-form}
		Consider the discriminant form from \autoref{ex:discr-symbol} with symbol $(3)^{-4} (7^2)^{+2}$. We already noted that we can split this discriminant form into $(3)^{-4}$ and $(7^2)^{+2}$ and then write down all possible orthogonal decompositions of these discriminant forms into discriminant forms of smaller rank (viewed as modules over $\modZ{3}$ and $\modZ{7^2}$, respectively).

		In this case, we obtained, up to reordering, the nine decompositions
		\begin{align*}
			(3)^{-4} & = (3)^{+3} (3)^{+1} = (3)^{-3} (3)^{-1} = (3)^{+2} (3)^{-2} \\
			& = (3)^{-2} (3)^{+1} (3)^{-1} = (3)^{+2} (3)^{-1} (3)^{-1} = (3)^{+2} (3)^{+1} (3)^{+1} \\
			& = (3)^{+1} (3)^{-1} (3)^{-1} (3)^{-1} = (3)^{-1} (3)^{+1} (3)^{+1} (3)^{+1}
		\end{align*}
		and the two decompositions
		\begin{align*}
			(7^2)^{+2} = (7^2)^{+1} (7^2)^{-1}.
		\end{align*}
		Among the discriminant forms appearing above, only $(3)^{-2}$, $(3)^{-4}$, and $(7^2)^{+2}$ have a conjugacy class of isometries~$\phi$ of order greater than $2$ and coprime to $3$ and~$7$, respectively, such that $\langle \phi \rangle$ is fixed-point free; see \autoref{sec:fpf-isometries}. We exclude order~$2$ since the only fixed-point free isometry of order~$2$ is $-\id$, which exists on every discriminant form of odd order.

		Namely, we obtain the following:
		\begin{gather*}
			\begin{tabularx}{12cm}{l | r | X | X | l}
				discr. & order & irreducible polynomial & reciprocal & type \\
				\hline
				$(3)^{-2}$ & $4$ & $T^2 + 1$ & $T^2 + 1$ & \RN{1} \\
				$(3)^{-4}$ & $5$ & $T^4 + T^3 + T^2 + T + 1$ & $T^4 + T^3 + T^2 + T + 1$ & \RN{1} \\
				& $10$ & $T^4 + 2 T^3 + T^2 + 2 T + 1$ & $T^4 + 2 T^3 + T^2 + 2 T + 1$ & \RN{1}\\
				$(7^2)^{+2}$ & $3$ & $T + 31$ & $T + 19$ & \RN{2} \\
				& $6$ & $T + 30$ & $T + 18$ & \RN{2}
			\end{tabularx}
		\end{gather*}
		That is:
		\begin{enumerate}
			\item $(3)^{-2}$ has an isometry of order~$4$ corresponding to the irreducible self-reciprocal polynomial $T^2 + 1 \in \modZ{3}[T]$.
			\item $(3)^{-4}$ has an isometry of order~$5$ corresponding to the irreducible self-reciprocal polynomial $T^4 + T^3 + T^2 + T + 1 \in \modZ{3}[T]$ and an isometry of order~$10$ corresponding to the irreducible self-reciprocal polynomial $T^4 + 2 T^3 + T^2 + 2 T + 1 \in \modZ{3}[T]$.
			\item $(7^2)^{+2}$ has an isometry of order~$3$ corresponding to the pair of non-self-reciprocal polynomials $T + 31, T + 19 \in \modZ{7^2}[T]$, the canonical Hensel lifts of $T + 3, T + 5 \in \finitefield{7}[T]$, and an isometry of order~$6$ corresponding to the pair of non-self-reciprocal polynomials $T + 30, T + 18 \in \modZ{7^2}[T]$, the canonical Hensel lifts of $T + 2, T + 4 \in \finitefield{7}[T]$.
		\end{enumerate}
	\end{ex}

	\medskip

	Furthermore, we can use results on self-reciprocal irreducible polynomials to count the polynomials $P_i$ and pairs $\{Q_j, \reciprocal{Q}_j\}$ in \autoref{thm:class-conj-classes} of a given degree and order, and thus to enumerate, for a regular quadratic module of a given rank, the conjugacy classes of its isometries of a given order. Since the polynomials $P_i$ and $Q_j$ in \autoref{thm:class-conj-classes} are the canonical Hensel lifts of their reductions, it suffices to consider irreducible polynomials over the finite residue field. To this end, we recall some facts.

	\smallskip

	For a polynomial $Q \in A[T]$ of degree $d$, we define its \textit{associated} self-reciprocal polynomial by
	\begin{align*}
		Q^{\ass}(T) \defeq T^d Q(T + T^{-1}).
	\end{align*}
	It is a polynomial of degree $2d$.
	\begin{lem}[{\cite[Lemma~2.1]{MR914197}}]
		Let $k$ be a finite field.
		\begin{enumerate}
			\item If $P$ is an irreducible, self-reciprocal polynomial in $k[T]$, then $P = T \pm 1$ or $P$ is the associate of another polynomial.
			\item If $Q$ is an irreducible polynomial in $k[T]$, its associate $Q^{\ass}$ is reducible if and only if $Q^{\ass} = R \reciprocal{R}$ for some irreducible polynomial $R$ where either $R = T \pm 1$ or $R$ is not self-reciprocal.
		\end{enumerate}
	\end{lem}

	The following lemma provides a simple criterion to check whether an irreducible polynomial of given degree and order is self-reciprocal.
	\begin{lem}[{\cite[Lemma~2.3]{MR914197}}]
		Let $Q$ be an irreducible polynomial in $\finitefield{q}[T]$ of degree $2d$. Then $Q$ is self-reciprocal if and only if its order $\ord{Q}$ divides $q^d+1$.
	\end{lem}

	We summarise some general results on irreducible polynomials over finite fields. We will use these in \autoref{ex:counting-self-reciprocal} to calculate for a given degree and order the number of self-reciprocal (or non-self-reciprocal) irreducible polynomials over finite fields. Note that, by Hensel's lemma, this will allow us to implement \autoref{thm:class-conj-classes} in a computer algebra system of our choice.
	\begin{prop}[{\cite[Theorem~2.4]{MR914197}}]
		Let $Q$ be an irreducible polynomial in $\finitefield{q}[T]$ of degree $d$ and level $N$, i.e.\ the roots of $Q^{\ass}$ have order~$N$. Then:
		\begin{enumerate}
			\item If $N = 1$, then $Q(T) = T - 2$, and if $N = 2$, then $q$ is odd and $Q(T) = T + 2$.
			\item If $N > 2$, then $N = \ord{Q^{\ass}}$.
			\item $d$ is the least positive integer satisfying $q^d \equiv \pm 1 \modulo{N}$.
			\item If $Q^{\ass}$ is reducible, then $d$ is the least integer with $q^d \equiv 1 \modulo{N}$. If $Q^{\ass}$ is irreducible, then $N > 2$ and $d$ is the least integer with $\smash{q^d \equiv -1 \modulo{N}}$.
		\end{enumerate}
		For $N > 2$ coprime to~$q$, there are $\inquo{\varphi(N)}{2d}$ irreducible polynomials $Q$ over $\finitefield{q}$ of degree $d$ with level $N$, where $\varphi$ denotes Euler's totient function.
	\end{prop}

	Using the above results, it is possible to calculate the number of self-reciprocal and non-self-reciprocal irreducible polynomials of a given degree and order over a given finite field.
	\begin{ex} \label{ex:counting-self-reciprocal}
		The following table enumerates irreducible polynomials over the finite fields~$\finitefield{p}$. By \autoref{thm:class-conj-classes}, their canonical Hensel lifts are precisely the possible polynomials $P_i$ (type~\RN{1}) and $Q_j$ (type~\RN{2}) over $\modZ{p^l}$ for every $l \in \Z_{\geq 1}$.

		The columns correspond to the degree $d$ of the polynomial; the rows to the order of the finite field $\finitefield{p}$. For a given prime, we list the self-reciprocal irreducible polynomials other than $T \pm 1$ (which correspond to $h_{\pm}$ in \autoref{thm:class-conj-classes}) in the first and the non-self-reciprocal irreducible polynomials in the second line. We list the order of the polynomials with the number of polynomials of that order as index; for type~\RN{2}, the index counts the unordered pairs $\{Q, \reciprocal{Q}\}$. For example, there are $20$ such pairs, i.e.\ $40$ non-self-reciprocal irreducible polynomials, of order~$366$ and degree~$3$ over $\finitefield{13}$.

		\begin{equation*}
			\begin{tabularx}{.95\textwidth}{| r | l |
					X | X | X |}
				\hline
				$p$ \textbackslash{} $d$ & & $1$ & $2$ & $3$ \\
				[0.5ex]
				\hline\hline
				$3$ & \RN{1} & & $4_{1}$ & \\
				\hline
				& \RN{2} & & $8_{1}$ & $13_{2}$, $26_{2}$ \\
				[0.5ex]
				\hline\hline
				$5$ & \RN{1} & & $3_{1}$, $6_{1}$ & \\
				\hline
				& \RN{2} & $4_{1}$ & $8_{1}$, $12_{1}$, $24_{2}$ & $31_{5}$, $62_{5}$,
				$124_{10}$ \\
				[0.5ex]
				\hline\hline
				$7$ & \RN{1} & & $4_{1}$, $8_{2}$ & \\
				\hline
				& \RN{2} & $3_{1}$, $6_{1}$ & $12_{1}$, $16_{2}$, $24_{2}$, $48_{4}$ & $9_{1}$,
				$18_{1}$, $19_{3}$, $38_{3}$, $57_{6}$, $114_{6}$, $171_{18}$,
				$342_{18}$ \\
				[0.5ex]
				\hline\hline
				$11$ & \RN{1} & & $3_{1}$, $4_{1}$, $6_{1}$, $12_{2}$ & \\
				\hline
				& \RN{2} & $5_{2}$, $10_{2}$ & $8_{1}$, $15_{2}$, $20_{2}$, $24_{2}$, $30_{2}$,
				$40_{4}$, $60_{4}$, $120_{8}$ & $7_{1}$, $14_{1}$, $19_{3}$, $35_{4}$,
				$38_{3}$, $70_{4}$, $95_{12}$,
				$133_{18}$, $190_{12}$, $266_{18}$, $665_{72}$, $1330_{72}$ \\ [0.5ex]
				\hline\hline
				$13$ & \RN{1} & & $7_{3}$, $14_{3}$ & \\
				\hline
				& \RN{2} & $3_{1}$, $4_{1}$, $6_{1}$, $12_{2}$ & $8_{1}$, $21_{3}$, $24_{2}$,
				$28_{3}$, $42_{3}$, $56_{6}$, $84_{6}$, $168_{12}$ & $9_{1}$, $18_{1}$,
				$36_{2}$, $61_{10}$, $122_{10}$, $183_{20}$, $244_{20}$,
				$366_{20}$, $549_{60}$, $732_{40}$, $1098_{60}$, $2196_{120}$ \\
				\hline
			\end{tabularx}
		\end{equation*}
	\end{ex}


	\section{Fixed-Point Free Isometries over Artinian Rings} \label{sec:fpf-isometries}

	In this section, we define fixed-point free automorphisms over local rings; see \autoref{defi:fpf-automorphism}. We will see that involution traces often have fixed-point free actions, i.e.\ actions by fixed-point free automorphisms; see \autoref{lem:Henselian-fpf}.

	We will apply the main theorem, \autoref{thm:class-conj-classes}, to classify all automorphisms of finite order coprime to the characteristic $\characteristic{k}$ that generate a fixed-point free group of isometries for quadratic modules over a finite local ring $A$ with maximal ideal~$\mathfrak{m}$ and residue field~$k$; see \autoref{thm:class-fpf-isos}. We will see that all of them correspond to lifts of such automorphisms over the finite residue field~$k$. These lifts have the same order since their minimal polynomials are canonical Hensel lifts (see \autoref{rem:realisations-conj}), and we will need the finiteness assumption on $A$ to guarantee that they are again fixed-point free; see \autoref{cor:fpf-reduction}. Even though every Artinian local ring with finite residue field is a finite local ring, we will explicitly speak of finite local rings here to make clear that this is important.

	In the end, we want to apply \autoref{thm:class-fpf-isos} to the finite local rings $\modZ{p^l}$ for an odd prime $p \in \Z$ and a positive integer $l \in \Z_{\geq 1}$.

	\medskip


	\subsection{Fixed-Point Free Automorphisms} \label{subsec:fpf-auto}
	In this section, we introduce \textit{fixed-point free} automorphisms of modules over local rings. Of course, every linear automorphism fixes $0$, but the terminology is nevertheless standard in group theory; see e.g.\ \cite[Section~10.5]{MR1357169} or \cite[Chapter~10]{MR569209}.

	\medskip

	We start by defining \textit{fixed-point free} automorphisms. We will see in \autoref{subsec:class-fpf-iso} that the classification of conjugacy classes in orthogonal groups can be refined to the conjugacy classes of isometries that are fixed-point free or those that generate fixed-point free groups; see \autoref{thm:class-fpf-isos}.
	\begin{defi}[Fixed-Point Free Automorphism] \label{defi:fpf-automorphism}
		Let $A$ be a local ring and $M$ a finitely generated projective $A$-module. Let $\phi \colon M \to M$ be an automorphism of finite order~$N$. Then $\phi$ is called \textit{fixed-point free} if $\phi(x) = x$ for $x \in M$ implies $x = 0$.
	\end{defi}
	Of course, one could define fixed-point free automorphisms in far greater generality, but the above definition suffices for our purposes.

	A subgroup of automorphisms $\Phi$ in $\Aut_{A}(M)$ is called \textit{fixed-point free} if every nontrivial element of $\Phi$ is fixed-point free.

	\medskip

	We immediately obtain many of the properties of fixed-point free automorphisms of groups, here stated for modules over rings; see e.g.\ \cite[Lemma~10.5.1]{MR1357169}.
	\begin{lem} \label{lem:fixed-point-properties}
		Let $A$ be a local ring and $M$ a finitely generated projective $A$-module. Let $\phi \colon M \to M$ be a fixed-point free automorphism of finite order~$N$. Then
		\begin{enumerate}
			\item If $i$ is coprime to $N$, then $\phi^i$ is fixed-point free.
			\item For all $x \in M$,
			\begin{align*}
				x + \phi(x) + \dots + \phi^{N-1}(x) = 0.
			\end{align*}
		\end{enumerate}
	\end{lem}
	In particular, if $N$ is prime, any fixed-point free automorphism generates a fixed-point free group of automorphisms where all elements but the identity are fixed-point free.\footnote{For composite~$N$, this is not true in general: if $2 \in \units{A}$, the automorphism $\phi = (-\id) \oplus t$ of $M = A \oplus \inquo{A[T]}{\ideal{T^2+1}}$, where $t$ acts by multiplication with the image of~$T$, is fixed-point free of order~$4$, but $\phi^2$ fixes $A \oplus 0$.}
	In this work, we will be interested in automorphisms that generate a fixed-point free group of automorphisms because they have the nice property that all nontrivial orbits have the same length $N$; see \autoref{lem:orbit-lengths}.

	\medskip

	In particular, groups of fixed-point free automorphisms have the property that they restrict nicely to invariant submodules. This is important in our classification below, \autoref{thm:class-fpf-isos}, as it ensures that the restrictions to all invariant submodules have the same order.
	\begin{lem}[{cf.\ \cite[Proposition~3.3]{mayr-thesis}}] \label{lem:restr-fpf}
		Let $A$ be a local ring and $M$ a finitely generated projective $A$-module. Let $\Phi$ be a finite fixed-point free group of automorphisms of~$M$ and $M^{\prime}$ a nontrivial $\Phi$-invariant submodule of $M$. Then the group of automorphisms restricted to $M^{\prime}$ is isomorphic to $\Phi$.
	\end{lem}

	\medskip

	Over finite local rings, many properties of fixed-point free automorphisms can be proved with the following equivalent definition. The analogous statement for finite groups is standard; see e.g.\ \cite[Lemma~10.5.1]{MR1357169}.

	\smallskip

	We will prove that an automorphism of a module over a finite local ring is fixed-point free if and only if its reduction to the residue field is fixed-point free; see \autoref{cor:fpf-reduction}. This is a corollary of the following lemma.
	\begin{lem}
		Let $A$ be a finite local ring with maximal ideal~$\mathfrak{m}$ and finite residue field~$k$ of characteristic~$p$. Let $M$ be a finitely generated projective $A$-module. Let $\phi \colon M \to M$ be an automorphism of finite order.
		Then $\phi$ is fixed-point free if and only if the endomorphism of $M$ defined by
		\begin{align*}
			W_\phi \colon M \to M, \quad x \mapsto x - \phi(x)
		\end{align*}
		is an automorphism.
	\end{lem}
	\begin{proof}
		Assume that $\phi$ is fixed-point free and that $x - \phi(x) = y - \phi(y)$ for some $x,y \in M$. Then $\phi$ fixes $x-y$, so $x = y$. Thus, $W_{\phi}$ is injective. Since $M$ is finite, $\smash{W_{\phi}}$ is also surjective, i.e.\ an automorphism. Conversely, if $\smash{W_{\phi}}$ is an automorphism and $\phi(x) = x$, then $\smash{W_{\phi}(x) = 0}$ and hence $x = 0$.
	\end{proof}

	\smallskip

	By applying the above lemma to $\phi$ and to its reduction $\smash{\overline{\phi}}$ (note that $\smash{\overline{W_{\phi}} = W_{\overline{\phi}}}$) and using the fact that a linear endomorphism over a local ring is an isomorphism if and only if its reduction is one (see \autoref{lem:isos-preserved-by-change-rings}), we obtain the corresponding criterion for fixed-point free automorphisms.
	\begin{cor} \label{cor:fpf-reduction}
		Let $A$ be a finite local ring with maximal ideal~$\mathfrak{m}$ and finite residue field~$k$ of characteristic~$p$. Let $M$ be a finitely generated projective $A$-module. Let $\phi \colon M \to M$ be an automorphism of finite order.
		Then $\phi$ is fixed-point free if and only if its reduction $\smash{\overline{\phi}}$ over $k$ is fixed-point free.
	\end{cor}
	This corollary will be useful later because it will allow us to lift fixed-point free automorphisms of modules over fields; see \autoref{lem:Henselian-fpf}. Here, we can again use the fact that the polynomial ring over a field is a principal ideal domain.


	\subsection{Classification of Fixed-Point Free Isometries over Artinian Rings} \label{subsec:class-fpf-iso}
	In this section, we classify conjugacy classes of fixed-point free isometries, and of isometries generating fixed-point free groups, of regular quadratic modules over Artinian local rings with finite residue fields whose orders are \textit{coprime to the characteristic} of the residue field; see \autoref{thm:class-fpf-isos}. As we already mentioned, these rings are finite local rings. For the proof, the finiteness assumption enters in several ways. Most importantly, we will use it to ensure that we can lift a fixed-point free automorphism of modules over the residue field to a \textit{fixed-point free} automorphism of modules over the finite local ring of the same order; see \autoref{lem:Henselian-fpf}.

	\medskip

	Now, we use the fact that the reduction map is surjective and that, for a \textit{finite ring} $A$, its kernel is a $p$-group where $p$ is the characteristic of the residue field $k$ of~$A$. We will use this in the proof of \autoref{thm:class-fpf-isos}.
	\begin{lem}[{see e.g.\ \cite{MR309963}}]
		Let $A$ be a finite local ring with maximal ideal~$\mathfrak{m}$ and finite residue field~$k$ of characteristic~$p$. Then the reduction map
		\begin{align*}
			\GL_{n}(A) \to \GL_{n}(k)
		\end{align*}
		is a surjective group homomorphism, whose kernel is a (solvable) $p$-group.
	\end{lem}
	We can thus conclude that an automorphism $\phi \in \GL_{n}(A)$ of order coprime to $p$ has the same order as its reduction $\smash{\overline{\phi}} \in \GL_{n}(k)$.

	The above lemma can be used to obtain (nonconstructive) results on the conjugacy classes in the general linear group over a finite local ring.
	\begin{rem} \label{rem:Pomfret}
		In fact, this way Pomfret proved that two elements of $\GL_{n}(A)$ of order coprime to~$p$ are conjugate if and only if their reductions are conjugate in $\GL_{n}(k)$; see \cite[Theorem~1]{MR309963}. As every element of $\GL_{n}(k)$ of order coprime to~$p$ lifts to an element of $\GL_{n}(A)$ of the same order by the Schur-Zassenhaus theorem, the conjugacy classes of elements of order coprime to~$p$ in $\GL_{n}(A)$ can thus be identified with those in $\GL_{n}(k)$.

		If $p$ is odd, the same argument applies to the reduction
		$\orthgrp{}{A}{M,q} \to \orthgrp{}{k}{\overline{M},\overline{q}}$ for a regular quadratic module $(M,q)$ over~$A$, which is surjective (see the beginning of \autoref{sec:orthogonal-grps}) and whose kernel is again a $p$-group. Hence, the conjugacy classes of isometries of $(M,q)$ of order coprime to~$p$ can be identified with those of $(\overline{M},\overline{q})$.
	\end{rem}

	\medskip

	Now, the only thing left to do is to show that irreducible polynomials over fields generate fixed-point free actions. The proof relies on the fact that automorphisms over finite local rings are fixed-point free if and only if their reductions are fixed-point free; see \autoref{cor:fpf-reduction}.
	\begin{lem} \label{lem:Henselian-fpf}
		Let $A$ be a finite local ring with maximal ideal~$\mathfrak{m}$ and finite residue field~$k$ of odd characteristic~$p$. Let $P \in A[T]$ be a polynomial such that its reduction $\smash{\overline{P} \in k[T]}$ is irreducible with $\smash{\ord{\overline{P}} = N}$ and $P$ is the canonical Hensel lift of~$\smash{\overline{P}}$. Let $E \defeq \inquo{A[T]}{\ideal{P}}$ and denote the image of $T$ in $E$ by $t$. Then multiplication with $t$ generates a fixed-point free cyclic group of automorphisms of $E$ of order~$N$.
	\end{lem}
	\begin{proof}
		Since $P$ divides $T^N - 1$, while $\smash{\overline{P}}$ does not divide $T^i - 1$ for $0 < i < N$, multiplication with~$t$ has order~$N$. By \autoref{cor:fpf-reduction}, it thus suffices to show that $\smash{\overline{t}^i}$ acts fixed-point freely on $\smash{\overline{E} = \inquo{k[T]}{\ideal{\overline{P}}}}$ for all $i = 1, \dots, N-1$.

		Let $\smash{\overline{x} \in \overline{E}}$ be a fixed point of $\smash{\overline{t}^i}$ for some $i = 1, \dots, N-1$. We pick an $\smash{\overline{S} \in k[T]}$ such that its image in $\smash{\overline{E}}$ equals $\overline{x}$. Then
		\begin{align*}
			\overline{t}^i \overline{x} = \overline{x}
		\end{align*}
		implies
		\begin{align*}
			\overline{P} \mid (T^i - 1) \overline{S}.
		\end{align*}
		Since the order of $\smash{\overline{P}}$ is $N$ and $\smash{\overline{P}}$ is irreducible, $\smash{\overline{P}}$ divides $\smash{\overline{S}}$, i.e.\ $\smash{\overline{x} = 0}$.
	\end{proof}
	For example, the isometries in \autoref{ex:isometries-type-one} generate fixed-point free groups.

	\medskip

	We conclude with the following classification. Here, we call an isometry of a regular quadratic module \textit{fixed-point free} if the underlying automorphism is fixed-point free. Note that an isometry~$\phi$ of order coprime to~$p$ is fixed-point free if and only if $h_{+} = 0$ in its invariants from \autoref{thm:class-conj-classes}, since $\phi$ acts as the identity on~$h_{+}$, while $T - 1$ is strongly coprime to the minimal polynomial of $\phi$ on all other components. Hence, \autoref{thm:class-conj-classes} directly describes the conjugacy classes of fixed-point free isometries. The following theorem concerns the stronger condition that the cyclic group $\langle \phi \rangle$ is fixed-point free. We exclude $N = 1$ and $N = 2$ as there is only one irreducible polynomial of order~$1$ and one of order~$2$ (namely $T \mp 1$).

	\smallskip

	The result follows from the main theorem, \autoref{thm:class-conj-classes}, and \autoref{lem:Henselian-fpf}.
	\begin{thm}[Conj.\ Classes of Isometries Generating Fixed-Point Free Groups] \label{thm:class-fpf-isos}
		Let $A$ be a finite local ring with maximal ideal~$\mathfrak{m}$ and finite residue field~$k$ of odd characteristic~$p$. Let $(M,q)$ be a regular quadratic module over $A$.

		Then conjugacy classes of isometries $\phi$ of $(M,q)$ of order $N > 2$ with $p \nmid N$ such that $\langle \phi \rangle$ is fixed-point free can be parametrised by
		\begin{align*}
			[(P_1, r_1^{\mathrm{\RN{1}}}), \dots, (P_m, r_m^{\mathrm{\RN{1}}}), (Q_1 \reciprocal{Q}_1, r_1^{\mathrm{\RN{2}}}), \dots, (Q_n\reciprocal{Q}_n, r_{n}^{\mathrm{\RN{2}}})]
		\end{align*}
		where $m, n \in \Z_{\geq 0}$, $\smash{r_i^{\mathrm{\RN{1}}}, r_j^{\mathrm{\RN{2}}} \in \Z_{\geq 1}}$, and
		\begin{enumerate}
			\item $P_i \in A[T]$ are pairwise strongly coprime polynomials such that
			\begin{enumerate}
				\item $P_i$ is monic and self-reciprocal with $\deg{P_i} \geq 2$,
				\item the reduction $\overline{P}_i \in k[T]$ is irreducible of order $\ord{\overline{P}_i} = N$ and $P_i$ is its canonical Hensel lift;
			\end{enumerate}
			\item $Q_j \reciprocal{Q}_j \in A[T]$ are pairwise strongly coprime polynomials such that
			\begin{enumerate}
				\item $Q_j$ is monic and non-self-reciprocal,
				\item the reduction $\overline{Q}_j \in k[T]$ is irreducible of order $\ord{\overline{Q}_j} = N$ and $Q_j$ is its canonical Hensel lift;
			\end{enumerate}
		\end{enumerate}
		up to reordering the polynomials $P_i$ and $Q_j \reciprocal{Q}_j$, respectively.
		These data are required to satisfy
		\begin{align*}
			\rank{A}{M} = \sum_{i} r_i^{\mathrm{\RN{1}}} \deg{P_i} + \sum_{j} 2 r_j^{\mathrm{\RN{2}}} \deg{Q_j}
		\end{align*}
		and the equality in the Witt ring $\witt{A}$
		\begin{align*}
			\wittcls{M} = \sum_{i} r_i^{\mathrm{\RN{1}}} \Lambda.
		\end{align*}
	\end{thm}
	\begin{proof}
		The parametrisation is the restriction of the one in \autoref{thm:class-conj-classes}; see \autoref{rem:realisations-conj}. Let $\phi$ be an isometry of order~$N > 2$ such that $\langle \phi \rangle$ is fixed-point free. By \autoref{lem:restr-fpf}, the restriction of $\phi$ to every nontrivial $\phi$-invariant submodule has order~$N$. Since $\phi$ acts on $h_{\pm}$ with order at most~$2 < N$, we obtain $h_{\pm} = 0$. Moreover, the restriction of $\phi$ to $M_{P_i,\phi}$ (respectively $M_{Q_j,\phi}$) from \autoref{lem:min-poly-orthgrp} has order coprime to $p$ and hence the same order as its reduction, namely $\smash{\ord{\overline{P}_i}}$ (respectively $\smash{\ord{\overline{Q}_j}}$). Thus, all these orders are equal to~$N$.

		Conversely, assume that $h_{\pm} = 0$ and $\smash{\ord{\overline{P}_i} = \ord{\overline{Q}_j} = \ord{\reciprocal{\overline{Q}}_j} = N}$ for all $i,j$. By \autoref{lem:min-poly-orthgrp} and (the proof of) \autoref{lem:irred-conjugated}, $M$ is then, as an $A[T]$-module with $T$ acting by $\phi$, a direct sum of copies of $\inquo{A[T]}{\ideal{P_i}}$, $\inquo{A[T]}{\ideal{Q_j}}$, and $\inquo{A[T]}{\ideal{\reciprocal{Q}_j}}$. On each of these summands, multiplication by~$T$ generates a fixed-point free group of order~$N$ by \autoref{lem:Henselian-fpf}, since $P_i$, $Q_j$, and $\reciprocal{Q}_j$ are canonical Hensel lifts by \autoref{lem:decomp-module-ann} and \autoref{lem:canonical-lift}. Hence, $\langle \phi \rangle$ is fixed-point free of order~$N$.
	\end{proof}
	Note that here the reductions of all polynomials are of order~$N$ as the restriction of the isometry to invariant submodules needs to be of the same order by \autoref{lem:restr-fpf}. Furthermore, since $N \geq 3$, we have $h_{\pm} = 0$, so that $\rank{A}{M}$ is even by the rank condition above and \autoref{rem:self-reciprocal-constant}; hence, the order $\abs{M} = \abs{A}^{\rank{A}{M}}$ of $M$ is a perfect square.

	\medskip

	\begin{rem}
		We have classified the isometries with the above invariants up to conjugacy, but not the cyclic subgroups they generate up to conjugacy. For automorphisms of finite abelian groups that generate fixed-point free groups, the latter is discussed in \cite[Corollary~6.13 and Proposition~6.17]{mayr-thesis}. For our applications, this question does not appear to be relevant.
	\end{rem}

	We end this section with an example.
	\begin{ex}
		We construct an example with order~$57$ on a regular quadratic module over the ring $A = \modZ{7^3}$ with residue field $\finitefield{7}$.
		We find that the least positive integer $d$ with $p^d \equiv \pm 1 \modulo{57}$ is given by $d = 3$ with $p^d \equiv 1 \modulo{57}$. Then there are $\inquo{\varphi(57)}{2d} = 6$ pairs $\{Q, \reciprocal{Q}\}$, i.e.\ $12$ irreducible non-self-reciprocal polynomials over $\finitefield{7}$ of degree~$3$ and order~$57$. For example, we can pick
		\begin{align*}
			\overline{Q} = T^3 + T^2 + 4 T + 6 \in \finitefield{7}[T]
		\end{align*}
		with reciprocal polynomial
		\begin{align*}
			\reciprocal{\overline{Q}} = T^3 + 3 T^2 + 6 T + 6 \in \finitefield{7}[T].
		\end{align*}
		Then, by Hensel's lemma, the canonical Hensel lift $Q \in A[T]$ of~$\overline{Q}$, i.e.\ the unique monic lift of~$\overline{Q}$ dividing $T^{57}-1 \in A[T]$, is given by
		\begin{align*}
			Q = T^3 + 64 T^2 + 32 T + 342.
		\end{align*}
		Thus, we can construct $\smash{\hat{M}_Q} = \Hyp{U}$ with $U \defeq \inquo{A[T]}{\ideal{Q}}$ as in \autoref{rem:realisations-conj}. In the basis $(1, t, t^2)$ of~$U$, followed by the dual basis of~$\rd{U}$, the isometry acts as
		\begin{align*}
			\begin{bmatrix}
				0 & 0 & 1 & 0 & 0 & 0 \\
				1 & 0 & 311 & 0 & 0 & 0 \\
				0 & 1 & 279 & 0 & 0 & 0 \\
				0 & 0 & 0 & 32 & 64 & 1 \\
				0 & 0 & 0 & 1 & 0 & 0 \\
				0 & 0 & 0 & 0 & 1 & 0
			\end{bmatrix}
		\end{align*}
		on~$\smash{\hat{M}_Q}$, with the companion matrix $C(Q)$ and its inverted transpose $\smash{\transpose{(\companion{Q}^{-1})}}$ on the diagonal. As $\smash{\overline{Q}}$ and $\smash{\reciprocal{\overline{Q}}}$ are irreducible of order~$57$, this isometry generates a fixed-point free group of order~$57$ by \autoref{lem:Henselian-fpf}, applied to $Q$ and to $\reciprocal{Q}$ (the minimal polynomial of the second block, which is the canonical Hensel lift of~$\smash{\reciprocal{\overline{Q}}}$).
	\end{ex}


	\subsection{Fixed-Point Free Isometries of Discriminant Forms} \label{subsec:fpf-discr-form}
	In this section, we discuss fixed-point free automorphisms of groups and discriminant forms. Mainly, we apply \autoref{thm:class-fpf-isos} to discriminant forms using the decomposition from \autoref{subsec:decomp-discr-form}; see \autoref{constr:decomp-discr-form1}.
	With \autoref{lem:orbit-lengths} and \autoref{lem:cohomology-fixed-point-free}, we highlight two possible reasons why we study fixed-point free actions; see \cite{paper2}.

	\medskip

	For an introduction to fixed-point free automorphisms of finite groups, we refer to \cite[Section~10.5]{MR1357169} or \cite[Chapter~10]{MR569209}. A very explicit treatment can be found in \cite{mayr-thesis}.

	\begin{defi}[Fixed-Point Free Automorphism] \label{defi:fixed-point-free}
		Let $\Pi$ be a finite group. Then an automorphism $\alpha \colon \Pi \to \Pi$ is called \textit{fixed-point free} if $\alpha(x) = x$ implies $x = 1$.
	\end{defi}
	A subgroup $\Phi$ of $\Aut(\Pi)$ is called \textit{fixed-point free} if every nontrivial element of~$\Phi$ is fixed-point free (as for modules; see below \autoref{defi:fpf-automorphism}). We call an action of a group $\Psi$ on $\Pi$ by automorphisms \textit{fixed-point free} if every nontrivial element of~$\Psi$ acts by a fixed-point free automorphism on~$\Pi$. If $\Pi$ is nontrivial, such an action is always faithful.

	\medskip

	We recall some basic properties of fixed-point free automorphisms of finite groups. These are the already mentioned analogues of \autoref{lem:fixed-point-properties} for finite groups.
	\begin{lem}
		Let $\Pi$ be a finite group, and let $\alpha$ be a fixed-point free automorphism of $\Pi$ of order~$N$. Then
		\begin{enumerate}
			\item If $i$ is coprime to $N$, then $\alpha^i$ is fixed-point free.
			\item For all $x \in \Pi$,
			\begin{align*}
				x \, \alpha(x) \cdots \alpha^{N-1}(x) = 1.
			\end{align*}
		\end{enumerate}
	\end{lem}

	\medskip

	Our primary motivations to study fixed-point free automorphisms of groups are expressed in the following two lemmas. We further elaborate on this in \cite{paper2}.
	\begin{lem}[{see e.g.\ \cite[Proposition~3.1]{mayr-thesis}}] \label{lem:orbit-lengths}
		Let $\Pi$ be a finite group with a fixed-point free action $\tau$ of $\modZ{N}$ on $\Pi$. Then the orbit of every nontrivial element $x \in \Pi$ is of length $N$.
	\end{lem}
	\begin{proof}
		Assume that there is an element $x \in \Pi$ such that its orbit under $\tau$ is not of length~$N$. Then there exist $i,j \in \modZ{N}$ with $\tau[i](x) = \tau[j](x)$ but $i \neq j$. Thus, $x$ is fixed by the nontrivial element $\tau[i-j]$. We obtain $x = 1$.
	\end{proof}
	In particular, the order $\abs{\Pi}$ is coprime to $N$. More precisely,
	\begin{align*}
		\abs{\Pi} \equiv 1 \modulo{N}.
	\end{align*}

	\medskip

	Furthermore, fixed-point free actions have the property that all corresponding cohomology groups vanish. Here, we prove it for actions by cyclic groups.
	\begin{lem} \label{lem:cohomology-fixed-point-free}
		Let $\Pi$ be a finite abelian group with a fixed-point free action $\tau$ of~$\modZ{N}$. Then all cohomology groups $H^n(\modZ{N}, \Pi)$ for $n \in \Z_{\geq 0}$ are trivial.
	\end{lem}
	\begin{proof}
		We define $\alpha \defeq \tau[1]$.
		By definition, $\Pi^{\modZ{N}}$ is trivial since $\alpha$ is fixed-point free. Similarly, for all elements $x \in \Pi$, the norm element, i.e.\ the formal sum of all powers of $\alpha$, sends $x$ to $0$, i.e.\
		\begin{align*}
			x + \alpha(x) + \dots + \alpha^{N-1}(x) = 0.
		\end{align*}
		Since $1 - \alpha$ is injective and $\Pi$ is finite, $1 - \alpha$ is a permutation of $\Pi$; it is in particular surjective, and hence
		\begin{align*}
			(1 - \alpha) \Pi = \Pi.
		\end{align*}
		This already implies, by using the standard periodic resolution, that every cohomology group is trivial. For details, we refer to \cite[Theorem~6.2.2]{MR1269324}.
	\end{proof}

	\medskip

	Let us describe how we can use the classification in \autoref{thm:class-fpf-isos} to obtain fixed-point free actions of cyclic groups on discriminant forms. This is a special case of \autoref{constr:decomp-discr-form2}, which builds on \autoref{constr:decomp-discr-form1}.
	\begin{constr} \label{constr:decomp-discr-form3}
		Let $(\Gamma, q)$ be a nontrivial odd-order discriminant form with a fixed-point free action $\tau$ of the cyclic group $\modZ{N}$ by isometries, where $N \geq 3$. We will denote $\tau[1]$ by $\phi$. Since $\tau$ is fixed-point free, the power $\phi^i = \tau[i]$ is fixed-point free for all $i = 1, \dots, N-1$. In particular, $\tau$ is faithful, and the group $\Phi$ generated by $\phi$ is a fixed-point free cyclic group of order~$N$. Moreover, $N$ is coprime to the order of $\Gamma$ since $\abs{\Gamma} \equiv 1 \modulo{N}$ by \autoref{lem:orbit-lengths}.

		We showed in \autoref{constr:decomp-discr-form1} how to obtain an orthogonal decomposition of~$(\Gamma, q)$ into discriminant forms,
		\begin{align*}
			(\Gamma, q) = \bigorthsum_{p \mid \abs{\Gamma}} (\Gamma^{(p)}, q^{(p)}) = \bigorthsum_{p \mid \abs{\Gamma}} \bigorthsum_{l \in \Z_{\geq 1}} (\Gamma^{(p^l)}, q^{(p^l)})
		\end{align*}
		where $\smash{\Gamma^{(p^l)}}$ is trivial or a homocyclic group of exponent $p^l$ for primes~$p$ and positive integers $l \in \Z_{\geq 1}$, and every $\smash{\Gamma^{(p^l)}}$ is $\Phi$-invariant.

		Using the change of coefficients in \autoref{lem:discr-forms-and-reg-quadr-form}, the restriction of $\tau$ to $\smash{(\Gamma^{(p^l)},q^{(p^l)})}$ yields an action of $\modZ{N}$ by isometries on a regular quadratic module $\smash{(M^{(p^l)}, \hat{q}^{(p^l)})}$ over $\modZ{p^l}$. Let us denote the generator of this action of $\modZ{N}$ on $\smash{(M^{(p^l)}, \hat{q}^{(p^l)})}$ by~$\smash{\phi^{(p^l)}}$. If $\smash{\Gamma^{(p^l)}}$ is nontrivial, every nontrivial element of $\Phi$ restricts to a fixed-point free, in particular nontrivial, automorphism of $\smash{\Gamma^{(p^l)}}$ (cf.\ \autoref{lem:restr-fpf}). Hence, $\smash{\phi^{(p^l)}}$ has order~$N$, and $\smash{\langle \phi^{(p^l)} \rangle}$ is fixed-point free.

		Now, \autoref{thm:class-fpf-isos} applies to the regular quadratic module $\smash{(M^{(p^l)}, \hat{q}^{(p^l)})}$ over $\smash{\modZ{p^l}}$ and the isometry $\smash{\phi^{(p^l)}}$ for every prime~$p$ and positive integer $l$ where $\smash{M^{(p^l)}}$ is nontrivial. We can decompose the minimal polynomial of $\smash{\phi^{(p^l)}}$ to obtain a finite set of tuples
		\begin{align*}
			\big[[(P_i^{(p^l)}, r_i^{\mathrm{\RN{1}}, {(p^l)}})]_i, [(Q_j^{(p^l)} \reciprocal{(Q_j^{(p^l)})}, r_j^{\mathrm{\RN{2}},(p^l)})]_j \big]_{p,l}
		\end{align*}
		where $\smash{r_i^{\mathrm{\RN{1}}, (p^l)}, r_j^{\mathrm{\RN{2}}, (p^l)} \in \Z_{\geq 1}}$ and
		\begin{enumerate}
			\item $\smash{P_i^{(p^l)} \in \modZ{p^l}[T]}$ are pairwise strongly coprime polynomials such that
			\begin{enumerate}
				\item $P_i^{(p^l)}$ is monic and self-reciprocal with $\deg{P_i^{(p^l)}} \geq 2$,
				\item the reduction $\smash{\overline{P}_i^{(p^l)} \in \finitefield{p}[T]}$ is irreducible with order $\smash{\ord{\overline{P}_i^{(p^l)}} = N}$ and $\smash{P_i^{(p^l)}}$ is its canonical Hensel lift;
			\end{enumerate}
			\item $\smash{Q_j^{(p^l)} \reciprocal{(Q_j^{(p^l)})} \in \modZ{p^l}[T]}$ are pairwise strongly coprime polynomials such that
			\begin{enumerate}
				\item $Q_j^{(p^l)}$ is monic and non-self-reciprocal,
				\item the reduction $\smash{\overline{Q}_j^{(p^l)} \in \finitefield{p}[T]}$ is irreducible with order $\smash{\ord{\overline{Q}_j^{(p^l)}} = N}$ and $\smash{Q_j^{(p^l)}}$ is its canonical Hensel lift;
			\end{enumerate}
		\end{enumerate}
		up to reordering the polynomials $\smash{P_i^{(p^l)}}$ and $\smash{Q_j^{(p^l)} \reciprocal{(Q_j^{(p^l)})}}$, respectively.
		These data are required to satisfy
		\begin{align*}
			\rank{\modZ{p^l}}{M^{(p^l)}} = \sum_{i} r_i^{\mathrm{\RN{1}}, {(p^l)}} \deg{P_i^{(p^l)}} + \sum_{j} 2 r_j^{\mathrm{\RN{2}}, {(p^l)}} \deg{Q_j^{(p^l)}}
		\end{align*}
		and the equalities
		\begin{align*}
			\wittcls{M^{(p^l)}} = \sum_{i} r_i^{\mathrm{\RN{1}}, {(p^l)}} \Lambda
		\end{align*}
		in the Witt ring $\smash{\witt{\modZ{p^l}}}$ for all primes $p$ and positive integers $l \in \Z_{\geq 1}$.
	\end{constr}
	Note that here the reductions of all polynomials involved have order~$N$, in contrast to \autoref{constr:decomp-discr-form2}. Conversely, by \autoref{thm:class-fpf-isos}, any such data yield isometries $\smash{\phi^{(p^l)}}$ of order~$N$ generating fixed-point free groups, and their orthogonal sum is an isometry of $(\Gamma, q)$ generating a fixed-point free group of order~$N$, i.e.\ a fixed-point free action of $\modZ{N}$ on $(\Gamma,q)$.

	\smallskip

	We construct an example of order $N = 12$.
	\begin{ex}
		We choose the monic non-self-reciprocal polynomials
		\begin{enumerate}
			\item $Q_1 \defeq T + 2 \in \modZ{13}[T]$ with multiplicity $r_1 = 2$,
			\item $Q_2 \defeq T^2 + 30 \in (\modZ{7^2})[T]$, the canonical Hensel lift of $T^2 + 2 \in \finitefield{7}[T]$, with multiplicity $r_2 = 1$.
		\end{enumerate}

		Then, there is a fixed-point free action of $\modZ{12}$ on the discriminant form with symmetric bilinear form on $(\modZ{13})^{4}$ and $(\modZ{7^2})^4$,
		\begin{align*}
			\begin{bmatrix}
				0 & \frac{7}{13} & 0 & 0 \\
				\frac{7}{13} & 0 & 0 & 0 \\
				0 & 0 & 0 & \frac{7}{13} \\
				0 & 0 & \frac{7}{13} & 0 \\
			\end{bmatrix}\!, \quad
			\begin{bmatrix}
				0 & 0 & \frac{25}{7^2} & 0 \\
				0 & 0 & 0 & \frac{25}{7^2} \\
				\frac{25}{7^2} & 0 & 0 & 0 \\
				0 & \frac{25}{7^2} & 0 & 0 \\
			\end{bmatrix}\!,
		\end{align*}
		respectively. Following \autoref{nota:discr-form-symbol}, we denote the invariants of this discriminant form by $(13)^{+4} (7^2)^{+4}$. The action is generated by the automorphism acting as
		\begin{align*}
			\begin{bmatrix}
				11 & 0 & 0 & 0\\
				0 & 6 & 0 & 0 \\
				0 & 0 & 11 & 0 \\
				0 & 0 & 0 & 6 \\
			\end{bmatrix}\!, \quad
			\begin{bmatrix}
				0 & 19 & 0 & 0 \\
				1 & 0 & 0 & 0 \\
				0 & 0 & 0 & 31 \\
				0 & 0 & 1 & 0\\
			\end{bmatrix}
		\end{align*}
		by multiplication from the left.
	\end{ex}


	\appendix

	\section{The Witt Class of an Involution Trace} \label{sec:appendix}

	In \autoref{subsec:conj-cls}, we determined conjugacy classes in the orthogonal group of a regular quadratic module $(M,q)$ over an Artinian local ring by decomposing the minimal polynomial and then obtaining an orthogonal decomposition of $(M,q)$.

	In this appendix, we prove that the quadratic module corresponding to an involution trace has even rank and nontrivial Arf invariant, thus describing its Witt class; see \autoref{lem:witt-cls-invol} and \autoref{lem:witt-cls-nontrivial}. We will use the fundamental theorem of Galois theory for Henselian local rings, \autoref{prop:Galois-Henselian-rings}, to study the étale algebra associated with an involution trace. In the end, our proof relies on explicitly calculating the discriminant for the extensions of the corresponding residue fields; see \autoref{cor:discr}.


	\subsection{Discriminants in Towers of Field Extensions}

	We will see that the Witt class of the quadratic module defined by an involution trace depends only on the discriminant of a field extension. Hence, we recall the following two lemmas.

	\begin{lem}[{\cite[Proposition~\RN{4}.6.7 and Proposition~\RN{4}.6.11]{MR643362}}] \label{lem:disc}
		Let $A$ be a ring and $P \in A[T]$ a monic polynomial of degree $d$. We define $E \defeq \inquo{A[T]}{\ideal{P}}$ and denote the image of $T$ in $E$ by $t$. Then the discriminant of $P$ is given by
		\begin{align*}
			\discr(P) = (-1)^{d(d-1)/2} \norm{\fieldextension{E}{A}}(\deriv{P}(t)).
		\end{align*}
		If $P = (T-\alpha_1) \dots (T-\alpha_d)$, then
		\begin{align*}
			\discr(P) = \prod_{i < j} \big(\alpha_i - \alpha_j \big)^2.
		\end{align*}
	\end{lem}

	\smallskip

	\begin{lem}
		Let $K$ be a finite field and $L$ a finite extension of $K$. Then for all $x \in L$, $\norm{\fieldextension{L}{K}}(x)$ is a perfect square if and only if $x$ is a perfect square.
	\end{lem}
	\begin{proof}
		The norm $\norm{\fieldextension{L}{K}} \colon \units{L} \to \units{K}$ is a surjective homomorphism of cyclic groups; see e.g.\ \cite[Theorem~2.28]{MR1429394}. Hence, the preimage of the subgroup of squares of~$\units{K}$ contains the subgroup of squares of~$\units{L}$, and both have index $\gcd(2, \abs{K} - 1) = \gcd(2, \abs{L} - 1)$, so they are equal. For $x = 0$, there is nothing to show.
	\end{proof}

	\medskip

	We use the fundamental theorem of Galois theory and calculate a discriminant using the following lemma.
	\begin{lem} \label{lem:tower-discr}
		Let $K \subseteq L \subseteq M$ be a tower of finite fields with $m = \extord{L}{K}$ and $n = \extord{M}{L}$. Let $(x_1, \dots, x_m)$ be a $K$-basis for $L$ and $(y_1, \dots, y_n)$ an $L$-basis for~$M$. Then $(x_1y_1,x_1y_2,\dots, x_m y_n)$ is a $K$-basis for $M$ and
		\begin{align*}
			\discr_{\fieldextension{M}{K}}&(x_1y_1,\dots,x_m y_n) \\
			& = \norm{\fieldextension{L}{K}}(\discr_{\fieldextension{M}{L}}(y_1,\dots,y_n)) \discr_{\fieldextension{L}{K}}(x_1,\dots,x_m)^n.
		\end{align*}
	\end{lem}
	\begin{proof}
		Let $\smash{t_{j,l} \defeq \tr_{\fieldextension{M}{L}}(y_j y_l)}$, so that $\smash{\det(t_{j,l})_{j,l} = \discr_{\fieldextension{M}{L}}(y_1,\dots,y_n)}$, and let $u$ be the $L$-linear endomorphism of~$L^n$ with matrix $(t_{j,l})_{j,l}$ with respect to the standard basis $(e_1, \dots, e_n)$. By the transitivity of traces \cite[Proposition~\RN{3}.9.6, Corollaire]{MR274237}, for all $i, k$ and $j, l$,
		\begin{align*}
			\tr_{\fieldextension{M}{K}}(x_i y_j x_k y_l) = \tr_{\fieldextension{L}{K}}(x_i x_k t_{j,l}) = \beta(x_i e_j, u(x_k e_l)),
		\end{align*}
		where $\beta$ is the $K$-bilinear form on~$L^n$ given by $\beta(v, w) \defeq \smash{\sum_{j} \tr_{\fieldextension{L}{K}}(v_j w_j)}$. Hence, with respect to the $K$-basis $(x_k e_l)_{k,l}$ of~$L^n$, the matrix $\smash{(\tr_{\fieldextension{M}{K}}(x_i y_j x_k y_l))_{(i,j),(k,l)}}$ is the product of the Gram matrix of~$\beta$, which consists of $n$ diagonal blocks $\smash{(\tr_{\fieldextension{L}{K}}(x_i x_k))_{i,k}}$, and the matrix of~$u$ considered as a $K$-linear endomorphism, whose determinant is $\smash{\norm{\fieldextension{L}{K}}(\det(t_{j,l})_{j,l})}$ by \cite[Proposition~\RN{3}.9.6]{MR274237}. Taking determinants yields the assertion.
	\end{proof}

	In particular, we obtain:
	\begin{cor} \label{cor:discr}
		Let $K$ be a finite field of characteristic not $2$. Let $N$ be a cyclic Galois extension of $K$ of even order~$d$. Let $\sigma$ be the field automorphism of $N$ over~$K$ of order~$2$ and $N^\sigma$ the corresponding fixed field. Then we get in $\inquo{\units{K}}{(\units{K})^2}$ the equality
		\begin{align*}
			\discr_{\fieldextension{N}{K}} \equiv \norm{\fieldextension{N^{\sigma}}{K}}(\discr_{\fieldextension{N}{N^{\sigma}}}) \modulo{(\units{K})^2}.
		\end{align*}
		If $t \in N \setminus N^{\sigma}$, we get
		\begin{align*}
			\discr_{\fieldextension{N}{K}} \equiv \norm{\fieldextension{N^{\sigma}}{K}}((t-\sigma t)^2) \modulo{(\units{K})^2}.
		\end{align*}
	\end{cor}
	\begin{proof}
		Let $P \in N^{\sigma}[T]$ be the minimal polynomial of $t$ over~$N^{\sigma}$, i.e.\
		\begin{align*}
			P = (T-t) (T-\sigma t).
		\end{align*}
		Then, the assertion follows from \autoref{lem:disc} and \autoref{lem:tower-discr}.
	\end{proof}

	\medskip

	Recall that we can identify the Witt ring $\witt{A}$ of a Henselian local ring $A$ with the Witt ring $\witt{k}$ of its residue field~$k$; see \autoref{lem:witt-iso-Henselian}. We freely use the notation from \autoref{subsec:quadr-Henselian-rings}.
	\begin{lem} \label{lem:witt-cls-invol}
		Let $A$ be a Henselian local ring with finite residue field~$k$ of odd characteristic and $P \in A[T]$ a self-reciprocal polynomial of degree at least $2$ such that $\overline{P}$ is irreducible over $k$. Let $E = \inquo{A[T]}{\ideal{P}}$ be the étale algebra over $A$ with involution $\sigma \colon E \to E,\ t \mapsto t^{-1}$ and involution trace $\tr_{\fieldextension{E}{A}} \colon E \to A$.

		Then the Witt class of the quadratic module defined by the involution trace $\lambda_{\fieldextension{E}{A}}$ for the involution $\sigma$ is given by $\Lambda$, i.e.\ it has even rank and nontrivial Arf invariant:
		\begin{align*}
			d(\wittcls{\lambda_{\fieldextension{E}{A}}}) = 0, \quad \legendre{\Arf(\wittcls{\lambda_{\fieldextension{E}{A}}})}{k} = -1.
		\end{align*}
	\end{lem}
	\begin{proof}
		Since $\overline{P}$ is irreducible over the finite field~$k$, it has $d \defeq \deg{P}$ simple roots in~$\overline{E}$, which lift to roots of~$P$ in the Henselian local ring~$E$ (see the proof of \autoref{prop:Galois-Henselian-rings}) by \autoref{lem:lifting-of-zeros}. Hence, the Frobenius automorphism of $\fieldextension{\overline{E}}{k}$ lifts to an automorphism $\alpha$ of the $A$-algebra~$E$ with $P = \prod_{r = 0}^{d-1} (T - \alpha^r(t))$ in $E[T]$, so that $\smash{\tr_{\fieldextension{E}{A}} = \sum_{r = 0}^{d-1} \alpha^r}$ and $\smash{\norm{\fieldextension{E}{A}} = \prod_{r = 0}^{d-1} \alpha^r}$.
		We calculate the determinant as
		\begin{align*}
			\det{\tr_{\fieldextension{E}{A}}(t^i \sigma(t^j))} & = \det{\tr_{\fieldextension{E}{A}}(t^{i-j})} \\
			& = \prod_{r < s} \big(\alpha^s(t) - \alpha^r(t)\big) \big(\alpha^s(t^{-1}) - \alpha^r(t^{-1})\big).
		\end{align*}
		Since the matrix $\smash{\big(\tr_{\fieldextension{E}{A}}(t^{i-j})\big)_{i,j} = \big(\sum_{r = 0}^{d-1} \alpha^r(t^{i-j})\big)_{i,j}}$ can be written as the product of the matrices
		\begin{align*}
			\begin{bmatrix}
				1 & 1 & \dots & 1 \\
				t & \alpha(t) & \dots & \alpha^{d-1}(t) \\
				t^2 & \alpha(t)^2 & \dots & \alpha^{d-1}(t)^2 \\
				& & \vdots & \\
				t^{d-1} & \alpha(t)^{d-1} & \dots & \alpha^{d-1}(t)^{d-1}
			\end{bmatrix}\!,
			\quad
			\begin{bmatrix}
				1 & t^{-1} & \dots & (t^{-1})^{d-1} \\
				1 & \alpha(t^{-1}) & \dots &\alpha(t^{-1})^{d-1} \\
				1 & \alpha^2(t^{-1}) & \dots & \alpha^2(t^{-1})^{d-1} \\
				& & \vdots & \\
				1 & \alpha^{d-1}(t^{-1}) & \dots & \alpha^{d-1}(t^{-1})^{d-1}
			\end{bmatrix}\!,
		\end{align*}
		we obtain that
		\begin{align*}
			\prod_{r < s} \big(\alpha^s(t^{-1}) - \alpha^r(t^{-1})\big) = \prod_{r < s} \big( \alpha^r(t) - \alpha^s(t) \big) \big(\alpha^r(t)\alpha^s(t)\big)^{-1}.
		\end{align*}
		Now, by using
		\begin{align*}
			\prod_{r < s} \big(\alpha^r(t)\alpha^s(t)\big) = \prod_{i = 0}^{d-1} \alpha^i(t^{d-1}) = \norm{\fieldextension{E}{A}}(t)^{d-1},
		\end{align*}
		we can conclude that
		\begin{align*}
			& \det{\tr_{\fieldextension{E}{A}}(t^i \sigma(t^j))} = \det{\tr_{\fieldextension{E}{A}}(t^{i-j})} \\
			& \quad = (-1)^{d(d-1)/2} \norm{\fieldextension{E}{A}}(t^{-(d-1)}) \prod_{r < s} \big(\alpha^s(t) - \alpha^r(t)\big)^2.
		\end{align*}
		Since $E = \inquo{A[T]}{\ideal{P}}$ and $\norm{\fieldextension{E}{A}}(t) = \det{C(P)}$, we obtain by \autoref{rem:self-reciprocal-constant}
		\begin{align*}
			\norm{\fieldextension{E}{A}}(t) = P(0) = 1.
		\end{align*}
		Thus, we get
		\begin{align*}
			(-1)^{d(d-1)/2} \det{\tr_{\fieldextension{E}{A}}(t^i \sigma(t^j))} = \discr_{\fieldextension{E}{A}}(1,t,\dots,t^{d-1}).
		\end{align*}
		By \autoref{lem:square-red}, the square class of the discriminant is determined by the square class on the corresponding residue fields. By \autoref{cor:discr}, the square class of the discriminant can be determined as
		\begin{align*}
			\discr_{\fieldextension{\overline{E}}{k}}(1,\overline{t},\dots,\overline{t}^{d-1}) \equiv \norm{\fieldextension{\overline{E}^{\overline{\sigma}}}{k}}((\overline{t} - \overline{\sigma} \overline{t})^2) \modulo{(\units{k})^2}.
		\end{align*}
		Since $\overline{t}-\overline{\sigma} \overline{t}$ is not in the fixed field $\smash{\overline{E}^{\overline{\sigma}}}$ of $\overline{\sigma}$, the assertion follows.
	\end{proof}

	Using the fact that Hermitian modules have orthonormal bases if the involution is nontrivial (see \autoref{thm:equivalence-Hermitian-forms}), we can generalise \autoref{lem:witt-cls-invol} and obtain:
	\begin{lem} \label{lem:witt-cls-nontrivial}
		Let $A$ be a Henselian local ring with finite residue field~$k$ of odd characteristic, $(M,q)$ a regular quadratic module over $A$, and $\phi \colon (M,q) \to (M,q)$ an isometry with self-reciprocal minimal polynomial $P$ of degree at least~$2$ such that $\overline{P}$ is irreducible over $k$. We define $E \defeq \inquo{A[T]}{\ideal{P}}$. Then
		\begin{align*}
			\wittcls{(M,q)} = (\rank{E}{M}) \Lambda.
		\end{align*}
	\end{lem}
	\begin{proof}
		By \autoref{lem:ex-hermitian-from}, there is a regular $1$-Hermitian form $(M_{\phi}, s_{\phi})$ over $E$ such that
		\begin{align*}
			\lambda_q(x,y) = \tr_{\fieldextension{E}{A}}(s_{\phi}(x,y))
		\end{align*}
		for all $x,y \in M$. We define $n \defeq \rank{E}{M}$.

		By (the proof of) \autoref{thm:equivalence-Hermitian-forms}, $(M_{\phi},s_{\phi})$ has an orthonormal basis. So, let $(x_k)_{k = 1}^{n}$ be an orthonormal basis of the $E$-module $M_{\phi}$. Then, $(t^i x_k)$ forms an $A$-basis of~$M$, and
		\begin{align*}
			\lambda_q(t^i x_k, t^j x_l) & = \tr_{\fieldextension{E}{A}} (t^i \sigma(t^j) s_{\phi}(x_k, x_l)) \\
			& = \tr_{\fieldextension{E}{A}}(t^{i-j}) \delta_{k,l}.
		\end{align*}
		Thus, the assertion follows from the previous lemma.
	\end{proof}


	\bibliographystyle{alpha}
	\bibliography{extracted.bib}

@Preamble{{\newcommand{\SortNoop}[1]{}}}

@article{MR0056295,
  author =        {Eilenberg, Samuel and {Mac Lane}, Saunders},
  journal =       {Ann. of Math.},
  pages =         {55--106},
  title =         {On the groups {$H(\Pi,n)$}. {I}},
  volume =        {58},
  year =          {1953},
}

@article{MR1250465,
  author =        {Joyal, Andr\'{e} and Street, Ross},
  journal =       {Adv. Math.},
  number =        {1},
  pages =         {20--78},
  title =         {Braided tensor categories},
  volume =        {102},
  year =          {1993},
}

@book{MR1233387,
  author =        {Dong, Chongying and Lepowsky, James~I.},
  publisher =     {Birkh{\"a}user},
  series =        {Progr. Math.},
  title =         {Generalized vertex algebras and relative vertex operators},
  volume =        {112},
  year =          {1993},
}

@unpublished{paper2,
  author =        {Kn\"otzele, Hannes and M\"oller, Sven},
  title =         {Constructing Modular $\mathbb{Z}/{N}$-Crossed {T}ambara-{Y}amagami Categories},
  year =          {2027},
  note =          {In preparation},
}

@article{MR0525944,
  author =        {Nikulin, Vyacheslav~V.},
  journal =       {Izv. Akad. Nauk SSSR Ser. Mat.},
  number =        {1},
  pages =         {111--177, 238},
  title =         {Integer symmetric bilinear forms and some of their geometric applications},
  volume =        {43},
  year =          {1979},
}

@unpublished{scheithauer,
  author =        {Scheithauer, Nils~R.},
  note =          {Unpublished manuscript},
  title =         {Discriminant forms and their automorphisms},
  year =          {2013},
}

@article{MR150210,
  author =        {Wall, Gordon~E.},
  journal =       {J. Austral. Math. Soc.},
  pages =         {1--62},
  title =         {On the conjugacy classes in the unitary, symplectic and orthogonal groups},
  volume =        {3},
  year =          {1963},
}

@article{MR591256,
  author =        {Shinoda, Ken-ichi},
  journal =       {J. Algebra},
  number =        {1},
  pages =         {251--280},
  title =         {The characters of {W}eil representations associated to finite fields},
  volume =        {66},
  year =          {1980},
}

@article{MR249519,
  author =        {Milnor, John},
  journal =       {Invent. Math.},
  pages =         {83--97},
  title =         {On isometries of inner product spaces},
  volume =        {8},
  year =          {1969},
}

@incollection{MR268192,
  author =        {Springer, Tonny~A. and Steinberg, Robert},
  booktitle =     {Seminar on {A}lgebraic {G}roups and {R}elated {F}inite {G}roups ({T}he {I}nstitute for {A}dvanced {S}tudy, {P}rinceton, {N}.{J}., 1968/69)},
  pages =         {167--266},
  publisher =     {Springer},
  series =        {Lecture Notes in Math.},
  title =         {Conjugacy classes},
  volume =        {131},
  year =          {1970},
}

@article{MR4968024,
  author =        {Gannon, Terry and Riesen, Andrew},
  journal =       {Adv. Math.},
  pages =         {110546},
  title =         {Orbifolds of pointed vertex operator algebras {I}},
  volume =        {482},
  year =          {2025},
}

@article{GLM26,
   author =       {Galindo, C\'esar and Lentner, Simon and M\"oller, Sven},
   title =        {Modular {$\mathbb{Z}_2$}-crossed {T}ambara-{Y}amagami-like categories for even groups},
   journal =      {Adv. Math.},
   volume =       {493},
   year =         {2026},
   pages =        {110902},
}

@unpublished{GLM24a,
  author =        {Galindo, C\'{e}sar and Lentner, Simon and M{\"o}ller, Sven},
  note =          {\href{http://arxiv.org/abs/2409.16357v2}{arXiv:2409.16357v2 [math.QA]}},
  title =         {Computing {$G$}-crossed extensions and orbifolds of vertex operator algebras},
  year =          {2024},
}

@article{MR252382,
  author =        {Knebusch, Manfred},
  journal =       {Math. Z.},
  pages =         {255--268},
  title =         {Isometrien \"uber semilokalen {R}ingen},
  volume =        {108},
  year =          {1969},
}

@article{MR1773561,
  author =        {Borcherds, Richard~E.},
  journal =       {Duke Math. J.},
  number =        {2},
  pages =         {319--366},
  title =         {Reflection groups of {L}orentzian lattices},
  volume =        {104},
  year =          {2000},
}

@article{MR3101818,
  author =        {Str\"{o}mberg, Fredrik},
  journal =       {Math. Z.},
  number =        {1--2},
  pages =         {509--527},
  title =         {Weil representations associated with finite quadratic modules},
  volume =        {275},
  year =          {2013},
}

@article{MR2609644,
  author =        {Drinfeld, Vladimir~G. and Gelaki, Shlomo and Nikshych, Dmitri~A. and Ostrik, Victor},
  journal =       {Selecta Math. (N.S.)},
  number =        {1},
  pages =         {1--119},
  title =         {On braided fusion categories. {I}},
  volume =        {16},
  year =          {2010},
}

@book{MR277519,
  author =        {Raynaud, Michel},
  publisher =     {Springer},
  series =        {Lecture Notes in Math.},
  title =         {Anneaux locaux hens\'eliens},
  volume =        {169},
  year =          {1970},
}

@book{MR559531,
  author =        {Milne, James~S.},
  publisher =     {Princeton University Press},
  series =        {Princeton Math. Ser.},
  title =         {\'Etale cohomology},
  volume =        {33},
  year =          {1980},
}

@article{MR40287,
  author =        {Azumaya, Goro},
  journal =       {Nagoya Math. J.},
  pages =         {119--150},
  title =         {On maximally central algebras},
  volume =        {2},
  year =          {1951},
}

@article{MR1578516,
  author =        {Sch\"onemann, Theodor},
  journal =       {J. Reine Angew. Math.},
  pages =         {93--105},
  title =         {Von denjenigen {M}oduln, welche {P}otenzen von {P}rimzahlen sind. ({F}ortsetzung)},
  volume =        {32},
  year =          {1846},
}

@book{MR3618889,
  author =        {Ford, Timothy~J.},
  publisher =     {Amer. Math. Soc.},
  series =        {Grad. Stud. Math.},
  title =         {Separable algebras},
  volume =        {183},
  year =          {2017},
}

@article{MR364218,
  author =        {Ganske, Gary~L. and McDonald, Bernard~R.},
  journal =       {Rocky Mountain J. Math.},
  number =        {4},
  pages =         {521--540},
  title =         {Finite local rings},
  volume =        {3},
  year =          {1973},
}

@article{MR1180736,
  author =        {Whelan, E. A.},
  journal =       {Rocky Mountain J. Math.},
  number =        {2},
  pages =         {757--759},
  title =         {A note on finite local rings},
  volume =        {22},
  year =          {1992},
}

@article{MR1778166,
  author =        {Corbas, Basil and Williams, Graham~D.},
  journal =       {J. Algebra},
  number =        {2},
  pages =         {691--704},
  title =         {Rings of order {$p^5$}. {II}. {L}ocal rings},
  volume =        {231},
  year =          {2000},
}

@book{MR1096299,
  author =        {Knus, Max-Albert},
  publisher =     {Springer},
  series =        {Grundlehren Math. Wiss.},
  title =         {Quadratic and {H}ermitian forms over rings},
  volume =        {294},
  year =          {1991},
}

@phdthesis{MR43785,
  author =        {Springer, Tonny~A.},
  school =        {University of Leiden},
  title =         {Over symplectische transformaties},
  year =          {1951},
}

@article{MR1957,
  author =        {Jacobson, Nathan},
  journal =       {Bull. Amer. Math. Soc.},
  pages =         {264--268},
  title =         {A note on hermitian forms},
  volume =        {46},
  year =          {1940},
}

@book{MR491773,
  author =        {Baeza, Ricardo},
  publisher =     {Springer},
  series =        {Lecture Notes in Math.},
  title =         {Quadratic forms over semilocal rings},
  volume =        {655},
  year =          {1978},
}

@misc{stacks-project,
  author =        {The {Stacks project authors}},
  howpublished =  {\url{https://stacks.math.columbia.edu}},
  title =         {Stacks project},
  year =          {2026},
}

@book{MR396410,
  author =        {Lam, Tsit-Yuen},
  publisher =     {W. A. Benjamin},
  series =        {Math. Lecture Note Ser.},
  title =         {The algebraic theory of quadratic forms},
  year =          {1973},
}

@book{MR0152507,
  author =        {O'Meara, O.~Timothy},
  publisher =     {Springer},
  series =        {Grundlehren Math. Wiss.},
  title =         {Introduction to quadratic forms},
  volume =        {117},
  year =          {1963},
}

@article{MR2496771,
  author =        {Scheithauer, Nils~R.},
  journal =       {Int. Math. Res. Not. IMRN},
  number =        {8},
  pages =         {1488--1545},
  title =         {The {W}eil representation of {$\operatorname{SL}_2(\mathbb{Z})$} and some applications},
  volume =        {2009},
  year =          {2009},
}

@book{MR1189139,
  author =        {Taylor, Donald~E.},
  publisher =     {Heldermann},
  series =        {Sigma Ser. Pure Math.},
  title =         {The geometry of the classical groups},
  volume =        {9},
  year =          {1992},
}

@book{MR770063,
  author =        {Scharlau, Winfried},
  publisher =     {Springer},
  series =        {Grundlehren Math. Wiss.},
  title =         {Quadratic and {H}ermitian forms},
  volume =        {270},
  year =          {1985},
}

@book{MR279159,
  author =        {Bass, Hyman},
  publisher =     {Tata Institute of Fundamental Research},
  series =        {Tata Inst. Fundam. Res. Lect. Math.},
  title =         {Lectures on topics in algebraic {$K$}-theory},
  volume =        {41},
  year =          {1967},
}

@article{MR0156890,
  author =        {Wall, Charles~T.~C.},
  journal =       {Topology},
  pages =         {281--298},
  title =         {Quadratic forms on finite groups, and related topics},
  volume =        {2},
  year =          {1963},
}

@book{MR4928289,
  author =        {Huppert, Bertram},
  publisher =     {Springer},
  series =        {Grundlehren Math. Wiss.},
  title =         {Finite groups. {I}},
  volume =        {364},
  year =          {2025},
}

@article{MR51838,
  author =        {Gasch\"utz, Wolfgang},
  journal =       {Math. Z.},
  pages =         {376--387},
  title =         {Über den {F}undamentalsatz von {M}aschke zur {D}arstellungstheorie der endlichen {G}ruppen},
  volume =        {56},
  year =          {1952},
}

@phdthesis{mayr-thesis,
  author =        {Mayr, Peter},
  type =          {Diplomarbeit},
  school =        {Universität Linz},
  title =         {Finite fixed point free automorphism groups},
  year =          {1998},
  note =          {\url{http://www.algebra.uni-linz.ac.at/~stein/papers/fpfaut.pdf}},
}

@phdthesis{hannes-thesis,
  author =        {Kn{\"o}tzele, Hannes},
  school = {Universit{\"a}t Hamburg},
  title =         {Extensions of fusion categories via fixed-point-free automorphisms},
  year =          {2026},
  note =          {\url{http://ediss.sub.uni-hamburg.de/handle/ediss/12547}},
}

@article{MR27759,
  author =        {Taunt, Derek~R.},
  journal =       {Proc. Cambridge Philos. Soc.},
  pages =         {24--42},
  title =         {On {$A$}-groups},
  volume =        {45},
  year =          {1949},
}

@article{MR0743731,
  author =        {Miranda, Rick},
  journal =       {Trans. Amer. Math. Soc.},
  number =        {2},
  pages =         {535--542},
  title =         {Nondegenerate symmetric bilinear forms on finite abelian {$2$}-groups},
  volume =        {284},
  year =          {1984},
}

@article{MR0480333,
  author =        {Durfee, Alan~H.},
  journal =       {Adv. Math.},
  number =        {2},
  pages =         {133--164},
  title =         {Bilinear and quadratic forms on torsion modules},
  volume =        {25},
  year =          {1977},
}

@article{MR1990929,
  author =        {M{\"u}ger, Michael},
  journal =       {Proc. London Math. Soc.},
  number =        {2},
  pages =         {291--308},
  title =         {On the structure of modular categories},
  volume =        {87},
  year =          {2003},
}

@article{MR3262438,
  author =        {Ferrand, Daniel and Laksov, Dan},
  journal =       {Enseign. Math.},
  number =        {1--2},
  pages =         {131--175},
  title =         {On the structure of endomorphisms of projective modules},
  volume =        {60},
  year =          {2014},
}

@book{MR643362,
  author =        {Bourbaki, Nicolas},
  title =          {\'{E}l\'{e}ments de math\'{e}matique. Alg\`ebre. Chapitres 4 \`a{} 7},
  publisher =     {Masson, Paris},
  year =          {1981},
}

@book{bourbaki-ac-1-4,
  author =        {Bourbaki, Nicolas},
  title =         {\'{E}l\'{e}ments de math\'{e}matique. {A}lg\`ebre commutative. {C}hapitres 1 \`a{} 4},
  publisher =     {Masson, Paris},
  year =          {1985},
}

@book{MR274237,
  author =        {Bourbaki, Nicolas},
  publisher =     {Hermann, Paris},
  title =         {\'El\'ements de math\'ematique. {A}lg\`ebre. {C}hapitres 1 \`a{} 3},
  year =          {1970},
}

@book{MR280479,
  author =        {DeMeyer, Frank~R. and Ingraham, Edward~C.},
  publisher =     {Springer},
  series =        {Lecture Notes in Math.},
  title =         {Separable algebras over commutative rings},
  volume =        {181},
  year =          {1971},
}

@article{MR220757,
  author =        {Davis, Ronald~W.},
  journal =       {Duke Math. J.},
  pages =         {49--59},
  title =         {Certain matrix equations over rings of integers},
  volume =        {35},
  year =          {1968},
}

@article{MR16094,
  author =        {Cohen, Irvin~S.},
  journal =       {Trans. Amer. Math. Soc.},
  pages =         {54--106},
  title =         {On the structure and ideal theory of complete local rings},
  volume =        {59},
  year =          {1946},
}

@book{MR266911,
  author =        {Matsumura, Hideyuki},
  publisher =     {W. A. Benjamin},
  title =         {Commutative algebra},
  year =          {1970},
}

@article{MR173675,
  author =        {Grothendieck, Alexander},
  journal =       {Inst. Hautes \'Etudes Sci. Publ. Math.},
  volume =        {20},
  pages =         {5--259},
  title =         {Éléments de géométrie algébrique. {IV}. {É}tude locale des schémas et des morphismes de schémas. {I}},
  year =          {1964},
}

@book{MR354618,
  author =        {Serre, Jean-Pierre},
  publisher =     {Hermann, Paris},
  series =        {Publications de l'Institut de Math\'ematique de l'Universit\'e{} de Nancago},
  title =         {Corps locaux},
  volume =        {VIII},
  year =          {1968},
  edition =       {2nd},
}

@article{MR103906,
  author =        {Auslander, Maurice and Buchsbaum, David~A.},
  journal =       {Proc. Nat. Acad. Sci. U.S.A.},
  pages =         {733--734},
  title =         {Unique factorization in regular local rings},
  volume =        {45},
  year =          {1959},
}

@book{MR242802,
  author =        {Atiyah, Michael~F. and Macdonald, Ian~G.},
  publisher =     {Addison-Wesley},
  title =         {Introduction to commutative algebra},
  year =          {1969},
}

@article{MR2099718,
  author =        {Frisch, Sophie},
  journal =       {Comm. Algebra},
  number =        {5},
  pages =         {2015--2017},
  title =         {Integrally closed domains, minimal polynomials, and null ideals of matrices},
  volume =        {32},
  year =          {2004},
}

@book{MR1181420,
  author =        {Adkins, William~A. and Weintraub, Steven~H.},
  publisher =     {Springer},
  series =        {Grad. Texts in Math.},
  title =         {Algebra. {A}n approach via module theory},
  volume =        {136},
  year =          {1992},
}

@book{MR1429394,
  author =        {Lidl, Rudolf and Niederreiter, Harald},
  publisher =     {Cambridge University Press},
  series =        {Encyclopedia Math. Appl.},
  title =         {Finite fields},
  volume =        {20},
  edition =       {2nd},
  year =          {1997},
}

@article{MR309963,
  author =        {Pomfret, James~C.},
  journal =       {Proc. Amer. Math. Soc.},
  pages =         {421--422},
  title =         {Similarity of matrices over finite rings},
  volume =        {37},
  year =          {1973},
}

@article{MR914197,
  author =        {Cohen, Stephen~D.},
  journal =       {Bull. Korean Math. Soc.},
  number =        {2},
  pages =         {83--96},
  title =         {Polynomials over finite fields with large order and level},
  volume =        {24},
  year =          {1987},
}

@book{MR1357169,
  author =        {Robinson, Derek~J.~S.},
  edition =       {2nd},
  publisher =     {Springer},
  series =        {Grad. Texts in Math.},
  title =         {A course in the theory of groups},
  volume =        {80},
  year =          {1996},
}

@book{MR569209,
  author =        {Gorenstein, Daniel},
  edition =       {2nd},
  publisher =     {Chelsea Publishing},
  title =         {Finite groups},
  year =          {1980},
}

@book{MR1269324,
  author =        {Weibel, Charles A.},
  publisher =     {Cambridge University Press},
  series =        {Cambridge Stud. Adv. Math.},
  title =         {An introduction to homological algebra},
  volume =        {38},
  year =          {1994},
}

\end{document}